\documentclass[fullpage]{article}

\usepackage{lineno,hyperref}
\modulolinenumbers[5]

\usepackage[utf8]{inputenc}

\usepackage{graphicx,color}
\usepackage{amssymb}
\usepackage{amsthm,amsmath}
\usepackage{appendix}
\usepackage{caption}
\usepackage{subcaption}
\usepackage{float}
\usepackage[normalem]{ulem}

\usepackage{algorithm}
\usepackage[noend]{algpseudocode}
\usepackage{authblk}
\usepackage{newtxtext,newtxmath}

\usepackage{doi}

\newcommand{\boldA}{\mathbf{A}}
\newcommand{\boldAij}{\mathbf{A_{ij}}}
\newcommand{\boldu}{\mathbf{u}}
\newcommand{\bolduij}{\mathbf{u_{ij}}}
\newcommand{\boldutilde}{\mathbf{\tilde{u}}}
\newcommand{\boldf}{\mathbf{f}}
\newcommand{\boldfij}{\mathbf{f_{ij}}}
\newcommand{\boldZero}{\mathbf{0}}
\newcommand{\boldU}{\mathbf{U}}
\newcommand{\boldOmegaOne}{\mathbf{\Omega_1}}
\newcommand{\boldOmegaTwo}{\mathbf{\Omega_2}}
\newcommand{\boldOmegaij}{\mathbf{\Omega_{ij}}}
\newcommand{\boldOmegaijepsilon}{\mathbf{\Omega_{ij}^\delta}}
\newcommand{\boldGammaOne}{\mathbf{\Gamma_1}}
\newcommand{\boldGammaTwo}{\mathbf{\Gamma_2}}
\newcommand{\boldGammaijB}{\mathbf{\Gamma_{ij}^B}}

\newcommand{\boldGammaijT}{\mathbf{\Gamma_{ij}^T}}
\newcommand{\boldGammaijL}{\mathbf{\Gamma_{ij}^L}}
\newcommand{\boldGammaijl}{\mathbf{\Gamma_{ij}^{\ell}}}

\renewcommand{\u}{\mathbf{u}}

\newcommand{\x}{\mathbf{x}}

\def\cO{{\cal O}}

\newcommand{\pluseq}{\mathrel{+}=}

\title{L-Sweeps: A scalable, parallel preconditioner for the high-frequency Helmholtz equation}

\author[1,2]{Matthias Taus}
\author[3,4,5]{Leonardo Zepeda-N\'u\~nez}
\author[6,7]{Russell J. Hewett}
\author[2]{Laurent Demanet}
\affil[1]{Institute of Analysis and Scientific Computing, Vienna University of Technology}
\affil[2]{Department of Mathematics, Massachusetts Institute of Technology}
\affil[3]{Computational Research Division, Lawrence Berkeley National Laboratory}
\affil[4]{Department of Mathematics, University of California Berkeley}
\affil[5]{Department of Mathematics, University of Wisconsin-Madison}
\affil[6]{Total E\&P Research and Technology, USA}
\affil[7]{Department of Mathematics, Virginia Tech}
\date{September 3, 2019}

\begin{document}



\maketitle

\begin{abstract}

We present the first fast solver for the high-frequency Helmholtz equation that scales optimally in parallel, for a single right-hand side. The L-sweeps approach achieves this scalability by departing from the usual propagation pattern, in which information flows in a $180^{\circ}$ degree cone from interfaces in a layered decomposition. Instead, with L-sweeps, information propagates in $90^{\circ}$ cones induced by a checkerboard domain decomposition (CDD). We extend the notion of accurate transmission conditions to CDDs and introduce a new sweeping strategy to efficiently track the wave fronts as they propagate through the CDD. The new approach decouples the subdomains at each wave front, so that they can be processed in parallel, resulting in better parallel scalability than previously demonstrated in the literature.  The method has an overall $\cO \big((N/p)\log \omega\big)$ empirical run-time for $N=n^d$ total degrees-of-freedom in a $d$-dimensional problem, frequency $\omega$, and $p=\cO(n)$ processors.  We introduce the algorithm and provide a complexity analysis for our parallel implementation of the solver. We corroborate all claims in several two- and three-dimensional numerical examples involving constant, smooth, and discontinuous wave speeds.

\end{abstract}

\section{Introduction}
\label{Sec::Introduction}
The Helmholtz equation, a time-harmonic form of the wave equation, arises in modeling many physical phenomena, including electromagnetic and subsurface wave propagation. Such applications are of interest when solving related inverse problems that require solutions at high frequency to recover fine-grained details, e.g., ultrasound and subsurface recovery.  In subsurface recovery, the propagation medium tends to be extremely complex and the asymptotic approximations ubiquitous in other modalities are not sufficiently accurate, thus the numerical solution of the wave equation is required.  Such computations are the backbone of the full-waveform inversion (FWI) method for subsurface recovery~\cite{VirieuxOperto2009,Pratt:Seismic_waveform_inversion_in_the_frequency_domain;_Part_1_Theory_and_verification_in_a_physical_scale_model}. In the context of time-harmonic wave equations, such as the Helmholtz equation, accurate reconstructions of the subsurface via FWI require solutions at a wide range of frequencies.  The finest recoverable detail is determined by the highest frequency for which the wave equation can be solved. Consequently, the efficient solution of time-harmonic wave equations at high-frequency is extremely important in scientific and industrial applications.

We consider the Helmholtz equation with variable wave speed and constant density on an open domain $\Omega_{\texttt{bulk}}$ with absorbing boundary conditions,
\begin{alignat}{1}
\label{HelmholtzBulk}
\begin{aligned}
-\Delta u-\omega^2mu = f\quad\text{in }\Omega_{\texttt{bulk}},\\
\text{+ A.B.C on }\partial\Omega_{\texttt{bulk}},
\end{aligned}
\end{alignat}
where for $x\in\Omega_{\texttt{bulk}}$, $m(x)=1/\,c^2(x)$ is the  squared slowness for the p-wave speed $c$,
$u$ is the solution wavefield, $\omega$ is the characteristic frequency, and $f$ is the source density.  We consider $\Omega_{\texttt{bulk}}$ to be a square or a cube, but this is not a limiting assumption.  While we consider only problems of this form in the ensuing developments, the method we propose is a framework that can be applied to other time-harmonic formulations that model more complex physics.

In particular, we consider discretizations of~\eqref{HelmholtzBulk} in the high-frequency regime. This means that the coarsest structure in the spatial discretization (mesh) has to scale as $1/\omega$.
In this regime the solution of the discrete system is notoriously difficult~\cite{Gander:why_is_difficult_to_solve_helmholtz_problems_with_classical_iterative_methods}. 
It is well-established that, independent of the discretization, the spectrum of the resulting system matrix  deteriorates as the frequency $\omega$ increases ~\cite{Moiola_Spence:Is_the_Helmholtz_Equation_Really_Sign_Indefinite,Spence:Wavenumber-Explicit_Bounds_in_Time-Harmonic_Acoustic_Scattering}. Therefore, it is not feasible to solve systems arising from high-frequency problems without the use of specialized solution strategies which directly address this issue.

A number of solvers for solving the resulting linear systems are available, including direct methods (e.g.,~\cite{Wang_de_Hoop:Massively_parallel_structured_multifrontal_solver_for_time-harmonic_elastic_waves_in_3D_anisotropic_media}), domain decomposition methods~\cite{Toselli:Domain_Decomposition_Methods_Algorithms_and_Theory}, and preconditioned iterative methods (including, e.g.,  multigrid and shifted-Laplacian methods)~\cite{Sheikh_Lahaye_Vuik:On_the_convergence_of_shifted_Laplace_preconditioner_combined_with_multilevel_deflation,Gander_Graham_Spence:largest_shift_for_complex_shitfted_Laplacian}. We consider a solver or a preconditioner to be {\it sequentially scalable} if, up to logarithmic factors, it can compute a solution in $\cO(N)$ run-time in a sequential computational environment where $N$ is the total number of degrees-of-freedom in the discrete problem. We consider a solver or a preconditioner to be {\it parallel scalable} if, up to logarithmic factors, it can compute a solution in $\cO(N/p)$ run-time in a parallel computational environment with $p=\cO(n)$ processors, for a single right-hand side. Here, $n$ is the number of degrees of freedom in one spatial direction, i.e., $n=\cO(N^{\frac{1}{d}})$ where $d$ is the problem dimension.  In this paper, we present what we believe to be the first \textit{parallel scalable} preconditioner for the high-frequency Helmholtz equation.

Currently, no scalable direct method, sequential or parallel, is available for the high-frequency regime. Standard domain decomposition methods (DDMs) localize the problem to subdomains and transfer information between subdomains.  Domain decomposition methods can be applied as a direct solver (e.g., in the context of optimal Schwarz methods~\cite{Gander_Nataf:Optimized_Schwarz_Methods_without_Overlap_for_the_Helmholtz_Equation,Gander:Optimized_Schwarz_Methods}), but the resulting solver is not scalable.  Domain decomposition methods can yield scalable preconditioners, both in parallel and sequentially.  However, the resulting preconditioned iterative solver is not scalable because it requires $\cO(\omega)$ iterations. This shows an important aspect for the construction of solvers based on preconditioning: a quality preconditioner has to be scalable but also the resulting iterative method has to converge with limited growth in iterations as $\omega$ increases. In this regard, we call a preconditioner {\it effective} if the resulting preconditioned system can be solved in $\cO(\log\omega)$ iterations.  Classical DDMs exhibit sub-optimal behavior for two reasons.  First, artificial and spurious reflections are induced by imprecise information transfer between subdomains.  Second, long-range wave-material interactions are not tracked consistently.

Sweeping preconditioners have been introduced~\cite{EngquistYing:Sweeping_PML,GeuzaineVion:double_sweep,Chen_Xiang:a_source_transfer_ddm_for_helmholtz_equations_in_unbounded_domain,ZepedaDemanet:the_method_of_polarized_traces,GanderZhang2019} to alleviate these drawbacks, while preserving the advantages of classical DDMs.
Sweeping preconditioners make use of layered domain decompositions, accurate transmission conditions, and a layer-by-layer sweeping strategy.
In particular, the layered domain decomposition provides scalability by controlling growth in computational cost and memory footprint.  The accurate transmission conditions allow information to flow between the subdomains without numerical artifacts, e.g., artificial reflections.  Finally, the sweeping strategy consistently tracks and propagates long-range wave-material interactions. The resulting approach can be interpreted as an approximate block-LU factorization, where the blocks correspond to the local problem in each subdomain. In particular, using sparse direct solvers on blocks that arise from sufficiently thin layers yields a method with quasi-linear (i.e., linear with poly-logarithmic factors) sequential complexity for a single right-hand side~\cite{EngquistYing:Sweeping_PML,EngquistYing:Sweeping_H,ZepedaDemanet:Nested_domain_decomposition_with_polarized_traces_for_the_2D_Helmholtz_equation,Liu_Ying:Recursive_sweeping_preconditioner_for_the_3d_helmholtz_equation,Stolk:An_improved_sweeping_domain_decomposition_preconditioner_for_the_Helmholtz_equation,Chen_Xiang:a_source_transfer_ddm_for_helmholtz_equations_in_unbounded_domain,GeuzaineVion:double_sweep,Modave2018}. In the presence of many right-hand sides, the layered domain decomposition allows for optimal parallelization~\cite{Zepeda2019}.

Efforts to improve the performance of sweeping preconditioners have focused on obtaining better and more accurate factorizations~\cite{GeuzaineVion:double_sweep}, reducing the over-all cost by restricting the unknowns to the interfaces~\cite{ZepedaDemanet:the_method_of_polarized_traces,ZepedaDemanet:Nested_domain_decomposition_with_polarized_traces_for_the_2D_Helmholtz_equation}, and accelerating the computation on local subproblems (blocks in the approximate factorization) with compression or parallelization~\cite{Poulson_Engquist:a_parallel_sweeping_preconditioner_for_heteregeneous_3d_helmholtz}.  Several approaches aim specifically to sparsify those blocks, thus decreasing the sequential costs \cite{Liu_Ying:Recursive_sweeping_preconditioner_for_the_3d_helmholtz_equation,ZepedaDemanet:Nested_domain_decomposition_with_polarized_traces_for_the_2D_Helmholtz_equation}. However, while leveraging parallelism to accelerate the solve for local blocks decreases the run-times, it does not result in a parallel scalable solver for a single right-hand side. 

The main bottleneck of current sweeping preconditioners is a lack of parallel scalability for a single right-hand side. The difficulty arises because the Helmholtz problem is inherently sequential, independent of domain decomposition strategy. This sequential nature is dictated by the hyperbolic nature of wave equations and is manifested in the need to accurately resolve the long-range interactions. The key to accurately resolving these interactions is a consistent information transfer between subdomains, which has precluded more general domain decompositions (i.e., beyond layered subdomains) and consequently inhibited parallelization.

In this work, we address these issues by departing from standard layered domain decompositions and introducing a checkerboard domain decomposition (CDD) along with a new sweeping strategy that is only viable on such domain decompositions.  Our novel sweeping preconditioner consistently and efficiently tracks the wavefield across subdomains.  The preconditioner can be interpreted as an approximate LU factorization, where parallelism arises because the diagonal blocks themselves have a block-diagonal structure. For a single right-hand side, the resulting algorithm has a $\cO(N/p)$ run-time, up to logarithmic factors, where, independent of geometric dimension, $p={\cal O}(n)$ is the number of processors. The new algorithm therefore provides the first parallel scalable solver for the Helmholtz equation at high-frequency for a single right-hand side.

Our approach can be applied to two- and three-dimensional problems.  In two dimensions, the algorithm exhibits good weak parallel scalability.  That is, as the frequency increases, and consequently the problem is refined, we refine the CDD so that the local problems in each subdomain have constant size. In three dimensions, we apply the same strategy as for the two-dimensional case in two of the spatial dimensions and extend the subdomains along the third spatial dimension resulting in beam-shaped quasi-one-dimensional local problems. In either case, we employ off-the-shelf direct solvers to solve the local problems and obtain parallel complexities of $\cO(n\log\omega)$ and $\cO(n^2\log\omega)$ in two and three dimensions, respectively. The 3D complexity can be reduced to ${\cal O}(n\log\omega)$ by employing existing parallel direct solvers~\cite{Rouet_Li_Ghysels:A_distributed-memory_package_for_dense_Hierarchically_Semi-Separable_matrix_computations_using_randomization,li_demmel03:SuperLU_DIST} for the quasi-one-dimensional problems in each subdomain, but we do not exploit these tools here.

The new algorithm therefore results in an $\cO(n\log\omega)$ complexity for a single right-hand side, regardless of the geometric dimension. This new algorithm requires sweeps across the CDD in both cardinal and diagonal directions, and we can exploit parallelism in all directions orthogonal to the current sweep direction.  While this new strategy increases parallelization, the sweeps are still inherently sequential and cannot be parallelized.  Thus, we do not anticipate further reduction in complexity, for a single right-hand, due exclusively to modification of the domain decomposition. In the presence of ${\cal O}(n\log\omega)$ right-hand-sides, it has been shown, however, that the sequential nature of the sweeps can be mitigated by pipelining multiple right-hand sides~\cite{Zepeda2019}. Using pipelining, solutions for all $\cO(n)$ right-hand-sides can  be computed with ${\cal O}(n\log\omega)$ parallel complexity, i.e., the average parallel complexity per right-hand-side is $\cO(\log\omega)$.

\subsection{Related Work} \label{section:related_work}
Our method is inspired by the method of polarized traces, which is well-established in the literature and, in contrast to other proposed preconditioners, has been proven applicable for problems with high-order discretizations~\cite{ZepedaZhao:Fast_Lippmann_Schwinger_solver,Taus_Demanet_Zepeda:HDG_Helmholtz} and highly heterogeneous (and even discontinuous) wave speed distributions~\cite{ZepedaDemanet:Nested_domain_decomposition_with_polarized_traces_for_the_2D_Helmholtz_equation}.

Standard linear algebra techniques such as nested dissection~\cite{GeorgeNested_dissection} and multi-frontal solvers~\cite{Duff_Reid:The_Multifrontal_Solution_of_Indefinite_Sparse_Symmetric_Linear,Xia:multifrontal}, coupled with $\mathcal{H}$-matrices~\cite{Hackbusch:Hierarchical_matrices}, have been applied to the Helmholtz problem~\cite{Gillman_Barnett_Martinsson:A_spectrally_accurate_solution_technique_for_frequency_domain_scattering_problems_with_variable_media,Wang:H_multifrontal,Wang_de_Hoop:Massively_parallel_structured_multifrontal_solver_for_time-harmonic_elastic_waves_in_3D_anisotropic_media,Amestoy:Fast_3D_frequency_domain_full_waveform_inversion_with_a_parallel_block_low-rank_multifrontal_direct_solver_Application_to_OBC_data_from_the_North_Sea}.  While these methods take advantage of compressed linear algebra to gain more efficiency (e.g.,~\cite{Bebendorf:2008}), in the high-frequency regime they still suffer from the same sub-optimal asymptotic complexity as standard multi-frontal methods (e.g.,~\cite{Demmel_Li:superlu,Amestoy_Duff:MUMPS,Davis:UMFPACK}).

Multigrid methods (e.g.,~\cite{LairdGiles:Preconditioned_iterative_solution_of_the_2D_Helmholtz_equation,Brandt_Livshits:multi_ray_multigrid_standing_wave_equations,Erlangga:shifted_laplacian,Sheikh_Lahaye_Vuik:On_the_convergence_of_shifted_Laplace_preconditioner_combined_with_multilevel_deflation,Mulder:A_multigrid_solver_for_3D_electromagnetic_diffusion,Erlangga:shifted_laplacian,AruliahAscher:Multigrid_Preconditioning_for_Krylov_Methods_for_Time-Harmonic_Maxwells_Equations_in_Three_Dimensions}), once thought to be inefficient for the Helmholtz problem, have been successfully applied~\cite{Calandra_Grattonn:an_improved_two_grid_preconditioner_for_the_solution_of_3d_Helmholtz,Hu_Zhang:Substructuring_Preconditioners_for_the_Systems_Arising_from_Plane_Wave_Discretization_of_Helmholtz_Equations,Stolk}.  Approaches stemming from the complex-shifted Laplacian~\cite{Erlangga:shifted_laplacian} can be advantageous if properly tuned. However, in general, they either require an expensive solver for the shifted problem or require a large number of iterations to reach convergence, depending on the scaling between the complex shift and the frequency~\cite{Gander_Graham_Spence:largest_shift_for_complex_shitfted_Laplacian}.  Although these algorithms do not result in a lower computational complexity, they are highly parallelizable, resulting in low run-times.

Within the geophysical community, the analytic incomplete LU (AILU) method was explored in~\cite{Plessix_Mulder:Separation_of_variable_preconditioner_for_iterativa_Helmholtz_solver,Plessix:A_Helmholtz_iterative_solver_for_3D_seismic_imaging_problems}.
A variant of Kaczmarz preconditioners~\cite{Gordon:A_robust_and_efficient_parallel_solver_for_linear_systems} has been studied and applied to time-harmonic wave equations by~\cite{Brossier:2D_and_3D_frequency-domain_elastic_wave_modeling_in_complex_media_with_a_parallel_iterative_solver}.
Another class of methods, called hybrid direct-iterative methods, have been explored by~\cite{Sourbier_2011_GEOP}.
Although these solvers have, in general, relatively low memory consumption they tend to require many iterations to converge, thus hindering practical run-times.

Domain decomposition methods for solving partial differential equations (PDEs) have a long history~\cite{Schwarz:Uber_einen_Grenzubergang_durch_alternierendes_Verfahren,Lions:on_the_Schwarz_alternating_method_I}.  The first application of domain decomposition to the Helmholtz problem was proposed by~\cite{Despres:domain_decomposition_hemholtz}, which inspired the development of various domain decomposition algorithms, which are now classified as Schwarz algorithms\footnote{For a review on classical Schwarz methods see~\cite{Chan:Domain_decomposition_algorithms,Toselli:Domain_Decomposition_Methods_Algorithms_and_Theory}; and for other applications of domain decomposition methods for the Helmholtz equations, see~\cite{Bourdonnaye_Farhat_Roux:A_NonOverlapping_Domain_Decomposition_Method_for_the_Exterior_Helmholtz_Problem,Ghanemi98adomain,McInnes_Keyes:Additive_Schwarz_Methods_with_Nonreflecting_Boundary_Conditions_for_the_Parallel_Computation_of_Helmholtz_Problems,Collino:Domain_decomposition_method_for_harmonic_wave_propagation_a_general_presentation,Magoules:Application_of_a_domain_decomposition_with_Lagrange_multipliers_to_acoustic_problems_arising_from_the_automotive_industry,Boubendir:An_analysis_of_the_BEM_FEM_non_overlapping_domain_decomposition_method_for_a_scattering_problem,Astaneh_Guddati:A_two_level_domain_decomposition_method_with_accurate_interface_conditions_for_the_Helmholtz_problem}.}.
However, the convergence rate of such algorithms is strongly dependent on the boundary conditions prescribed at the interfaces between subdomains~\cite{Gander_Nataf:Optimized_Schwarz_Methods_without_Overlap_for_the_Helmholtz_Equation}.  The subsequent introduction of the optimized Schwarz framework in~\cite{Gander:Optimized_Schwarz_Methods}, which uses optimized boundary conditions to obtain good convergence, has inspired several competing approaches, including, but not limited to~\cite{Gander_Kwok:optimal_interface_conditiones_for_an_arbitrary_decomposition_into_subdomains,Geuzaine:A_quasi-optimal_nonoverlapping_domain_decomposition_algorithm_for_the_Helmholtz_equation,Gander_Zhang:Domain_Decomposition_Methods_for_the_Helmholtz_Equation:_A_Numerical_Investigation,Gander:Optimized_Schwarz_Methods_with_Overlap_for_the_Helmholtz_Equation,Gander:Optimized_Schwarz_Method_with_Two_Sided_Transmission_Conditions_in_an_Unsymmetric_Domain_Decomposition}.

Absorbing boundary conditions for domain decomposition schemes for elliptic problems are introduced by~\cite{Engquist_Zhao:Absorbing_boundary_conditions_for_domain_decomposition} and the first application of such techniques to the Helmholtz problem traces back to the AILU factorization~\cite{GanderNataf:ailu_for_hemholtz_problems_a_new_preconditioner_based_on_an_analytic_factorization}.
The sweeping preconditioner, introduced in~\cite{EngquistYing:Sweeping_H,EngquistYing:Sweeping_PML}, was the first algorithm to show that those ideas could yield algorithms with quasi-linear complexity, leading to several related algorithms with similar claims such as the source transfer preconditioner~\cite{Chen_Xiang:a_source_transfer_ddm_for_helmholtz_equations_in_unbounded_domain}, the rapidly converging domain decomposition~\cite{CStolk_rapidily_converging_domain_decomposition} and its extensions~\cite{Stolk:An_improved_sweeping_domain_decomposition_preconditioner_for_the_Helmholtz_equation}, the double sweep preconditioner~\cite{GeuzaineVion:double_sweep}, and the method of polarized traces~\cite{ZepedaDemanet:the_method_of_polarized_traces}. For an extensive review on sweeping-type methods we direct the reader to~\cite{GanderZhang2019}. 

To our knowledge, all domain decomposition methods and sweeping preconditioners have been based on layered domain decompositions, with one exception. In~\cite{Leng2019}, a $q\times q$ CDD is considered in the context of the source transfer method.  Their proposed strategy requires $p=q^2$ processors to obtain 
the sub-optimal parallel complexity of $\cO\big((N/\sqrt{p}) \log N\big)$, in both 2D and 3D. The primary difference between this method and our method is the sweeping strategy. The strategy employed in~\cite{Leng2019} is inspired by classical domain decomposition methods, while our method is inspired by sweeping preconditioners.  Consequently, our method inherits a superior sweeping strategy.


\subsection{Model Problem and Discretization}
\label{Sec::ModelProblem}

\begin{figure}[htp]
\centering
\begin{subfigure}[t]{0.45\textwidth}
\centering
\includegraphics[scale=0.15]{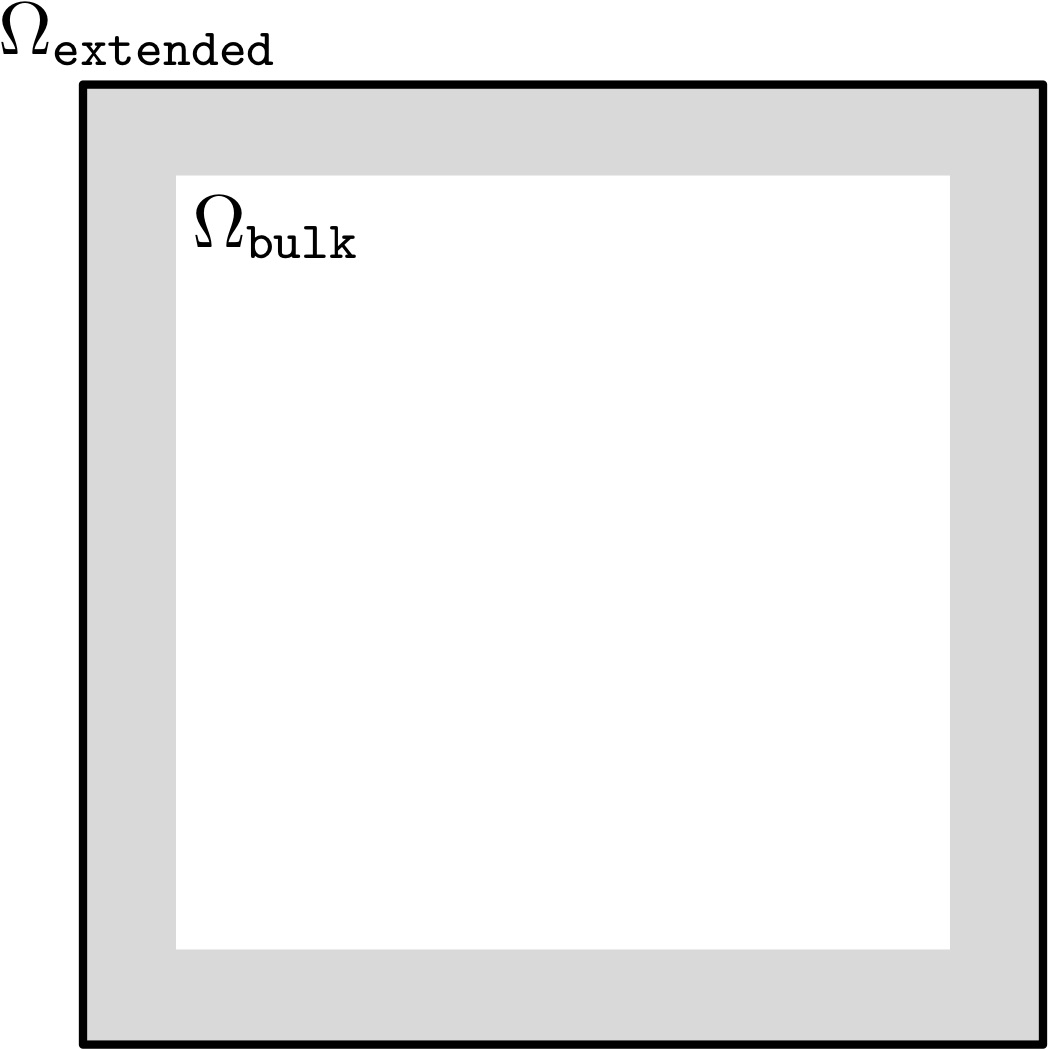}
\subcaption{}
\label{Fig::Domain::Extended}
\end{subfigure}
     \hfill
\begin{subfigure}[t]{0.45\textwidth}
\centering
\includegraphics[scale=0.15]{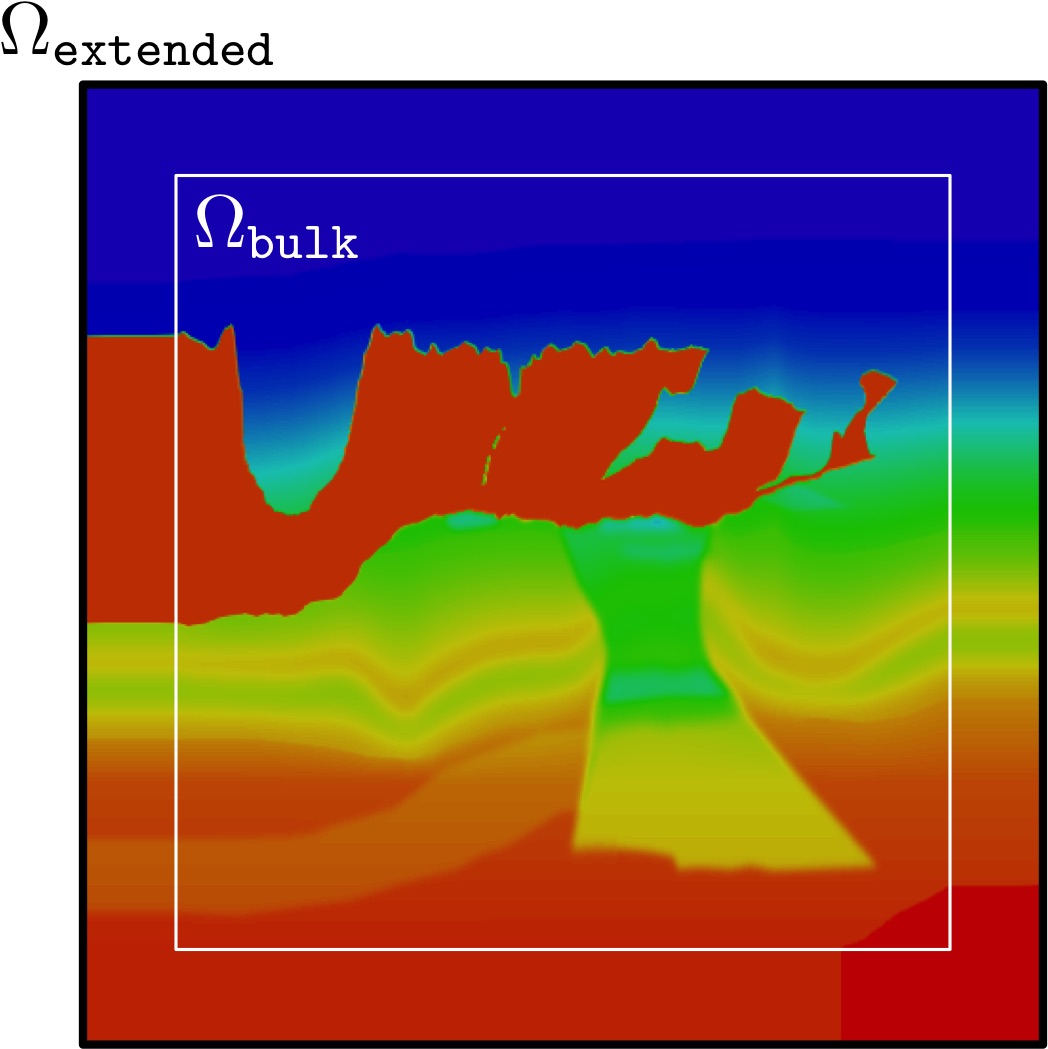}
\subcaption{}
\label{Fig::Domain::WaveSpeedExtension}
\end{subfigure}
\caption{(\subref{Fig::Domain::Extended}) $\Omega_{\texttt{bulk}}$ (white) and $\Omega_{\texttt{extended}}$ (white and gray) and  (\subref{Fig::Domain::WaveSpeedExtension}) part of the BP 2004~\cite{BP_model} wave speed ($\Omega_{\texttt{bulk}}$) extended to $\Omega_{\texttt{extended}}$ via the normal extension.}
\label{Fig::Domain}
\end{figure}

We consider problems formulated within the framework of~\eqref{HelmholtzBulk}, where we choose to model absorbing boundary conditions using perfectly matched layers (PMLs)~\cite{Berenger:PML,Johnson:PML}.  To preserve the solution in all of $\Omega_\texttt{bulk}$, using a PML requires the domain to be extended, which in turn also implies that the solution, material property, and source spaces must also be extended.  The extended domain, $\Omega_{\texttt{extended}}$, contains all of $\Omega_{\texttt{bulk}}$, as illustrated in Figure~\ref{Fig::Domain::Extended} and a new boundary value problem is formulated on $\Omega_{\texttt{extended}}$,
\begin{alignat}{1}
\label{HelmholtzExtended}\begin{aligned}-\operatorname{div}\left(\Lambda \nabla u\right) -\omega^2 m u &= f \quad \text{in }\Omega_{\texttt{extended}},\\
u &= 0 \quad \text{on }\partial \Omega_{\texttt{extended}},\end{aligned}
\end{alignat}
where, for brevity, we re-use the symbols $u$, $m$, and $f$ to represent the extended solution wavefield, slowness, and source distribution.  As is customary for PMLs in domains with varying wave speed, we assume that the source density, $f$, is extended by zero into $\Omega_{\texttt{extended}}$ and
the squared slowness, $m$, is extended into $\Omega_{\texttt{extended}}$ along the normal direction of $\partial\Omega_{\texttt{bulk}}$ in a constant fashion, as illustrated in Figure~\ref{Fig::Domain::WaveSpeedExtension}.
The PML-extended Helmholtz equation~\eqref{HelmholtzExtended} is reduced to~\eqref{HelmholtzBulk} in $\Omega_{\text{bulk}}$ by imposing that $\Lambda$ is the identity matrix in $\Omega_{\texttt{bulk}}$.  In the PML region, $\Omega_{\texttt{extended}}\backslash\Omega_{\texttt{bulk}}$, $\Lambda$ is a complex valued diagonal matrix which depends on the PML formulation, and $m$ and $f$ are 
complex-valued functions, obtained from imposing the PML.  Details of the precise formulation used in our developments are provided in Appendix~\ref{App::PML}.  Equation~\eqref{HelmholtzExtended} is effectively a complex-valued boundary-value problem in $\Omega_{\texttt{extended}}$ with homogeneous Dirichlet boundary conditions.  In the remainder of the discussion, we consider~\eqref{HelmholtzExtended} to be the canonical problem and therefore, for simplicity of notation, denote $\Omega_{\texttt{extended}}$ as $\Omega$.

When~\eqref{HelmholtzExtended} is discretized, we obtain the linear algebraic system
\begin{alignat}{1}
\label{LS} \boldA\boldu=\boldf
\end{alignat}
where $\boldA$ is the model-dependent system matrix, $\boldu$ is the solution vector, and $\boldf$ is the vector of the source density $f$.  In this paper we restrict our discussion of absorbing boundary conditions to PMLs and discretization to finite-difference methods.  These restrictions are merely to simplify the exposition: other transparent boundary conditions, such as absorbing layers or sponge layers, and other discretizations, such as higher-order finite differences and those derived from finite element methods, may be used in this framework.  For 2D problems we use standard 5-point finite difference stencils and in 3D, we use 9-point stencils. As a result, the discretization has $N$ global degrees-of-freedom and $n=N^{\frac{1}{2}}$  and $n=N^{\frac{1}{3}}$ degrees-of-freedom in each spatial dimension for two- and three-dimensional problems respectively.

For the PMLs we use a cubic PML profile function. As we increase the frequency of the problem, we do not increase the width of the PML. Instead, we choose the PML width so that the number of wavelengths is constant in the PML region, and increase the absorption constant logarithmically with the frequency. This is motivated by an analysis of the PML~\cite{Bramble2007}, where this choice is rigorously justified
for a more complex PML-profile. In our work, we employ the same strategy for the cubic PML-profile and obtain satisfactory results. Details on the construction of the linear system and the PMLs are given in Appendix~\ref{App::PML}. 

\subsection{Continuous Polarization}
\label{Sec::Background::Continuous}

The method of polarized traces was introduced as a solver~\cite{ZepedaDemanet:the_method_of_polarized_traces} and then as a preconditioner~\cite{Zepeda_Hewett_Demanet:Preconditioning_the_2D_Helmholtz_equation_with_polarized_traces} for the linear system~\eqref{LS}. At its core, the method of polarized traces spatially subdivides the discrete degrees-of-freedom in $\Omega$ into layers and computes an approximate solution to the global wavefield by sweeping over the layers and solving a local discrete half-space problem in each layer. Following~\cite{ZepedaDemanet:the_method_of_polarized_traces}, the solutions of the half-space problems are called polarized wavefields. In this work, we make extensive use of this concept and therefore provide a brief review, in the continuous setting, in this section.  In Section~\ref{Sec::Background::Discrete}, we present a similar treatment in the discrete context.

\begin{figure}
\begin{tabular}{c|c}
\begin{subfigure}[b]{0.45\textwidth}
\centering\includegraphics[scale=0.3]{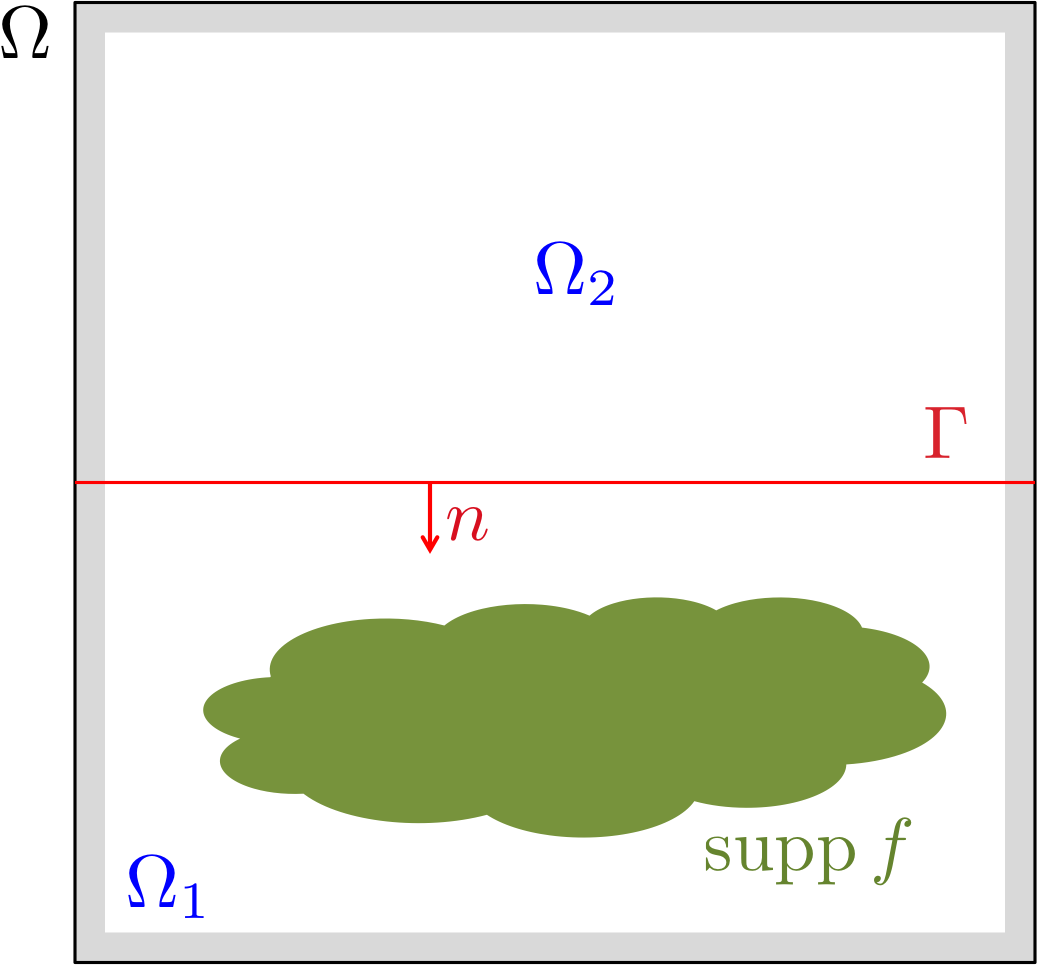}
\subcaption{}
\label{Fig::HalfSpace}
\end{subfigure} & \begin{subfigure}[b]{0.45\textwidth}
\centering\includegraphics[scale=0.3]{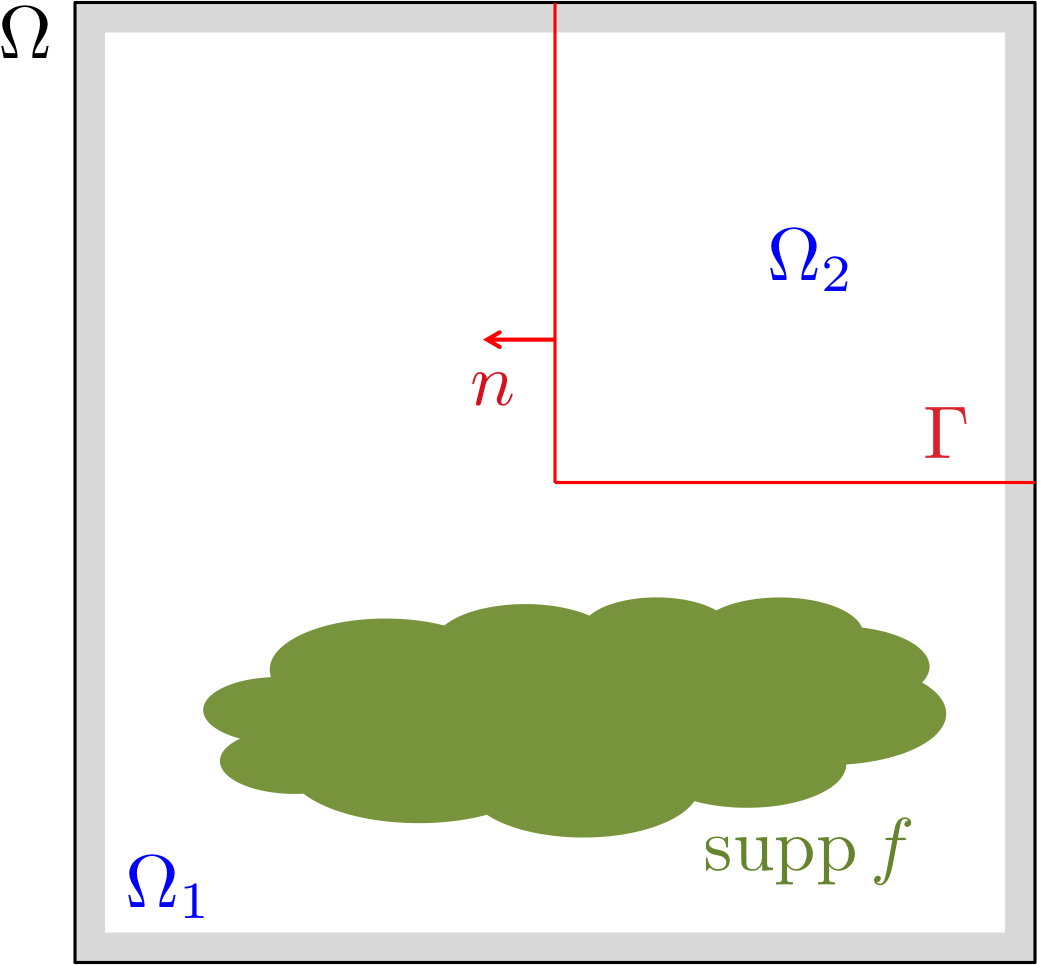}
\subcaption{}
\label{Fig::QuarterSpace}
\end{subfigure}
\end{tabular}
\caption{A schematic representation of the truncated (\subref{Fig::HalfSpace}) half- and (\subref{Fig::QuarterSpace}) quarter-space problem.}
\end{figure}

Consider the boundary-value problem~\eqref{HelmholtzExtended} with a source density function $f$ supported only in $\Omega_{\texttt{bulk}}$. Let $\Gamma$ be an interface dividing $\Omega$ into two regions, $\Omega_1$ and $\Omega_2$, such that the support of $f$ lies entirely within $\Omega_1$. For example, $\Gamma$ could be a straight line (Figure~\ref{Fig::HalfSpace}) or an L-shaped line (Figure~\ref{Fig::QuarterSpace}).

The wavefield in $\Omega_2$ can be computed from the representation formula
\begin{alignat}{1}\label{RepForm}
u(x) = \int_{\Gamma}\mu(y)G(x,y)ds_y-\int_{\Gamma}\lambda(y)\left[n(y)\cdot\left(\Lambda(y)\nabla_yG(x,y)\right)\right]ds_y \quad x\in\Omega_2,
\end{alignat}
where
\begin{alignat*}{1}
\lambda=u|_{\Gamma}, \quad\text{and}\quad \mu=\left.\left[n\cdot\left(\Lambda\nabla u\right)\right]\right|_{\Gamma}
\end{alignat*}
are the Dirichlet and Neumann traces on $\Gamma$, and $G(x,y)$ is the Green's function corresponding to the problem~\eqref{HelmholtzExtended}, i.e., for $x\in\Omega$
\begin{alignat*}{2}
-\operatorname{div}_y\left(\Lambda(y) \nabla_y G(y,x)\right) -\omega^2 m(y) G(y,x) &= \delta(x-y) \quad &y\in\Omega,\\
G(y,x) &= 0 \quad &y\in\partial \Omega.
\end{alignat*}
This formula directly follows from the divergence theorem and properties of the Green's function, assuming that $\Lambda$ is a diagonal matrix.

Equation~\eqref{RepForm} requires knowledge of the Dirichlet and Neumann traces of the solution almost everywhere on $\Gamma$.
Thus, corners in $\Gamma$ are admissible, as in the quarter-space problem, even though the normal is not uniquely defined.  Using~\eqref{RepForm}, the solution $u$ can then be computed on $\Omega_2$. Extending equation~\eqref{RepForm} to all of $\Omega$,
\begin{alignat}{1}\label{UnboundedPolWF}
U(x)= \int_{\Gamma}\mu(y)G(x,y)ds_y-\int_{\Gamma}\lambda(y)\left[n(y)\cdot\left(\Lambda(y)\nabla_yG(x,y)\right)\right]ds_y \quad x\in\Omega,
\end{alignat}
it can be shown that $U$ vanishes on $\Omega_1$,
\begin{alignat}{1}
\label{AnnihilationCondition}
U(x)=0\quad \text{for }x\in\Omega_1.
\end{alignat}
Following the terminology of~\cite{ZepedaDemanet:the_method_of_polarized_traces}, we call~\eqref{AnnihilationCondition} the annihilation condition, $\Gamma$ the polarization interface, and $U$ a polarized wavefield. A proof of the annihilation condition~\eqref{AnnihilationCondition} is provided in Appendix~\ref{App::AnCond}.

\subsection{Discrete Polarization}
\label{Sec::Background::Discrete}
\begin{figure}[htp]
\centering
\includegraphics[scale=0.3]{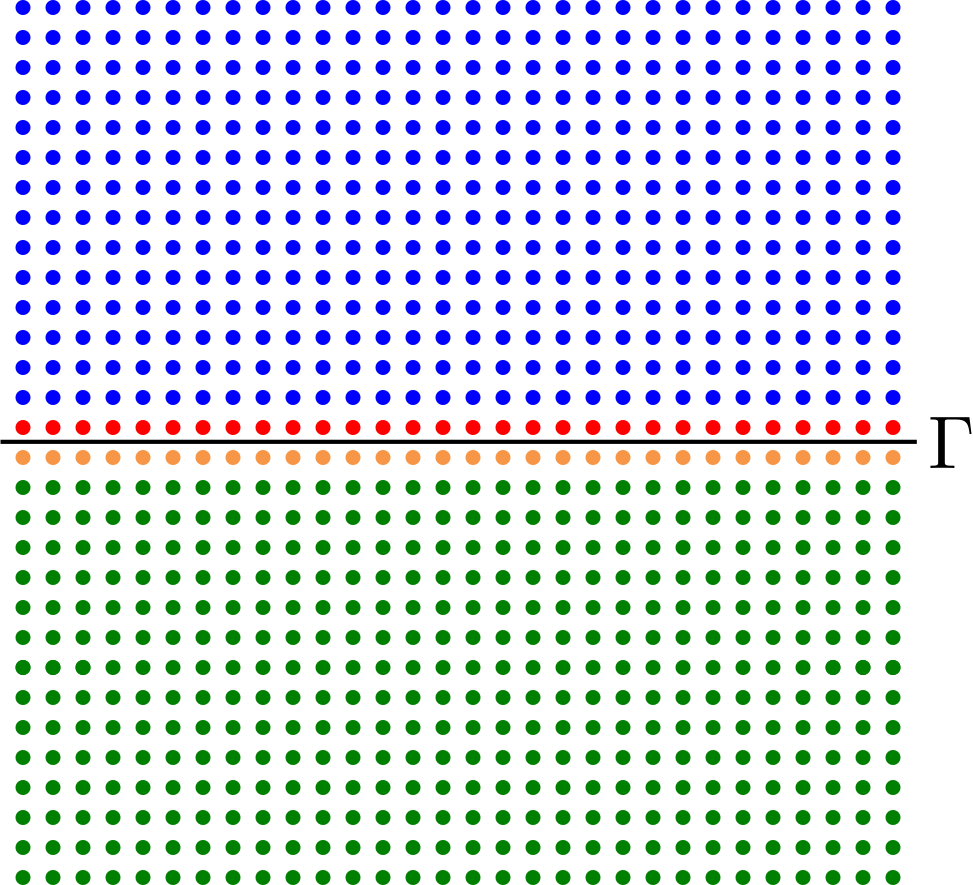}
\quad
\includegraphics[scale=0.3]{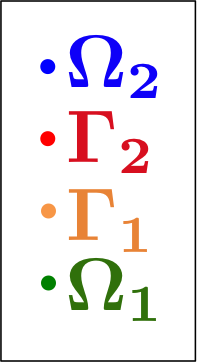}
\quad
\includegraphics[scale=0.3]{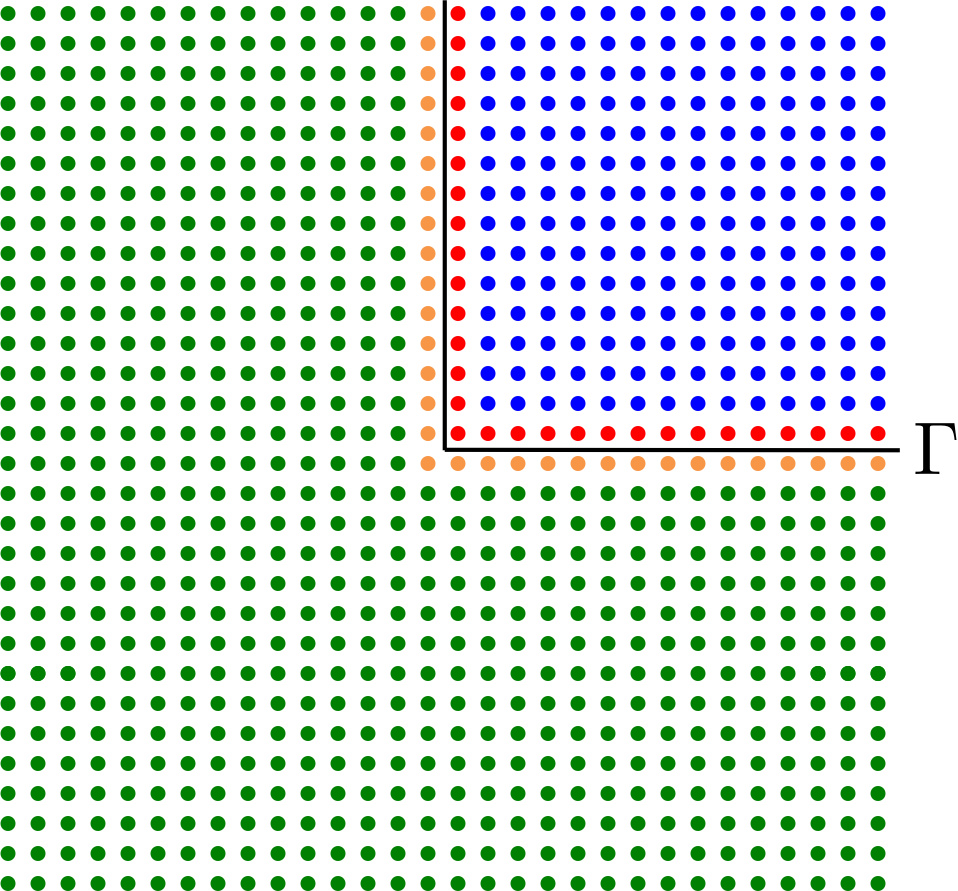}
\caption{A schematic representation of four sets of degrees-of-freedom $\boldOmegaOne$, $\boldOmegaTwo$, $\boldGammaOne$ and $\boldGammaTwo$.}
\label{Fig::DiscHalfQuarterSpace}
\end{figure}

A discrete counterpart of the polarized wavefield $U$ can be derived by constructing a discrete solution that satisfies the discrete analogue to the annihilation condition in~\eqref{AnnihilationCondition}. Consider the discretization points corresponding to degrees-of-freedom in~\eqref{LS}, as well as an interface $\Gamma$ which does not intersect with any discretization point, as illustrated in Figure~\ref{Fig::DiscHalfQuarterSpace} for both the half-space and quarter-space subdomains.  Given a discretization-dependent distance $\delta$ (e.g., $\delta = 1$ for a classical 5-point finite difference stencil in 2D), such an interface divides the degrees-of-freedom into four sets, as labeled in
Figure~\ref{Fig::DiscHalfQuarterSpace}: 
\begin{enumerate}
    \item $\boldGammaOne$, the set of all degrees-of-freedom physically contained in $\Omega_1$ and $\delta$-adjacent to $\Gamma$;
    \item $\boldOmegaOne$, the set of all degrees-of-freedom physically contained in $\Omega_1$, excluding $\boldGammaOne$;
    \item $\boldGammaTwo$, the set of all degrees-of-freedom physically contained in $\Omega_2$ and $\delta$-adjacent to $\Gamma$;
    \item $\boldOmegaTwo$, the set of all degrees-of-freedom physically contained in $\Omega_2$, excluding $\boldGammaTwo$.
\end{enumerate}

Upon reordering the degrees-of-freedom with respect to these four sets, the discrete system~\eqref{LS} can be rewritten as
\begin{alignat}{1}\label{ReorderedLS}
\begin{pmatrix}
\boldA_{\boldOmegaOne\boldOmegaOne} & \boldA_{\boldOmegaOne\boldGammaOne} & \boldZero & \boldZero\\
\boldA_{\boldGammaOne\boldOmegaOne} & \boldA_{\boldGammaOne\boldGammaOne} & \boldA_{\boldGammaOne\boldGammaTwo} & \boldZero\\
\boldZero & \boldA_{\boldGammaTwo\boldGammaOne} & \boldA_{\boldGammaTwo\boldGammaTwo} & \boldA_{\boldGammaTwo\boldOmegaTwo}\\
\boldZero & \boldZero & \boldA_{\boldOmegaTwo\boldGammaTwo} & \boldA_{\boldOmegaTwo\boldOmegaTwo}
\end{pmatrix}
\begin{pmatrix}
\boldu_{\boldOmegaOne}\\
\boldu_{\boldGammaOne}\\
\boldu_{\boldGammaTwo}\\
\boldu_{\boldOmegaTwo}
\end{pmatrix}
=
\begin{pmatrix}
\boldf_{\boldOmegaOne}\\
\boldf_{\boldGammaOne}\\
\boldf_{\boldGammaTwo}\\
\boldf_{\boldOmegaTwo}
\end{pmatrix},
\end{alignat}
where $\boldA_{\boldGammaOne\boldOmegaOne}$ denotes the slice of the matrix $\boldA$ generated from the rows corresponding to $\boldGammaOne$ and the columns corresponding to $\boldOmegaOne$.  The remaining submatrices are similarly generated.  In the same way, the vector $\boldu_\boldOmegaTwo$ (and similarly all other vectors) denotes the slice of the vector $\boldu$ with respect to $\boldOmegaTwo$. As in the continuous case, we consider the case without sources in $\Omega_2$, i.e., $\boldf_{\boldOmegaTwo}=\boldZero$ and $\boldf_{\boldGammaTwo}=\boldZero$. Then, due to the invertibility of $\boldA$, it is easy to see that the solution $\boldU$ of the linear system
\begin{alignat}{1}\label{DiscPol2}
\begin{pmatrix}
\boldA_{\boldOmegaOne\boldOmegaOne} & \boldA_{\boldOmegaOne\boldGammaOne} & \boldZero & \boldZero\\
\boldA_{\boldGammaOne\boldOmegaOne} & \boldA_{\boldGammaOne\boldGammaOne} & \boldA_{\boldGammaOne\boldGammaTwo} & \boldZero\\
\boldZero & \boldA_{\boldGammaTwo\boldGammaOne} & \boldA_{\boldGammaTwo\boldGammaTwo} & \boldA_{\boldGammaTwo\boldOmegaTwo}\\
\boldZero & \boldZero & \boldA_{\boldOmegaTwo\boldGammaTwo} & \boldA_{\boldOmegaTwo\boldOmegaTwo}
\end{pmatrix}
\boldU
=
\begin{pmatrix}
\boldZero\\
\boldA_{\boldGammaOne\boldGammaTwo}\boldu_{\boldGammaTwo}\\
-\boldA_{\boldGammaTwo\boldGammaOne}\boldu_{\boldGammaOne}\\
\boldZero
\end{pmatrix}
\end{alignat}
satisfies 
\begin{alignat*}{1}
\boldU=
\begin{pmatrix}
\boldZero\\
\boldZero\\
\boldu_{\boldGammaTwo}\\
\boldu_{\boldOmegaTwo}
\end{pmatrix}.
\end{alignat*}
Consequently, the discrete counterpart of the polarized wavefield $\boldU$ is
\begin{alignat}{1}
\label{DiscPol}
\boldU=\begin{pmatrix}
\boldA_{\boldOmegaOne\boldOmegaOne} & \boldA_{\boldOmegaOne\boldGammaOne} & \boldZero & \boldZero\\
\boldA_{\boldGammaOne\boldOmegaOne} & \boldA_{\boldGammaOne\boldGammaOne} & \boldA_{\boldGammaOne\boldGammaTwo} & \boldZero\\
\boldZero & \boldA_{\boldGammaTwo\boldGammaOne} & \boldA_{\boldGammaTwo\boldGammaTwo} & \boldA_{\boldGammaTwo\boldOmegaTwo}\\
\boldZero & \boldZero & \boldA_{\boldOmegaTwo\boldGammaTwo} & \boldA_{\boldOmegaTwo\boldOmegaTwo}
\end{pmatrix}^{-1}
\begin{pmatrix}
\boldZero\\
\boldA_{\boldGammaOne\boldGammaTwo}\boldu_{\boldGammaTwo}\\
-\boldA_{\boldGammaTwo\boldGammaOne}\boldu_{\boldGammaOne}\\
\boldZero
\end{pmatrix}.
\end{alignat}

From the right-hand side of~\eqref{DiscPol2}, one can easily see that knowledge of both $\boldu_\boldGammaOne$ and $\boldu_\boldGammaTwo$ is required in order to compute $\boldU$. By construction, these sets contain information about the discrete wavefield and its normal derivative in the vicinity of $\Gamma$. Thus, as with the continuous case, the discrete case also requires information about the Dirichlet and Neumann traces to compute the polarized wavefield $\boldU$ in $\boldOmegaTwo$.  In fact,~\cite{ZepedaDemanet:the_method_of_polarized_traces} demonstrates, using similar techniques, that a discrete counterpart of the representation formula~\eqref{RepForm} can be derived.

Finally, our technique is easily extended to more general discretization techniques, as long as they allow for a reordering of the degrees-of-freedom to the block-tridiagonal system~\eqref{ReorderedLS}.  Depending on the discretization, change in the selection of the sets $\boldGammaOne$ and $\boldGammaTwo$ may be required.  Similar schemes have been applied  for many different discretizations such as high-order finite difference methods~\cite{Zepeda2019}, finite element methods~\cite{ZepedaDemanet:A_short_note_on_the_nested-sweep_polarized_traces_method_for_the_2D_Helmholtz_equation}, enriched finite element methods~\cite{Fang_Qian_Zepeda_Zhao:Learning_Dominant_Wave_Directions_For_Plane_Wave_Methods_For_High_Frequency_Helmholtz_Equations}, discontinuous Galerkin methods~\cite{Taus2016}, and integral representations~\cite{ZepedaZhao:Fast_Lippmann_Schwinger_solver}.  For example, for higher-order finite difference methods the stencils centered at the points in $\boldOmegaTwo$ (respecting $\boldGammaTwo$, $\boldGammaOne$, or $\boldOmegaOne$) cannot involve discretization points in $\boldGammaOne$ (respecting $\boldOmegaOne$, $\boldOmegaTwo$, $\boldGammaTwo$). This is easily enforced by defining an appropriate $\delta$, e.g., $\delta=2$ for a 9-point, 5x5, stencil in 2D.  Similarly, in finite element or discontinuous Galerkin methods, the sparsity of the system matrices can be exploited in order to obtain suitable sets of degrees-of-freedom.

\subsection{Organization}
In Section~\ref{Sec::LSweeps}, we introduce the algorithm to compute an approximate global solution $u$ of the boundary value problem~\eqref{HelmholtzExtended} with constant squared slowness $m$. We introduce the algorithm first on the continuous level in Section~\ref{Sec::LSweeps::Cont} and then extend it to the discrete level in Section~\ref{Sec::LSweeps::Disc}. In Section~\ref{Sec::Complexity}, we show how the algorithm can be used to precondition the linear system~\eqref{LS}. This opens the possibility of using the algorithm as part of an optimally parallel scaling solver based on a preconditioned GMRES method. We conclude Section~\ref{Sec::Complexity} with a complexity analysis of this solver with regards to computational and communication effort. In Section~\ref{Sec::LSweeps::Heterogeneities}, we analyze the effects of heterogeneous wave speeds on the effectiveness of the preconditioner. In Section~\ref{Sec::NumericalExamples}, we provide several numerical examples in two- and three-dimensions for constant and non-constant wave speeds to corroborate all claims. The paper is concluded by a discussion where we briefly summarize our results and discuss possible extensions. Additional details and pseudo-code for the proposed algorithms are provided in the Appendices.

\section{L-sweeps: Reconstruction of wavefields}
\label{Sec::LSweeps}
In this section, we introduce the algorithm to compute the global solution $u$ of the boundary value problem~\eqref{HelmholtzExtended}. This solution is obtained from local solutions of local problems defined over a CDD. The concept of polarization introduced in Section~\ref{Sec::Introduction} plays a crucial role in this procedure. We introduce the algorithm for constant wave speeds at the continuous and the discrete level in Sections~\ref{Sec::LSweeps::Cont} and~\ref{Sec::LSweeps::Disc}. The procedure can be applied to problems with non-constant wave speeds in an analogous fashion. The effects of these heterogeneous wave speeds, in particular discontinuous ones, are addressed in Section~\ref{Sec::LSweeps::Heterogeneities}.

\subsection{Continuous formulation}
\label{Sec::LSweeps::Cont}
\begin{figure}[htp]
\centering
\includegraphics[scale=0.5]{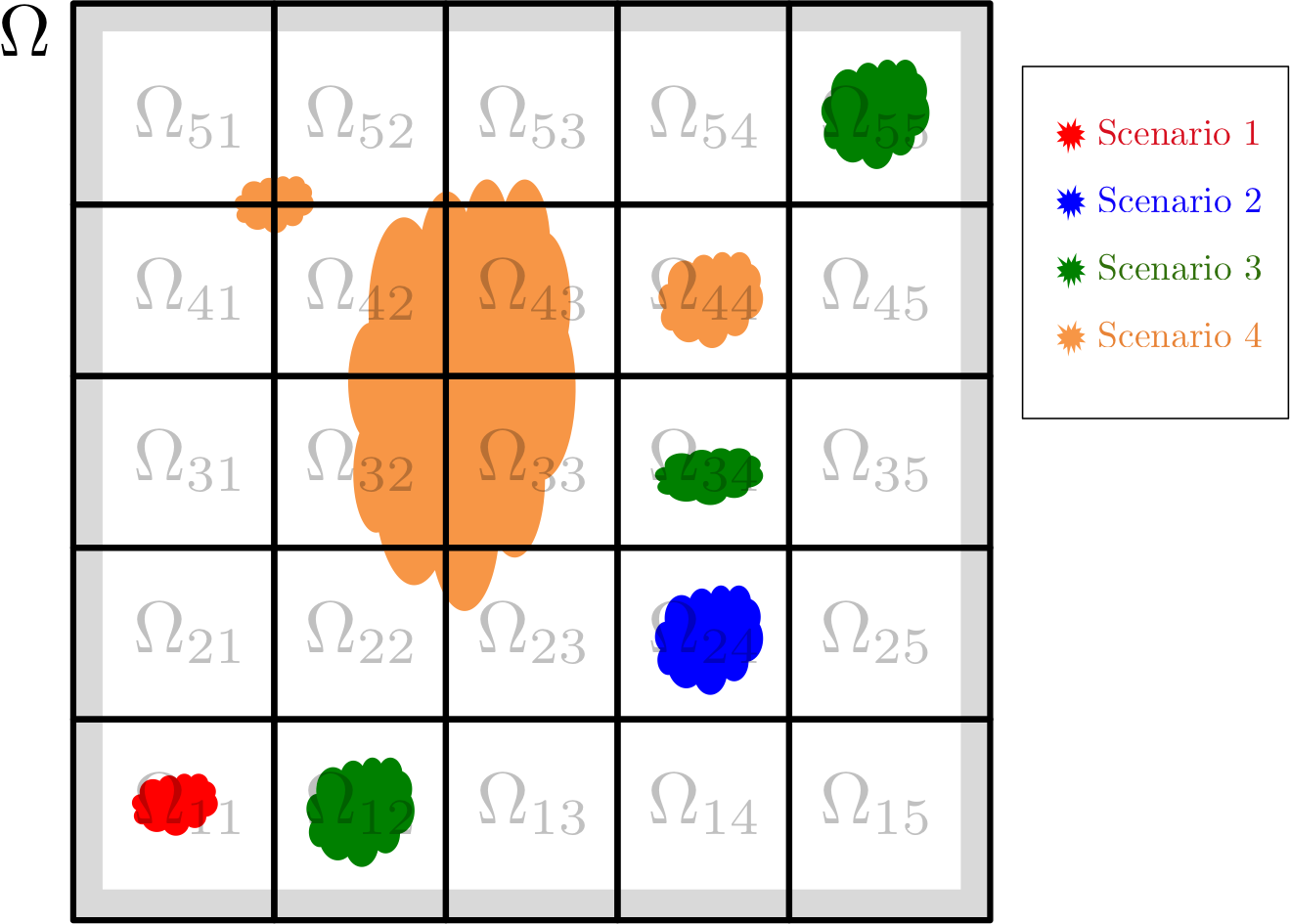}
\caption{The domain decomposition, where the colored regions show examples of source distributions for each of the four scenarios.}
\label{Fig::CDD}
\end{figure}

Consider a decomposition of $\Omega$ into a CDD\footnote{Note that the CDD is chosen so that neglecting the PML regions (shown in gray in Figure~\ref{Fig::CDD}) each subdomain has the same size.}, with $q$ rows and $r$ columns, of non-overlapping open subdomains.  For example, see Figure~\ref{Fig::CDD}, where $q=r=5$. We first define the local problems associated with each subdomain, $\Omega_{ij}$. To this end, we define extended domains $\Omega_{ij,\texttt{bulk}}^\varepsilon$ and $\Omega_{ij,\texttt{extended}}^\varepsilon$. The domain $\Omega_{ij,\texttt{bulk}}^\varepsilon$ is obtained from extending $\Omega_{ij}$ by a $\varepsilon$-layer along interior edges of the CDD, and $\Omega_{ij,\texttt{extended}}^\varepsilon$ is an extension of $\Omega_{ij,\texttt{bulk}}^\varepsilon$ to impose absorbing boundary conditions via PMLs. In $\Omega_{ij,\texttt{extended}}^\varepsilon$ we define the local squared slowness $m_{ij}:=m|_{\Omega_{ij,\texttt{extended}}^\varepsilon}$. For a  given source density $f_{ij}$ in $\Omega_{ij,\texttt{extended}}^\varepsilon$, we define the local problem 
\begin{alignat*}{1}
-\operatorname{div}\left(\Lambda_{ij} \nabla u_{ij}\right) -\omega^2 m_{ij} u_{ij} &= f_{ij} \quad \text{in }\Omega_{ij,\texttt{extended}}^\varepsilon,\\
u_{ij} &= 0 \quad \text{on }\partial \Omega_{ij,\texttt{extended}}^\varepsilon, 
\end{alignat*}
where $u_{ij}$ denotes the local solution. Due to the local PML, $\Lambda_{ij}$ is a complex-valued diagonal matrix, and $m_{ij}$ and $f_{ij}$ are complex-valued functions in the PML region $\Omega_{ij,\texttt{extended}}^\varepsilon\backslash\Omega_{ij,\texttt{bulk}}^\varepsilon$.  Further details of the PML formulation are provided in Appendix~\ref{App::PML}. For the sake of brevity we denote the PML-adjusted squared slowness and source density still as $m_{ij}$ and $f_{ij}$, respectively.  Note that in contrast to~\eqref{HelmholtzExtended}, we do not set the squared slowness to be constant using a normal extension on $\partial \Omega_{ij,\texttt{bulk}}^\varepsilon$ in the PML region, rather we use the squared slowness inherited from the global problem in these regions. This definition makes our approach more accurate, since we use the local problem to compute sections of the global problem.

Note that the construction of the local problems consists of two subsequent extensions of $\Omega_{ij}$: we first add an additional $\varepsilon$-layer around $\Omega_{ij}$, and then further extend the subdomain by a PML region. While the latter is clearly needed to avoid artificial reflections in the local solutions, the former appears to be ad-hoc at this stage, but for reasons that will be explained in the sequel, the first extension by an $\varepsilon$-layer is crucial for a consistent exchange of information between subdomains.
The resulting local problems are Dirichlet boundary-value problems defined on the extension $\Omega_{ij,\texttt{extended}}^\varepsilon$ of $\Omega_{ij}$.
As before, we denote the domain $\Omega_{ij,\texttt{extended}}^\varepsilon$ as $\Omega_{ij}^\varepsilon$. Examples of these local problems are illustrated in Figure~\ref{Fig::LocProbs}.

In what follows, we introduce our method for computing global solutions to~\eqref{HelmholtzExtended} by considering four scenarios, illustrated in Figure~\ref{Fig::CDD}, each increasingly more general: 
\begin{enumerate}
    \item the source density is supported in the interior of a corner subdomain of the CDD,
    \item the source density is supported in the interior of an arbitrary subdomain,
    \item the source density is supported in the interior of an arbitrary number of subdomains such that its support does not intersect the skeleton of the CDD, and
    \item the source density has arbitrary support in $\Omega_{\texttt{bulk}}$.
\end{enumerate}
In the subsequent developments, we consider these scenarios using source distributions constructed from unions of point sources. However, the developments do not depend on an assumption that the source distributions are point-sources -- any source density is admitted as long as it satisfies the conditions on its support. In particular, scenario 4 allows for arbitrary source distributions, including those which intersect the CDD skeleton, making our approach widely applicable.

\subsubsection{Scenario 1: A source density supported in a corner subdomain}
\label{Sec::LSweeps::Cont::CornerSrc}
Without loss of generality, we introduce the algorithm for a source density supported in $\Omega_{11}$.  Source densities supported in any other corner subdomain can be constructed in an analogous way.  The computation of the global solution is performed in three stages:
\begin{enumerate}
	\item compute an approximation of the global solution in $\Omega_{11}$ using the local problem associated with $\Omega_{11}$,
	\item compute the approximate global solution in the first row and column of the CDD, and,
	\item compute the global solution in the rest of the subdomains.
\end{enumerate}
Ultimately, using these three stages, we are able to compute the global solution, up to PML induced errors, using three sweeps: one vertical, one horizontal, and one diagonal.
\begin{figure}[htp]
\centering
\begin{tabular}{c|c}
\begin{subfigure}[b]{0.55\textwidth}
\centering\includegraphics[scale=0.6]{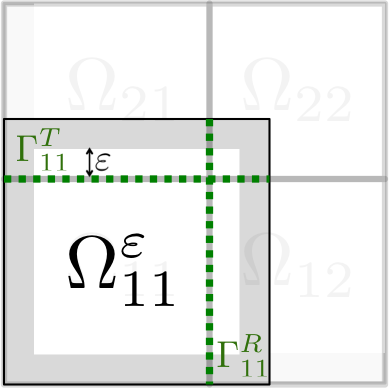}
\subcaption{Example of a local problem in stage 1.}
\label{Fig::LocProbs::LocProb11}
\end{subfigure} & \begin{subfigure}[b]{0.35\textwidth}
\centering\includegraphics[scale=0.125]{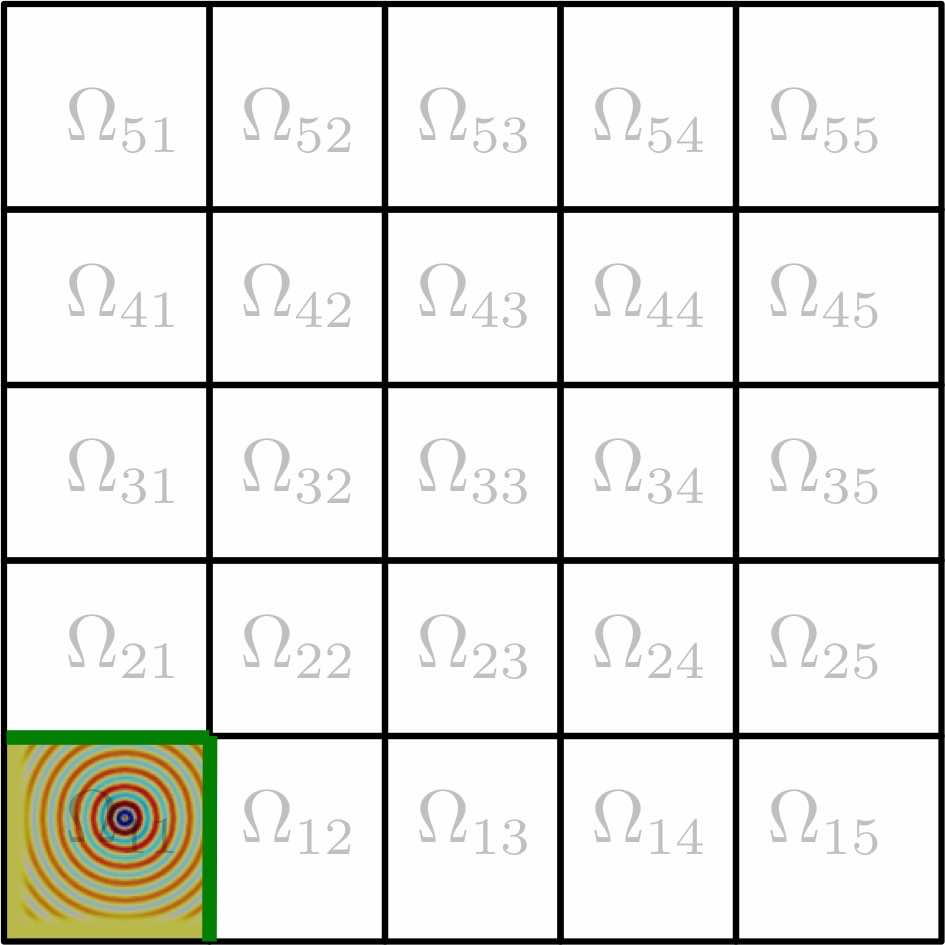}
\subcaption{Solution after stage 1.}
\label{Fig::LocProbs::WaveField1}
\end{subfigure} \\
\vspace{0.5cm} & \\
\begin{subfigure}[b]{0.55\textwidth}
\centering\includegraphics[scale=0.6]{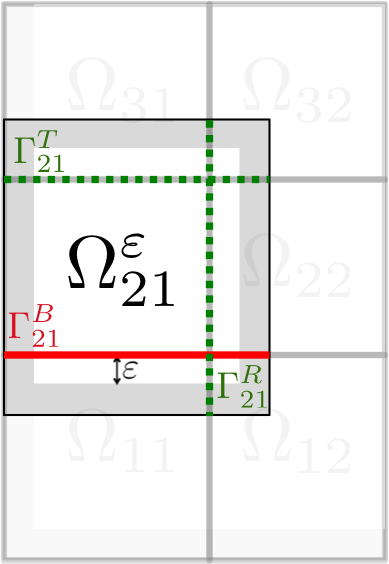}
\subcaption{Example of a local problem in stage 2.}
\label{Fig::LocProbs::LocProb21}
\end{subfigure} & \begin{subfigure}[b]{0.35\textwidth}
\centering\includegraphics[scale=0.125]{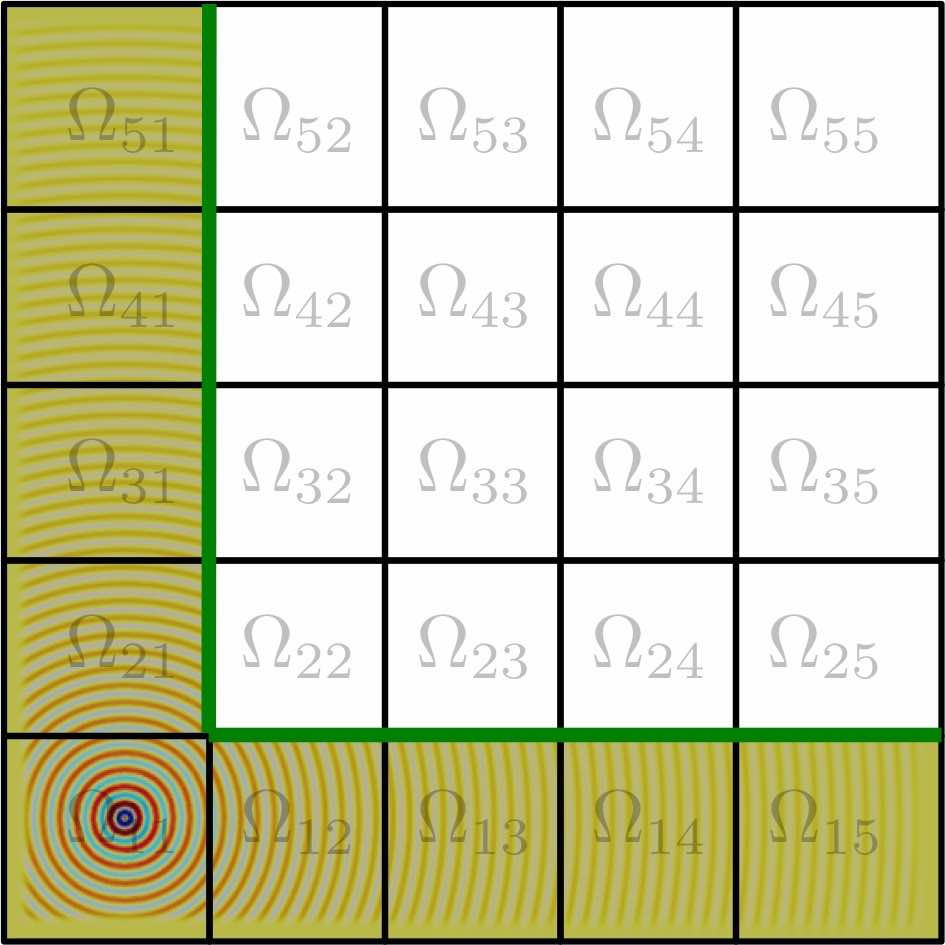}
\subcaption{Solution after stage 2.}
\label{Fig::LocProbs::WaveField2}
\end{subfigure} \\
\vspace{0.5cm} & \\
\begin{subfigure}[b]{0.55\textwidth}
\centering\includegraphics[scale=0.6]{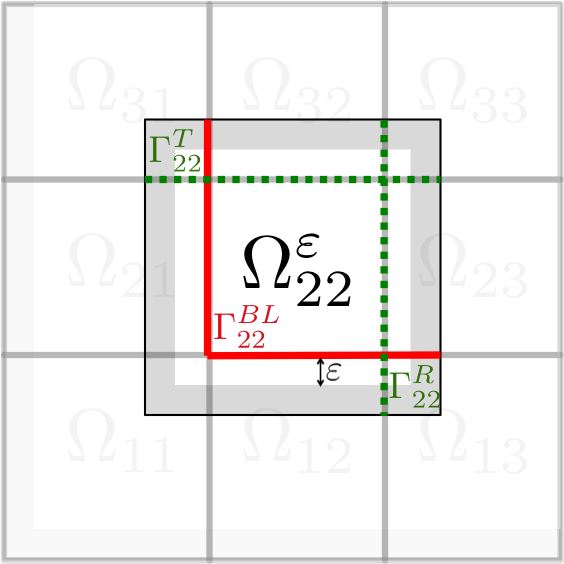}
\subcaption{Example of a local problem in stage 3.}
\label{Fig::LocProbs::LocProb22}
\end{subfigure} & \begin{subfigure}[b]{0.35\textwidth}
\centering\includegraphics[scale=0.125]{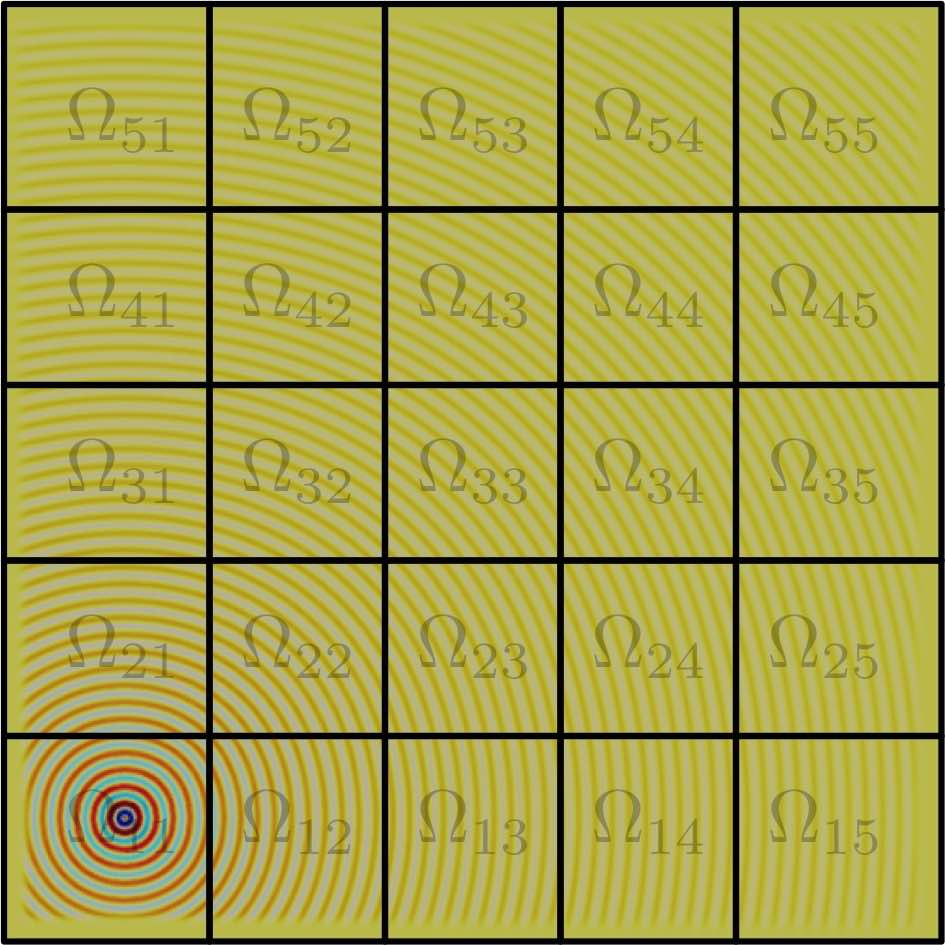}
\subcaption{Solution after stage 3.}
\label{Fig::LocProbs::WaveField3}
\end{subfigure}
\end{tabular}
\caption{Illustration of local problems and the computed solution after each stage of a sweep. The green lines in the wavefields depict extracted trace information used in the next stage.}
\label{Fig::LocProbs}
\end{figure}

\bigskip
\noindent
{\it Stage 1: Local solution}\\
The local solution $u_{11}$ defined over $\Omega_{11}^\varepsilon$ can be computed by solving the local problem associated with $\Omega_{11}$ with the source density $f_{11}:=f|_{\Omega_{11}^\varepsilon}$. Since $f$ is supported in $\Omega_{11}$ and the wave speed is constant, the global solution $u$ restricted to $\Omega_{11}$ coincides with $u_{11}$ up to errors induced by the PMLs. We therefore simply set the global solution $u$ to $u_{11}$ in $\Omega_{11}$: $u|_{\Omega_{11}}:=u_{11}|_{\Omega_{11}}$.

We define the straight lines containing the top and right boundary of $\Omega_{11}$ to be $\Gamma^T_{11}$ and $\Gamma^R_{11}$, as depicted in Figure~\ref{Fig::LocProbs::LocProb11}. On $\Gamma_{11}^T$ and $\Gamma_{11}^R$, we extract the Dirichlet and Neumann traces of $u_{11}$:
\begin{alignat}{3}
\label{Traces}\begin{aligned}
\lambda_{11}^{T}&:=u_{11}|_{\Gamma_{11}^T},& \mu_{11}^{T}&:=\left[n\cdot\left(\Lambda_{11}\nabla u_{11}\right)\right]|_{\Gamma_{11}^T}, \\
\lambda_{11}^{R}&:=u_{11}|_{\Gamma_{11}^R},& \mu_{11}^{R}&:=\left[n\cdot\left(\Lambda_{11}\nabla u_{11}\right)\right]|_{\Gamma_{11}^R}.
\end{aligned}
\end{alignat}
These traces provide the necessary information to compute good approximations of the global solution in the neighboring subdomains $\Omega_{21}$ and $\Omega_{12}$ in stage 2. The solution computed after stage 1 is shown in Figure~\ref{Fig::LocProbs::WaveField1}, where the trace information that is used in stage 2 is shown in green. 

\bigskip
\noindent
{\it Stage 2: Global solution in the same row/column}\\
By construction, $\Gamma_{11}^T$ lies entirely inside $\Omega_{21}^\varepsilon$ and the source density $f$ is zero on $\Gamma_{11}^T$. Thus, the traces $\lambda_{11}^{T}$ and $\mu_{11}^{T}$ can be used to compute a polarized wavefield in $\Omega_{21}^\varepsilon$ using $\Gamma_{21}^B:=\Gamma_{11}^T$ as the polarization interface, and
\begin{alignat}{1}
\label{Eq:ProbStraightTrace}
u_{21}(x)=\int_{\Gamma_{21}^B}\mu_{11}^{T}(y)G_{21}(x,y)ds_y-\int_{\Gamma_{21}^B}\lambda_{11}^{T}(y)\left[n(y)\cdot\left(\Lambda_{21}\nabla_y G_{21}(x,y)\right)\right]ds_y.
\end{alignat}
Here, $G_{21}(x,y)$ is the Green's function corresponding to the local problem defined on $\Omega_{21}^\varepsilon$.

Following the same reasoning used for $u_{11}$ in stage 1, the global solution restricted to $\Omega_{21}$ coincides with $u_{21}$, up to errors induced by the PMLs, and we set $u|_{\Omega_{21}}:=u_{21}|_{\Omega_{21}}$. We also extract the Dirichlet and Neumann traces $\lambda_{21}^{T}$ and $\mu_{21}^{T}$ of $u_{21}$ on $\Gamma_{21}^T$ in the same way as in~\eqref{Traces} and repeat the process to compute an approximation of the global solution in the entire first column of the CDD. In addition, for each of these polarized wavefields, we extract the Dirichlet and Neumann traces $\lambda_{i1}^{R}$ and $\mu_{i1}^{R}$ on $\Gamma_{i1}^R$. These traces are needed in stage 3 to compute approximations of the global solution in the rest of the subdomains.
Similarly, we compute the global solution in the first row of the CDD and extract the traces $\lambda_{1j}^{T}$ and $\mu_{1j}^{T}$ on $\Gamma_{1j}^T$ for further use in stage 3. The solution after stage 2 and the trace information extracted for stage 3 are shown in Figure~\ref{Fig::LocProbs::WaveField2}. 

For the extension of the solution into the neighboring subdomains to be accurate, the local problems associated with $\Omega_{i1}$ and $\Omega_{(i+1)1}$ need to coincide in an $\varepsilon$-tube around $\Gamma_{i1}^T=\Gamma_{(i+1)1}^B$. This is ensured by the extension of the local problem by the additional $\varepsilon$-layer.

\bigskip
\noindent
{\it Stage 3: Global solution in the remaining subdomains}\\
Consider $\Omega_{22}$, where, by construction, $\Gamma_{22}^B:=\Gamma_{12}^T$ and $\Gamma_{22}^L:=\Gamma_{21}^R$ lie entirely inside $\Omega_{22}^\varepsilon$. We combine these lines to form an L-shaped line, $\Gamma_{22}^{BL}$, so that $\Omega_{22}$ is entirely contained in the quadrant defined by $\Gamma_{22}^{BL}$ as shown in Figure~\ref{Fig::LocProbs::LocProb22}. In addition, we combine the trace information on $\Gamma_{21}^R$ and $\Gamma_{12}^T$ to define Dirichlet and Neumann traces on $\Gamma_{22}^{BL}$ as
\begin{alignat}{1}
\label{CombinedTracesDirichlet}
\lambda_{22}^{BL}(x)&:=\left\{\begin{array}{ll}
\lambda_{21}^{R}(x) & \quad x\in\Gamma_{22}^{BL}\cap\Gamma_{22}^L,\\
\lambda_{12}^{T}(x) & \quad x\in\Gamma_{22}^{BL}\cap\Gamma_{22}^B,
\end{array},\right.\\
\label{CombinedTracesNeumann}
\mu_{22}^{BL}(x)&:=\left\{\begin{array}{ll}
\mu_{21}^{R}(x) & \quad x\in\Gamma_{22}^{BL}\cap\Gamma_{22}^L,\\
\mu_{12}^{T}(x) & \quad x\in\Gamma_{22}^{BL}\cap\Gamma_{22}^B.
\end{array}\right.
\end{alignat}
The traces can be used to compute a polarized wavefield given by
\begin{alignat}{1}
\label{LShapedPol}
u_{22}(x)=\int_{\Gamma_{22}^{BL}}\mu_{22}^{BL}(y)G_{22}(x,y) ds_y-\int_{\Gamma_{22}^{BL}}\lambda_{22}^{BL}(y)\left[n\cdot\left(\Lambda\nabla_y G_{22}(x,y)\right)\right] ds_y, 
\end{alignat}
where $G_{22}(x,y)$ is the Green's function corresponding to the local problem defined on $\Omega^\varepsilon_{22}$. In particular,~\eqref{LShapedPol} holds due to the fact that the local problems defined on $\Omega_{12}^\varepsilon$ and $\Omega_{21}^\varepsilon$ coincide with the local problem in $\Omega_{22}^\varepsilon$ in the vicinity of $\Gamma_{22}^B$ and $\Gamma_{22}^L$, respectively.

\begin{figure}[htp]
\centering
\includegraphics[scale=0.25]{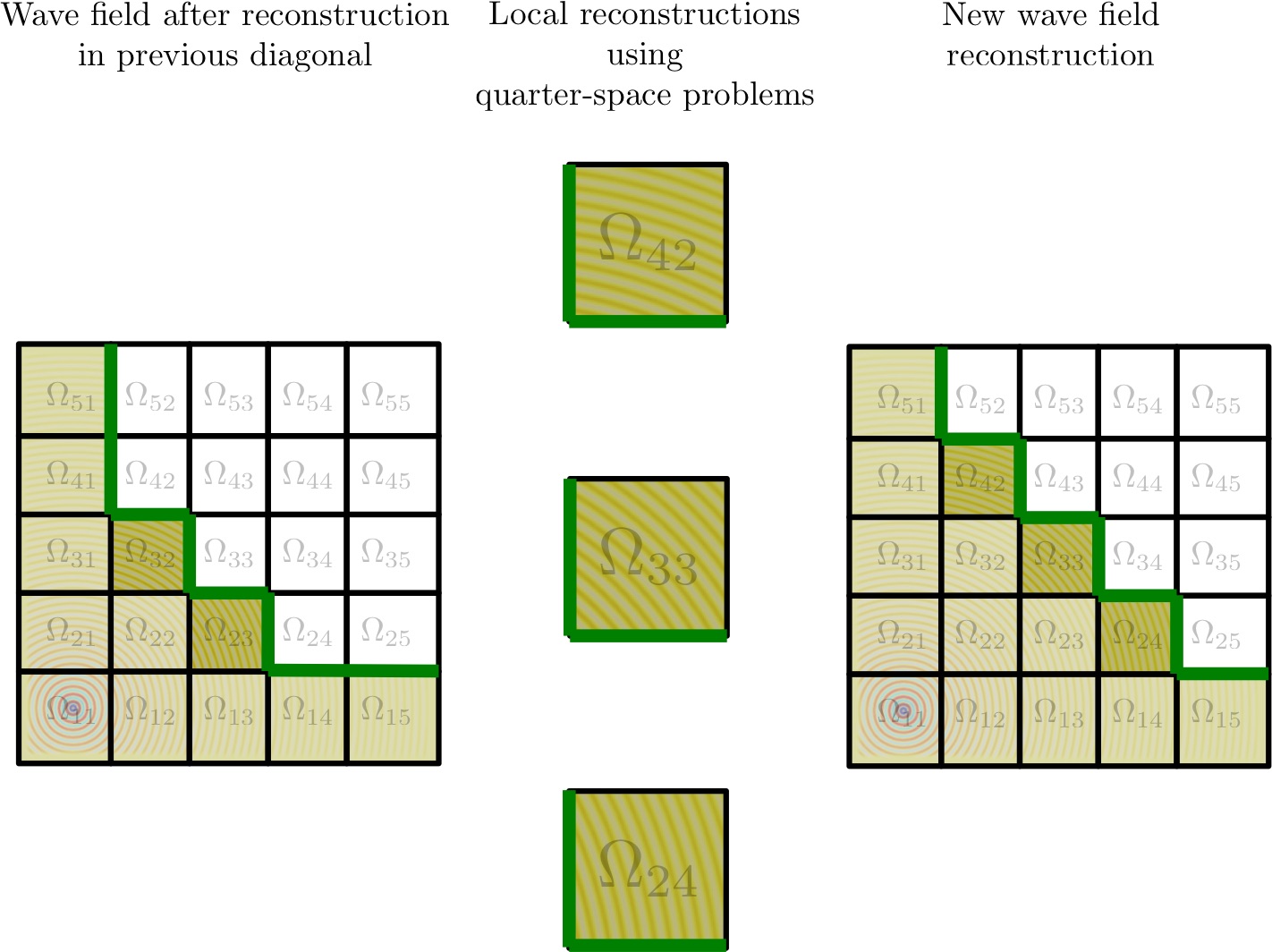}
\caption{Summary of the computed solution in one diagonal sweeping step.  Updated subdomains are diagonally adjacent in a direction perpendicular to the bottom-left to top-right sweep direction.}
\label{Fig::LSweepsSummary}
\end{figure}

It is clear that the polarized wavefield $u_{22}$ coincides with the global wavefield $u$ in $\Omega_{22}$ up to errors induced by the PMLs, and we simply set $u|_{\Omega_{22}}:=u_{22}|_{\Omega_{22}}$. Similarly to~\eqref{Traces}, we extract the Dirichlet and Neumann traces of $u_{22}$ on $\Gamma_{22}^T$ and $\Gamma_{22}^R$ and we use them to compute the approximations of the global wavefield in the other subdomains.
Following this pattern, the wavefield is computed in the remaining subdomains by sweeping diagonally from the bottom-left corner to the top-right corner. At each step of the sweep, a diagonal perpendicular to the direction of the sweep is updated. This diagonal consists of subdomains only touching each other at a corner.
For one diagonal, this procedure is summarized in Figure~\ref{Fig::LSweepsSummary}.

\subsubsection{Scenario 2: A source density supported in an arbitrary subdomain}
\label{Sec::ArbSrcPoint}
\begin{figure}[htp]
\centering
\includegraphics[scale=0.4]{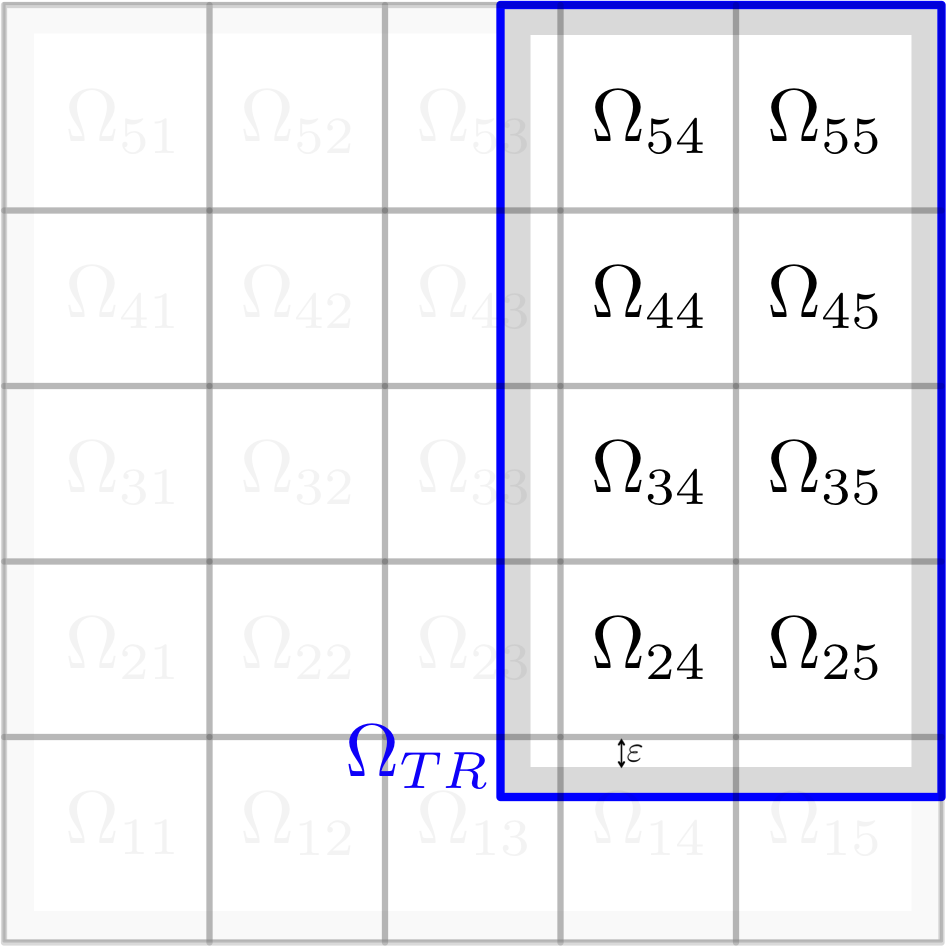}
\caption{The definition of the subproblem by truncation with absorbing boundary conditions for scenario 2.}
\label{Fig::DefSubproblem}
\end{figure}

\begin{figure}[htp]
\centering
\includegraphics[scale=0.125]{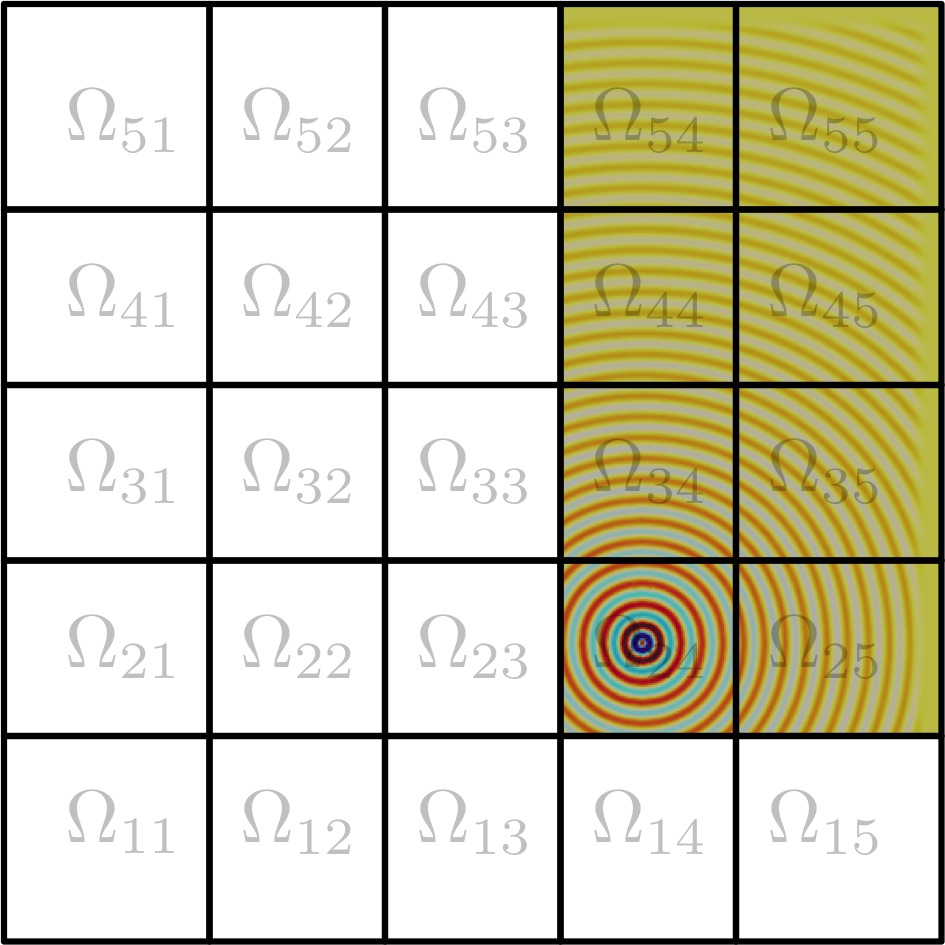}\hspace{1cm}\includegraphics[scale=0.125]{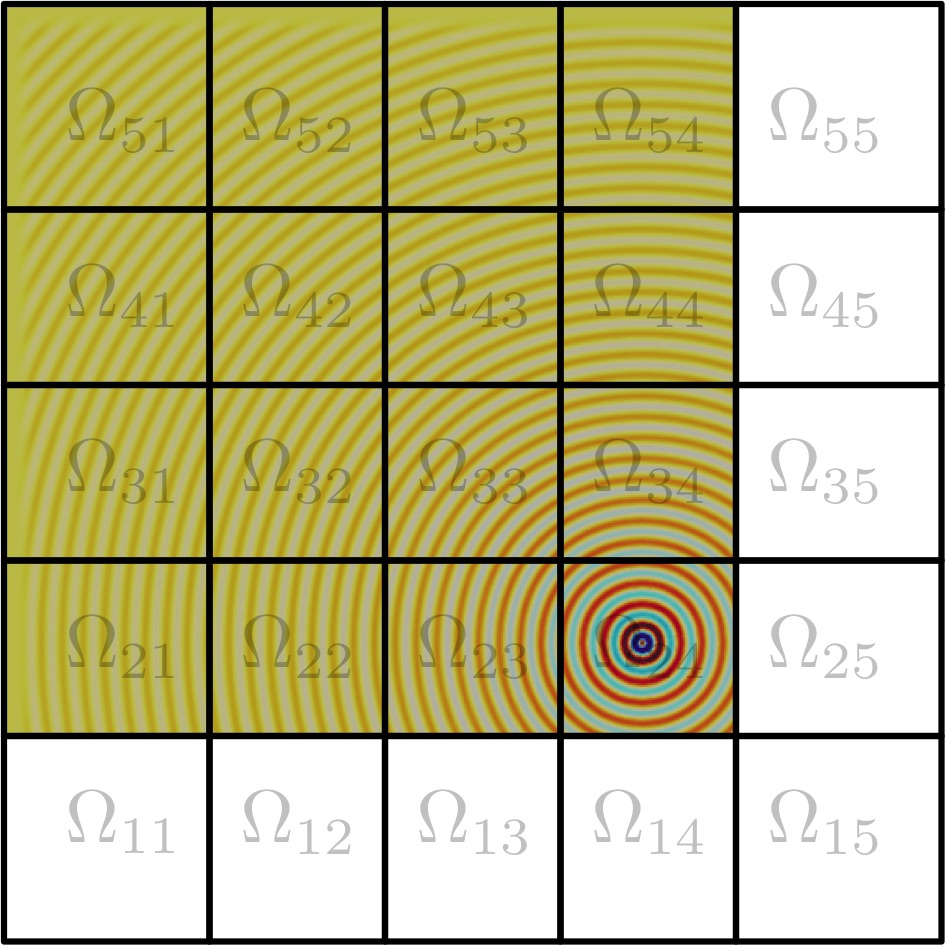}\\
\vspace{1cm}
\includegraphics[scale=0.125]{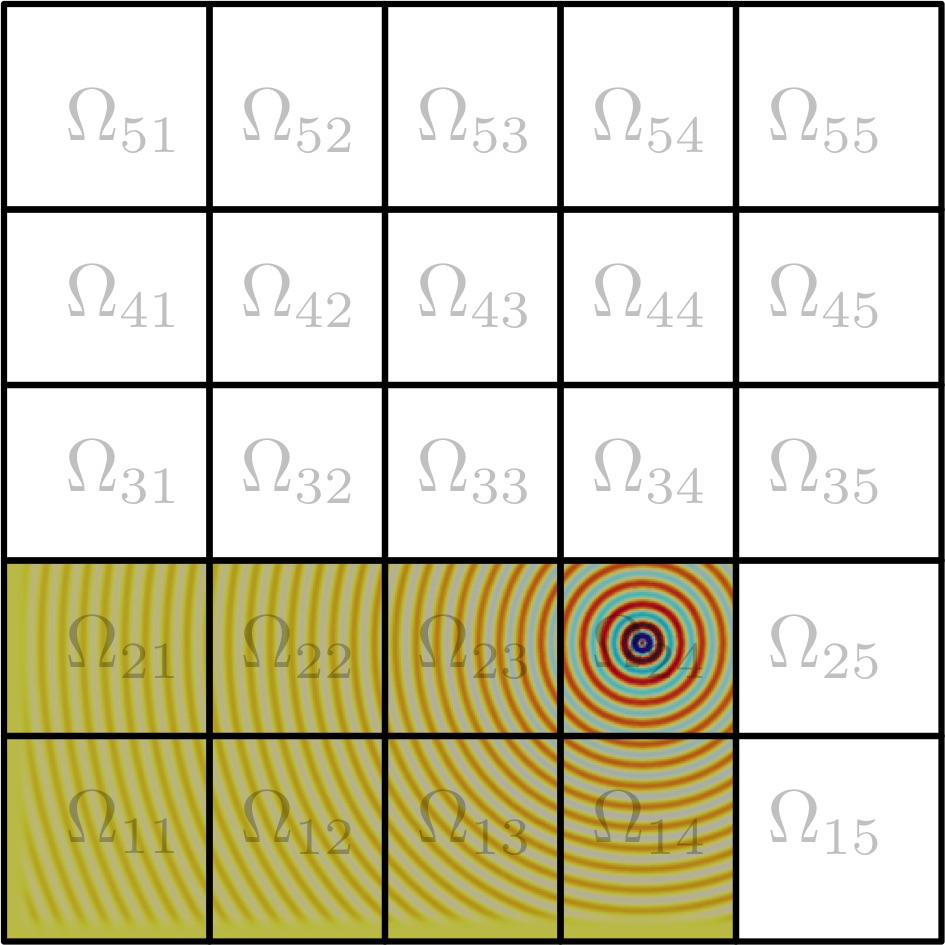}\hspace{1cm}\includegraphics[scale=0.125]{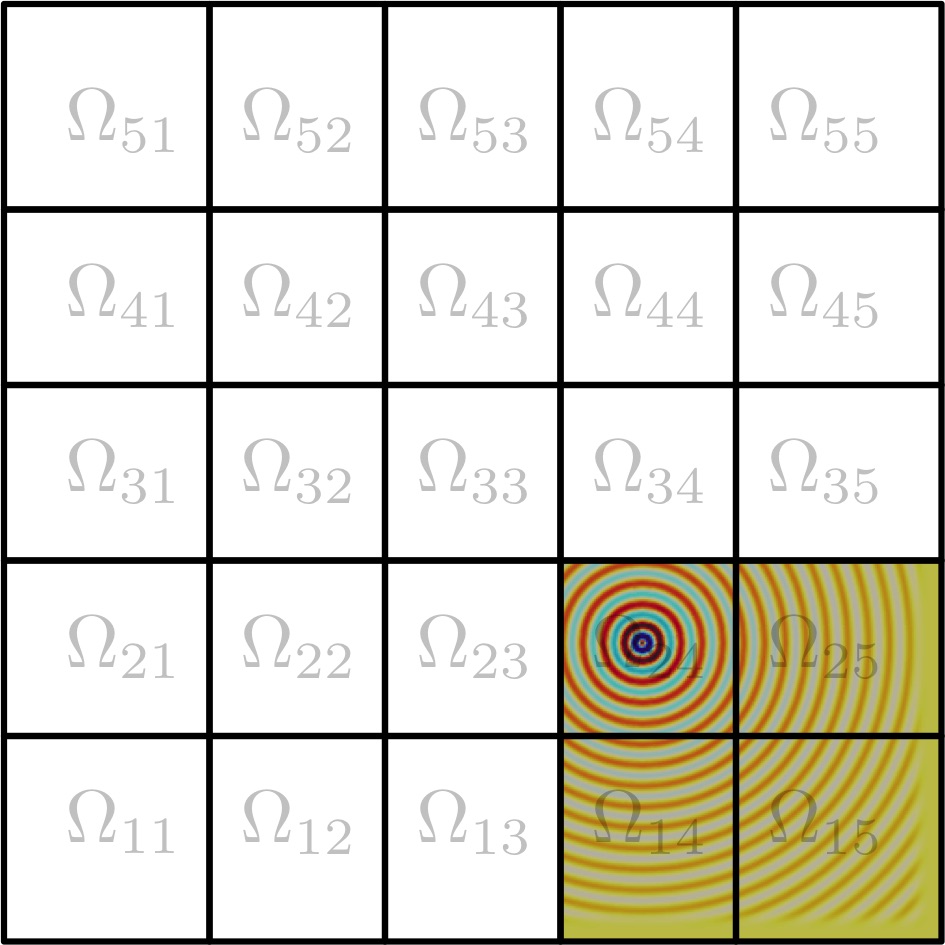}
\caption{The solutions in each quadrant in scenario 2.}
\label{Fig::ArbSrc}
\end{figure}

In this scenario, we consider a point source supported in an arbitrary subdomain.  Without loss of generality, for illustrative purposes we select the subdomain $\Omega_{24}$.  Source densities supported in any other subdomain can be treated analogously.  First we restrict the problem to the quadrant in the top right corner of the domain, inclusive of $\Omega_{24}$, which we denote $\Omega_{TR}$ and illustrate in Figure~\ref{Fig::DefSubproblem}, by truncating the global problem with PMLs.  The interior boundary of $\Omega_{TR}$ is extended by a layer of thickness $\varepsilon$ so that the local problems in $\Omega_{TR}$ are precisely the same as in Section~\ref{Sec::LSweeps::Cont::CornerSrc}. Similarly, we define domains $\Omega_{TL}$, $\Omega_{BL}$, and $\Omega_{BR}$ in the top-left, bottom-left, and bottom-right corner.
These definitions reduce the problem posed on $\Omega$ to four sub-problems, each of which have a point source supported in a corner subdomain. Thus, we can readily apply the procedure introduced in Section~\ref{Sec::LSweeps::Cont::CornerSrc} to each of them. The computed solutions are shown in Figure~\ref{Fig::ArbSrc}.

The global solution may be obtained from summing these four solutions, taking care to avoid counting the contributions from the overlapped region multiple times. 
Executing this procedure as described requires many redundant computations.  We achieve an equivalent method by generalizing the algorithm presented in  Section~\ref{Sec::LSweeps::Cont::CornerSrc}, which we describe in three stages:
\begin{enumerate}
\item compute the local solution in the subdomain where the source density is supported and extract the traces on all interior boundaries of the domain,
\item extend the solution into the subdomains in the same row/column following stage 2 of Section~\ref{Sec::LSweeps::Cont::CornerSrc}, and extract the traces required for stage 3, and,
\item extend the solution into the rest of the subdomains by computing the local solutions in diagonal rows, perpendicular to the direction of the sweep, similarly to stage 3 of Section~\ref{Sec::LSweeps::Cont::CornerSrc}. 
\end{enumerate}
\begin{figure}[htp]
\centering
\begin{subfigure}[b]{0.45\textwidth}
\centering
\includegraphics[scale=0.15]{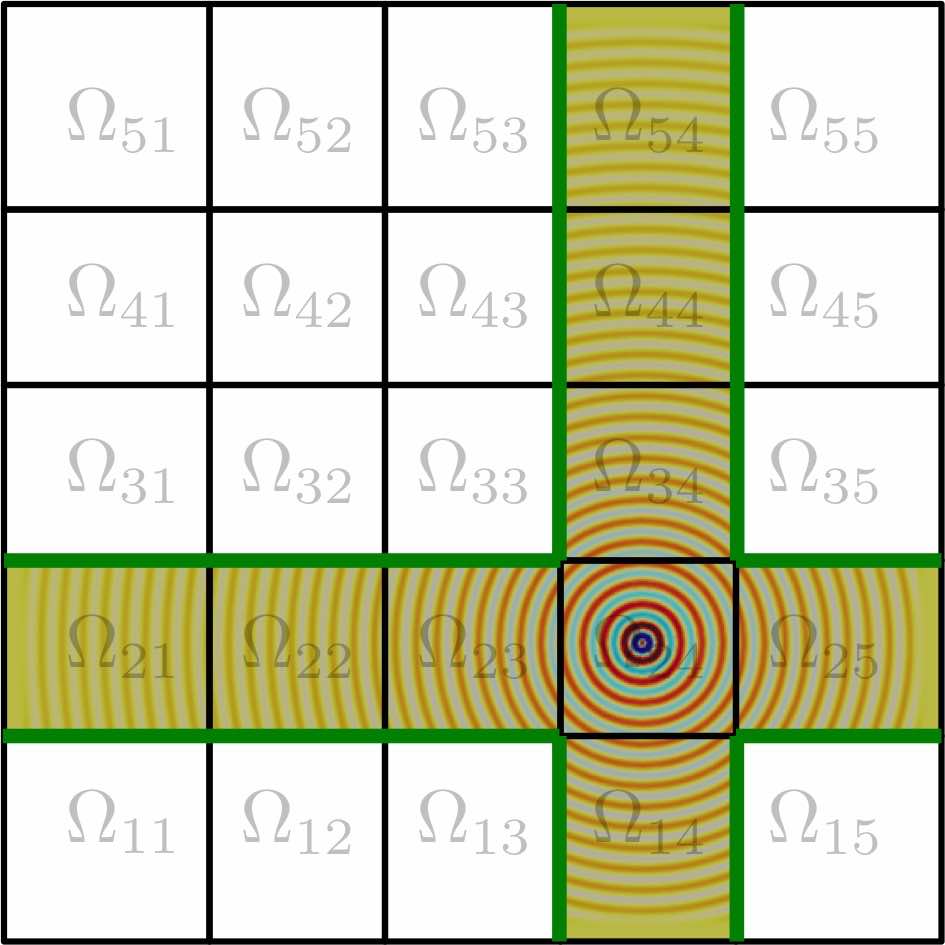}
\caption{Before stage 3.}
\end{subfigure}\\
\vspace{1cm}
\begin{subfigure}[b]{0.45\textwidth}
\centering
\includegraphics[scale=0.15]{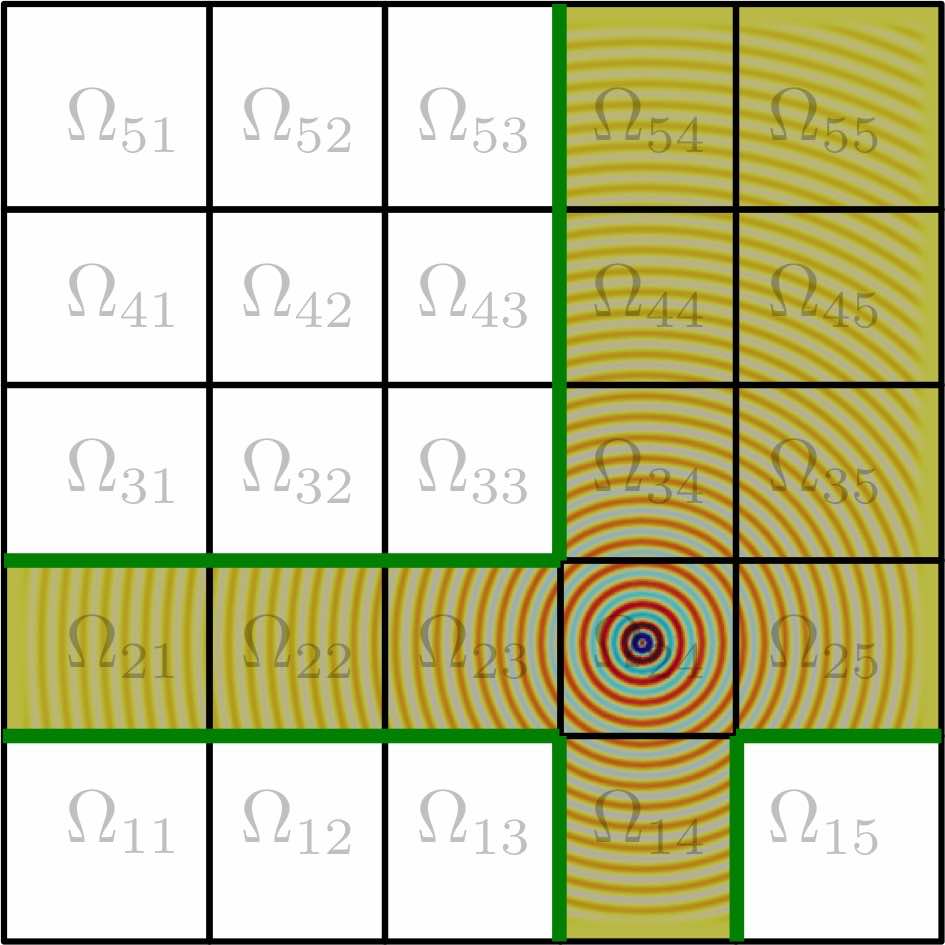}
\caption{After bottom-left to top-right sweep.}
\end{subfigure}
\begin{subfigure}[b]{0.45\textwidth}
\centering
\includegraphics[scale=0.15]{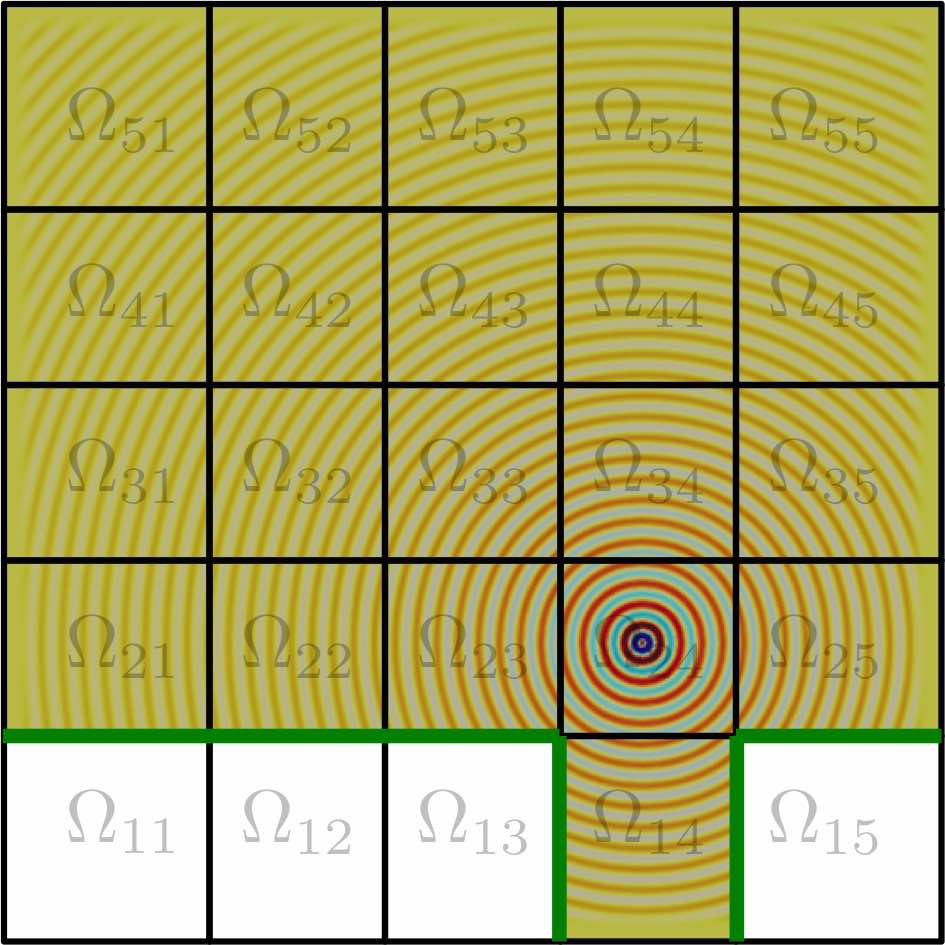}
\caption{After bottom-right to top-left sweep.}
\end{subfigure}\\
\vspace{1cm}
\begin{subfigure}[b]{0.45\textwidth}
\centering
\includegraphics[scale=0.15]{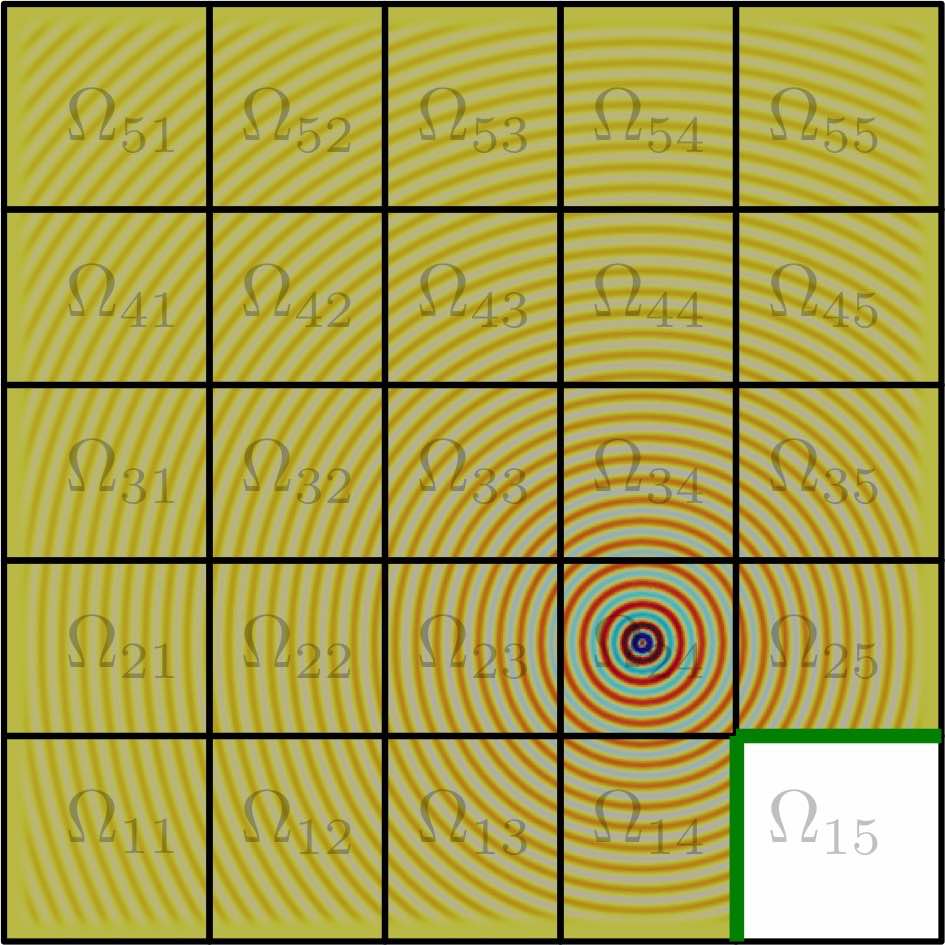}
\caption{After top-right to bottom-left sweep.}
\end{subfigure}
\begin{subfigure}[b]{0.45\textwidth}
\centering
\includegraphics[scale=0.15]{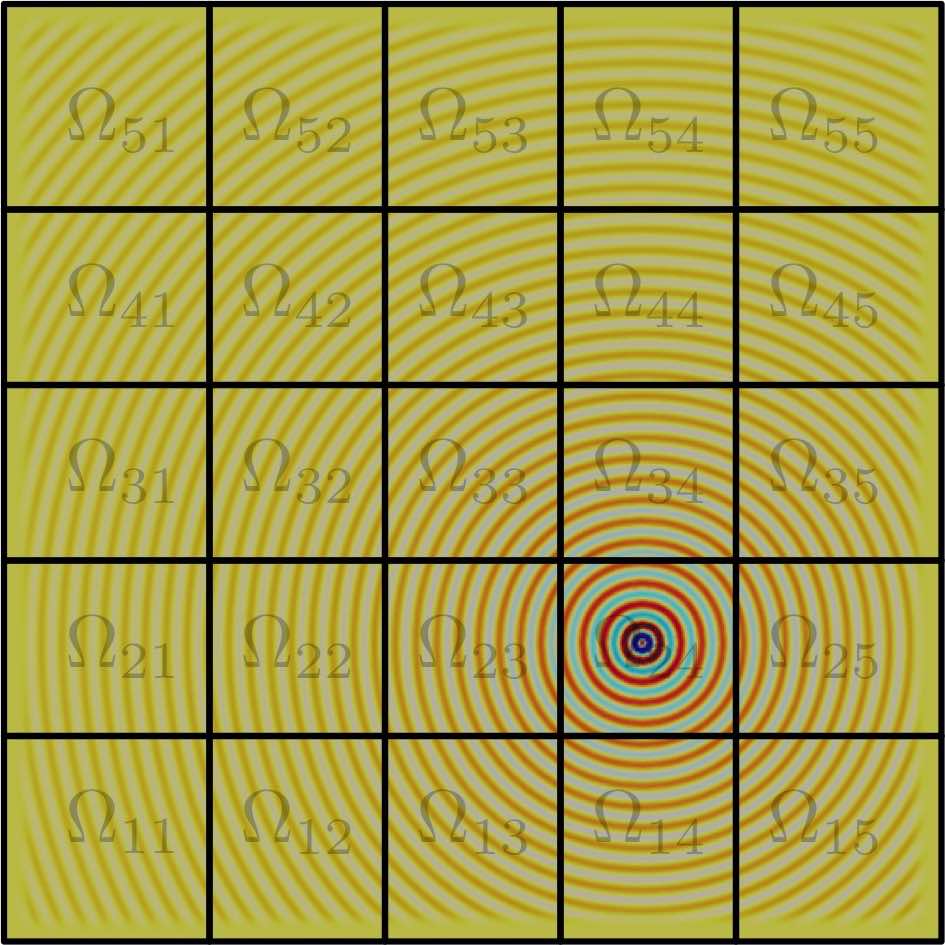}
\caption{After top-left to bottom-right sweep.}
\end{subfigure}
\caption{Summary of stage 3 for scenario 2. The plots show the computed wave field after stage 2, and the computed wave field after each sweep from corner to corner.}
\label{Fig::Cross}
\end{figure}
Stage 3 is illustrated in Figure~\ref{Fig::Cross}.  In this scenario, the full algorithm requires a total of eight sweeps. Stage 2 requires four sweeps over the domain: up, down, left, right. Stage 3 also requires four sweeps: bottom-left to top-right, top-right to bottom-left, bottom-right to top-left, top-left to bottom-right.

\subsubsection{Scenario 3: Arbitrary source distributions not intersecting the CDD skeleton}
\label{Sec::ArbSrcSkeleton}
Global solutions for arbitrary source distributions that do not intersect the skeleton of the CDD can be computed from a union of distinct source densities, each supported in a single subdomain. Therefore, the global solution can be na\"{\i}vely computed by applying the procedure from scenario 2 to each of the localized source densities and summing the results.  This approach is not computationally efficient because much of the work is redundant. 

In this section, we show how to compute the global solution without redundancy, by applying each of the eight sweeps of scenario 2 only once. This allows for an efficient computation of the global wavefield, regardless of the number of subdomains containing components of the source. The algorithm can still be performed in three stages, which are detailed below. For a succinct summary in pseudo-code we direct the reader to Appendix~\ref{App::Code}.

\begin{figure}[htp]
\centering
\includegraphics[scale=0.125]{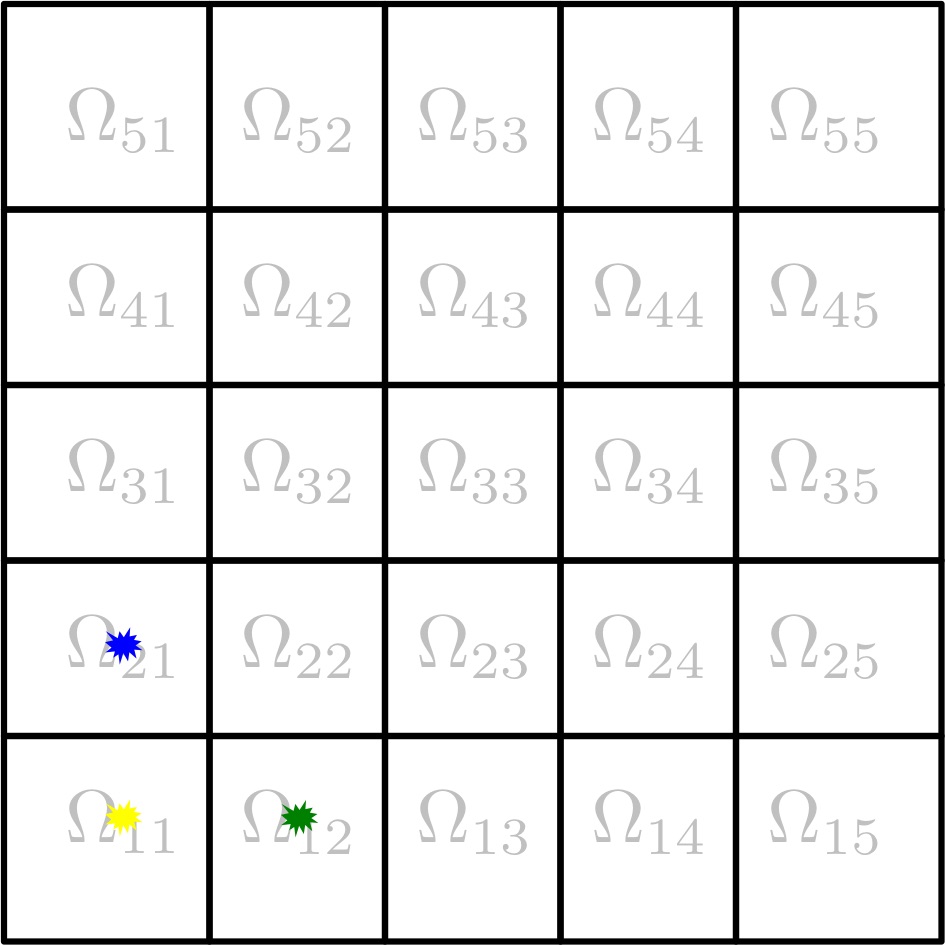}\quad\quad\quad\includegraphics[scale=0.125]{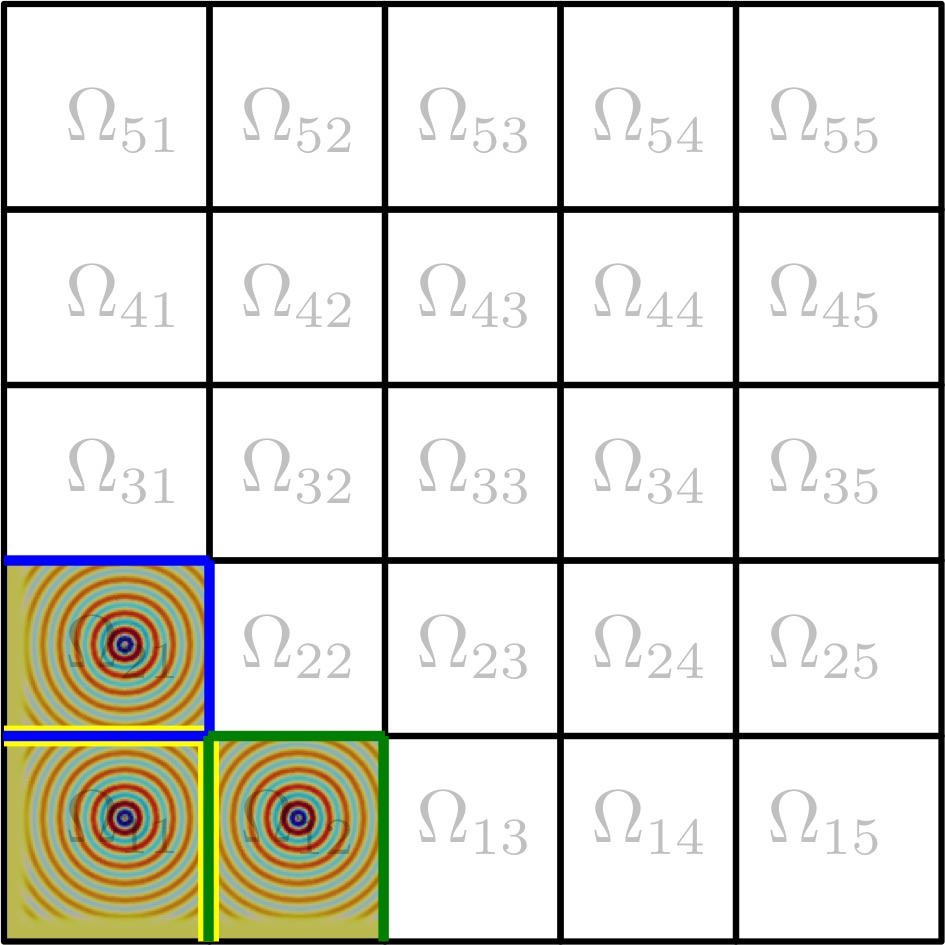}
\caption{Left: The setup of the problem in scenario 3, with multiple point sources shown in red, yellow and blue. Right: The corresponding local solutions in each subdomain. The extracted trace information is shown using the same colors as their corresponding point sources.}\label{Fig::SetupLocSols}
\end{figure}

In the first stage, we restrict the source density to each subdomain. Within each subdomain, we use this restricted source density to compute the corresponding local solution, from which we extract the traces on all interior interfaces $\Gamma_{ij}^{l}$, $l=B,R,T,L$. This stage is illustrated in Figure~\ref{Fig::SetupLocSols}.

\begin{figure}[htp]
\centering
\includegraphics[scale=0.25]{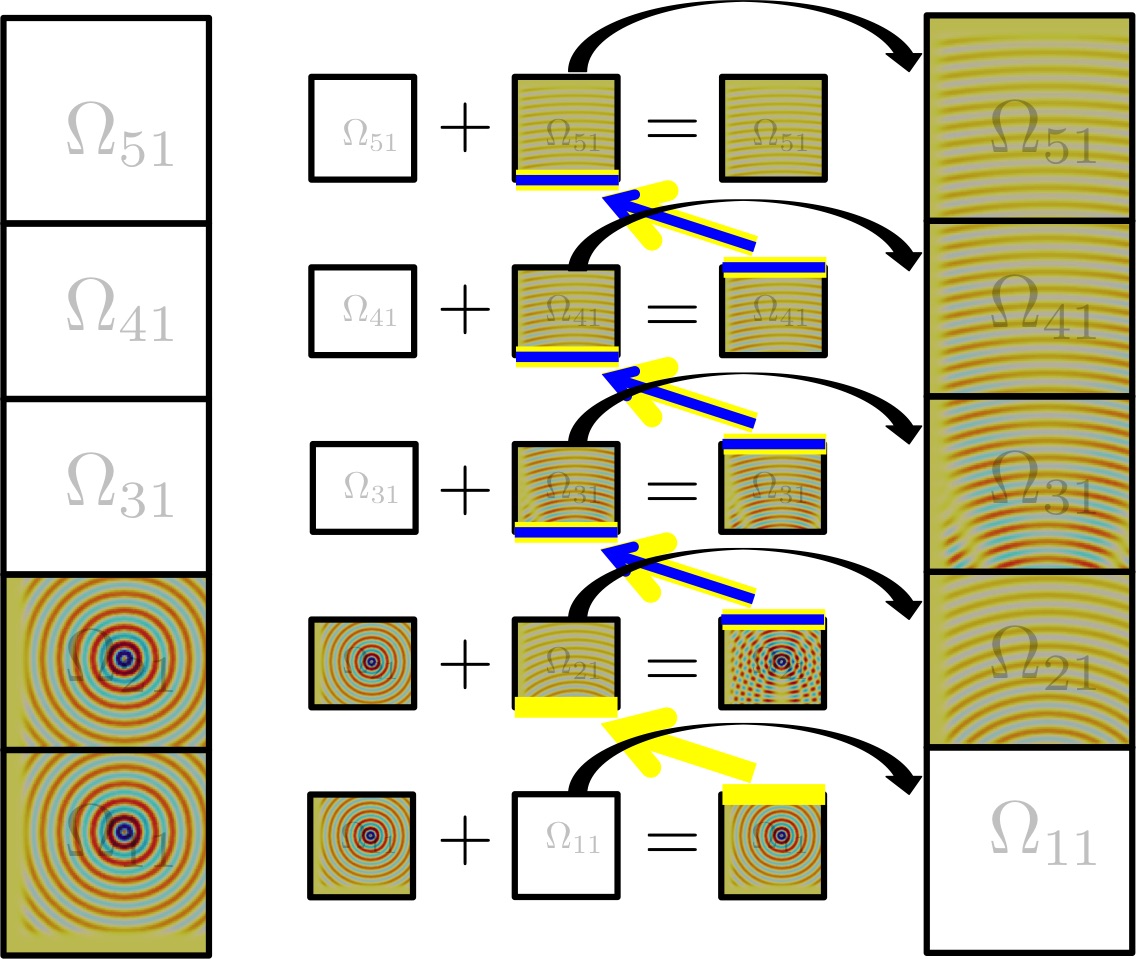}
\caption{Summary of the upwards sweep, for scenario 3, in one column in the presence of several point sources.}
\label{Fig::UpSweepSummary}
\end{figure}

In the second stage, we use the traces extracted in stage 1 to extend the local solutions into the same column and row as subdomains containing sources.  This can be done using four sweeps: up, down, left, and right. We consider the upwards sweep in the first column in detail, all other sweeps are performed analogously. In the upward sweep, we use traces on the bottom of each element to compute polarized wavefields and update local solutions. For example, in Figure~\ref{Fig::SetupLocSols}, the subdomain $\Omega_{11}$ has no incoming bottom trace. The outgoing top trace is therefore simply taken from stage one and transferred to $\Omega_{21}$. In $\Omega_{21}$, this trace is the incoming bottom trace and used to compute a polarized wavefield. We then update the local wavefield in $\Omega_{21}$ by adding the computed polarized wavefield. The outgoing trace on the top is then extracted from the updated local wavefield and used as the incoming trace in $\Omega_{31}$. Continuing this procedure, we update the local wavefields in the entire first column which concludes the upwards sweep. The procedure is illustrated in Figure~\ref{Fig::UpSweepSummary}.  Then, we perform the remaining three cardinal sweeps.

\begin{figure}[htp]
\centering
\includegraphics[scale=0.25]{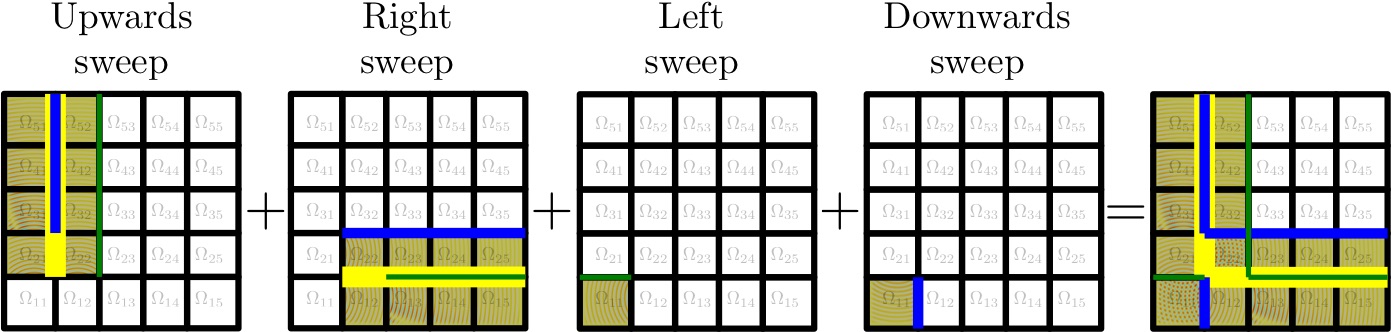}
\caption{The contributions to the intermediate wavefield from the sweeps along each row and column in scenario 3. The generated trace information that still needs to be propagated is shown by the colored lines.}
\label{Fig::StraightSweepContribution}
\end{figure}

\begin{figure}[htp]
\centering
\includegraphics[scale=0.25]{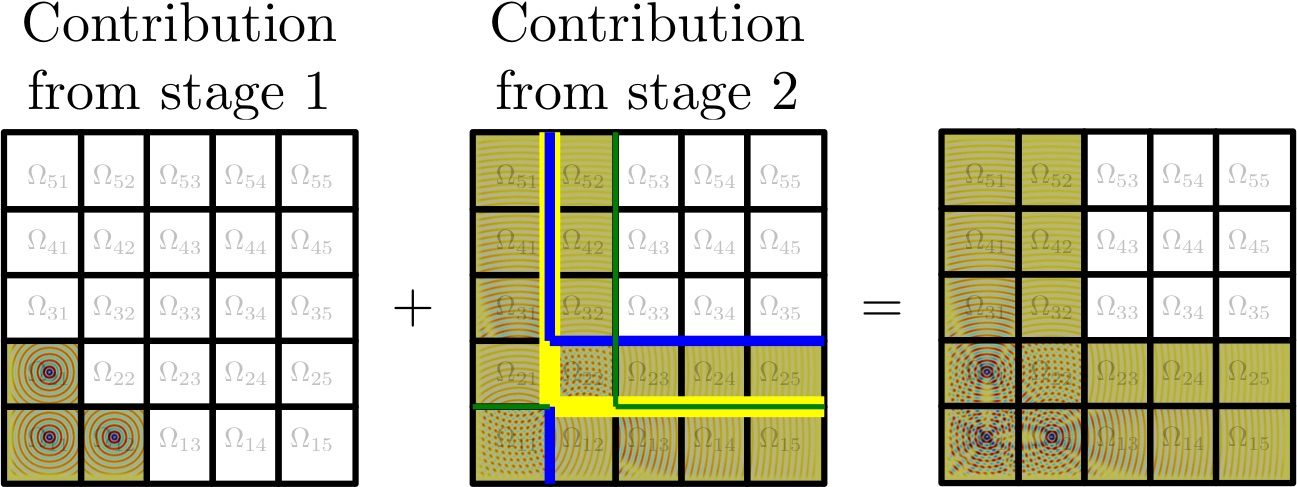}
\caption{The reconstructed wavefield composed of the contributions from stage 1 and 2 of scenario 3, for the problem in Figure~\ref{Fig::SetupLocSols}.}
\label{Fig::StraightSweepsResult}
\end{figure}

Applying these four sweeps computes an intermediate wavefield. This wavefield is updated in stage 3. As part of the vertical and horizontal sweeps we also extract the vertical and horizontal traces of the subdomains needed in stage 3. The individual contributions from each of the four sweeps are shown in Figure~\ref{Fig::StraightSweepContribution}. The intermediate solution obtained from the sum of the contributions from stages 1 and 2 is illustrated in Figure~\ref{Fig::StraightSweepsResult}. In both figures the traces extracted in preparation for stage 3 are shown by the colored lines.

We generalize the third stage of the algorithm in a similar way by computing polarized wavefields from the trace information coming into a subdomain, and updating the local solutions by adding these polarized wavefields. We then extract trace information from the updated local wavefields and use them in the neighboring subdomains as incoming traces.

\begin{figure}[htp]
\centering
\includegraphics[scale=0.2]{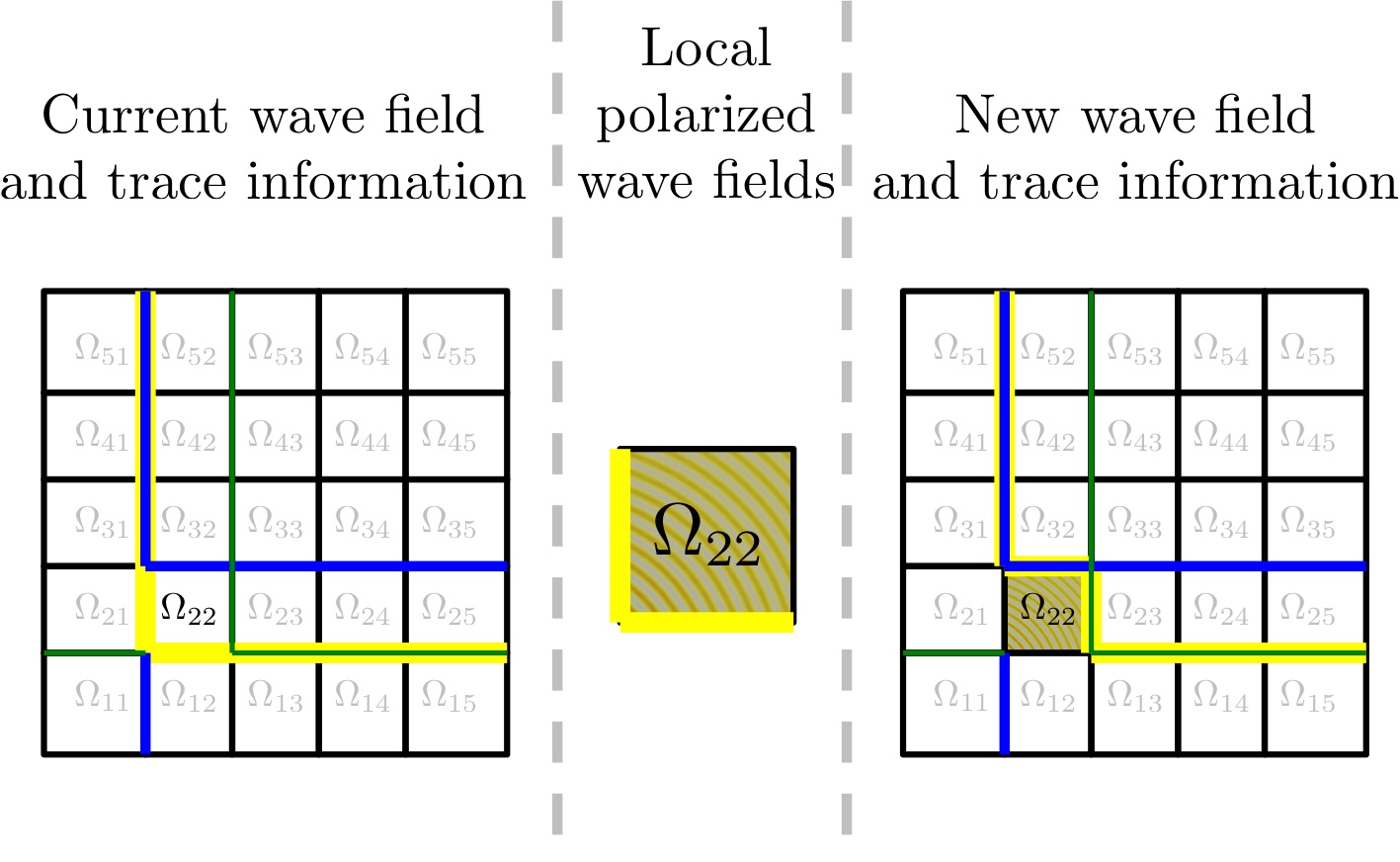}\\
\includegraphics[scale=0.2]{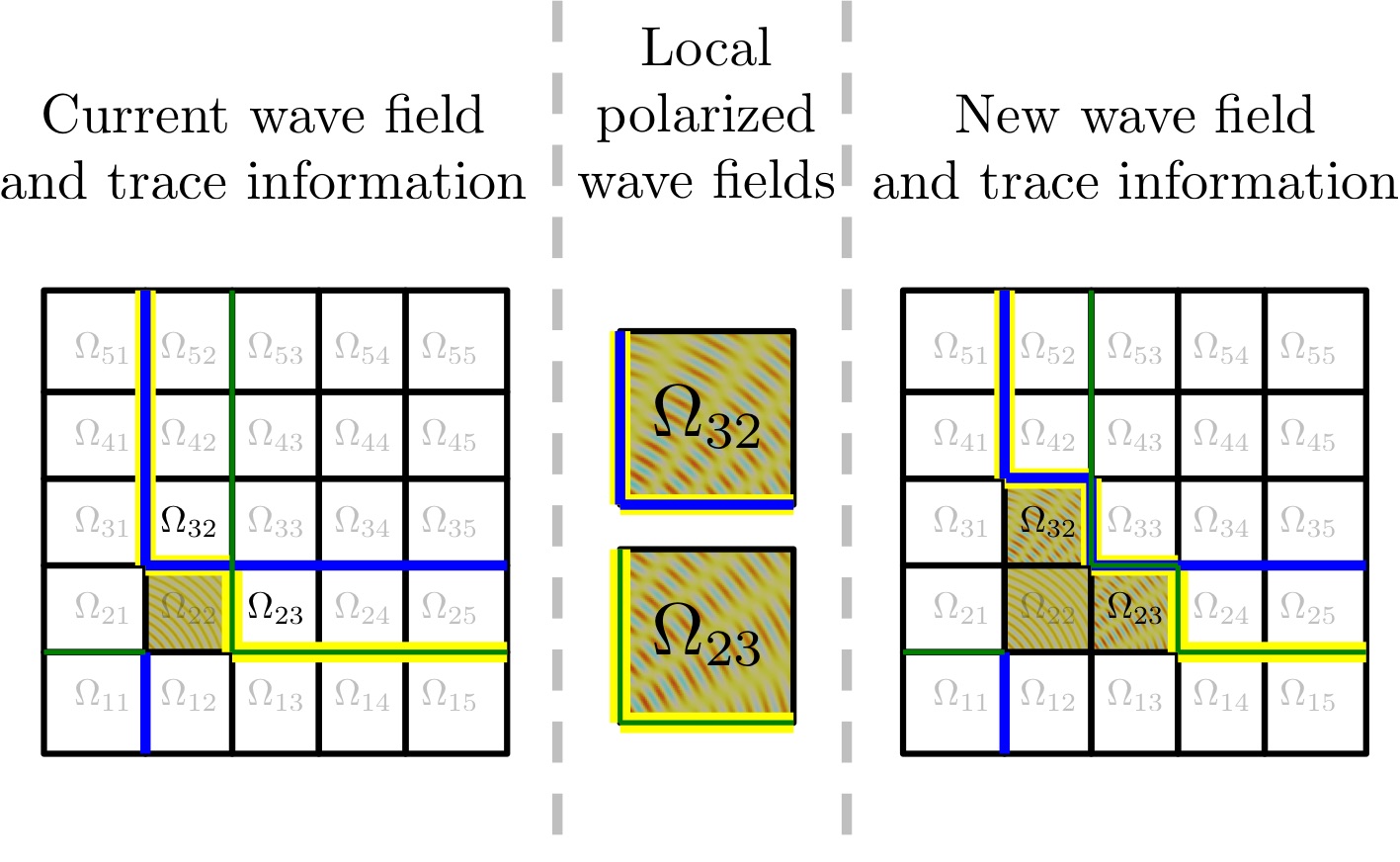}\\
\includegraphics[scale=0.2]{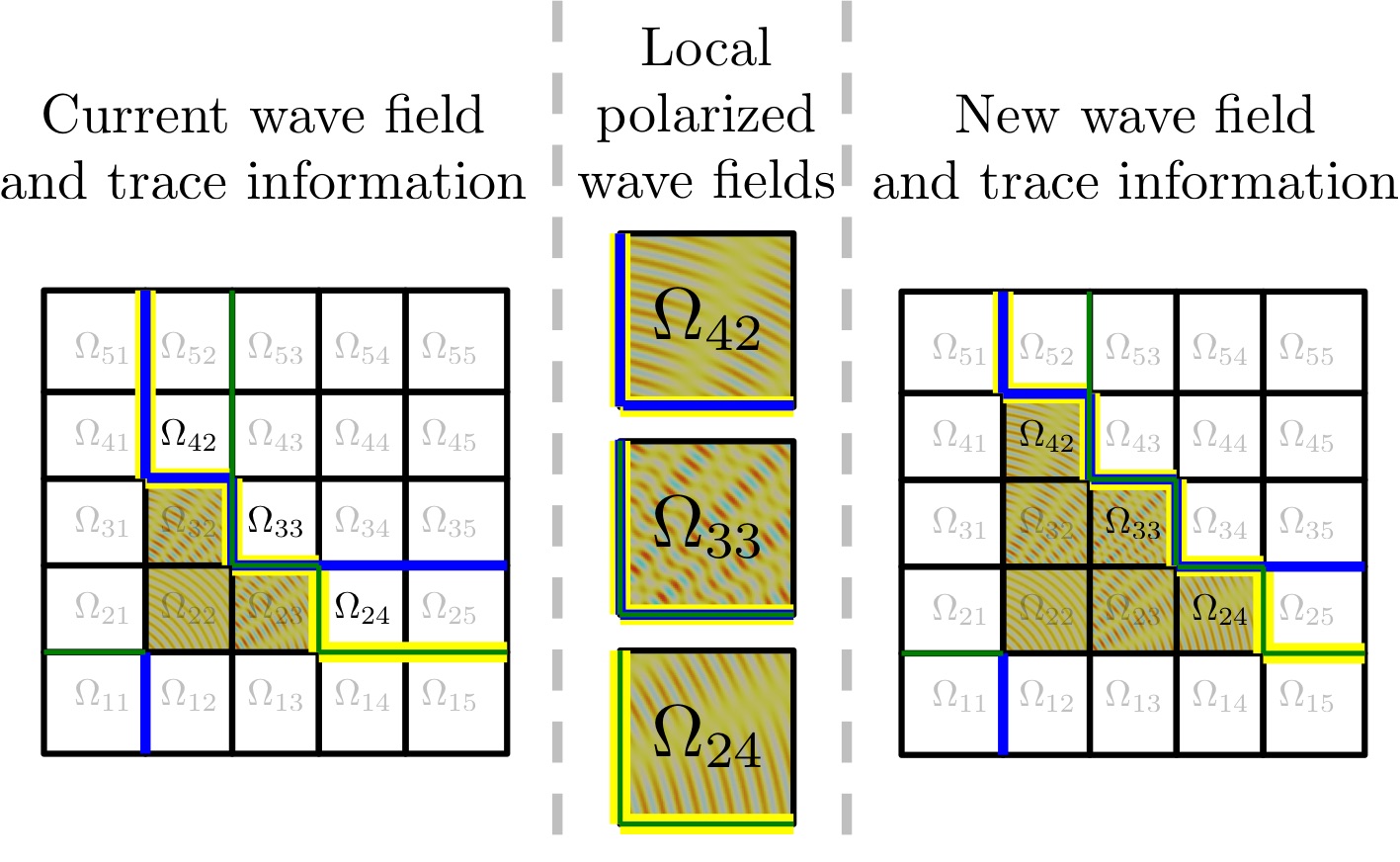}
\caption{Illustration of the sweep from the bottom-left to the top-right corner in scenario 3.}
\label{Fig::DiagonalSweep}
\end{figure}

\begin{figure}[htp]
\begin{subfigure}[b]{\textwidth}
\centering
\includegraphics[scale=0.25]{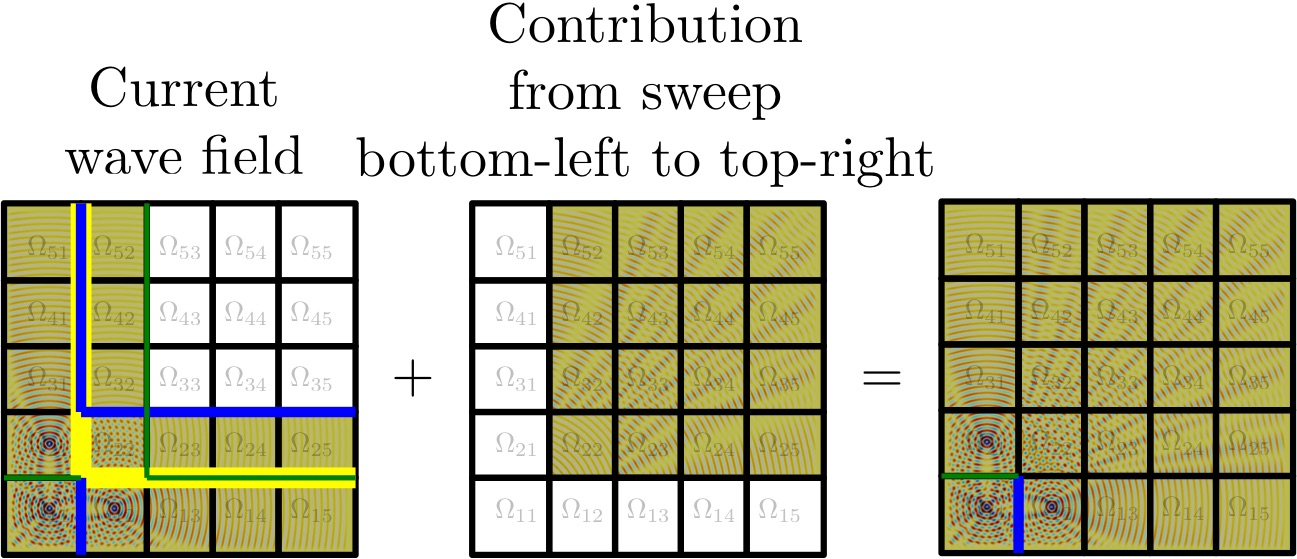}
\subcaption{The sweep from the bottom-left to the top-right corner.}
\end{subfigure}\\
\vspace{1cm}\\
\begin{subfigure}[b]{\textwidth}
\centering
\includegraphics[scale=0.25]{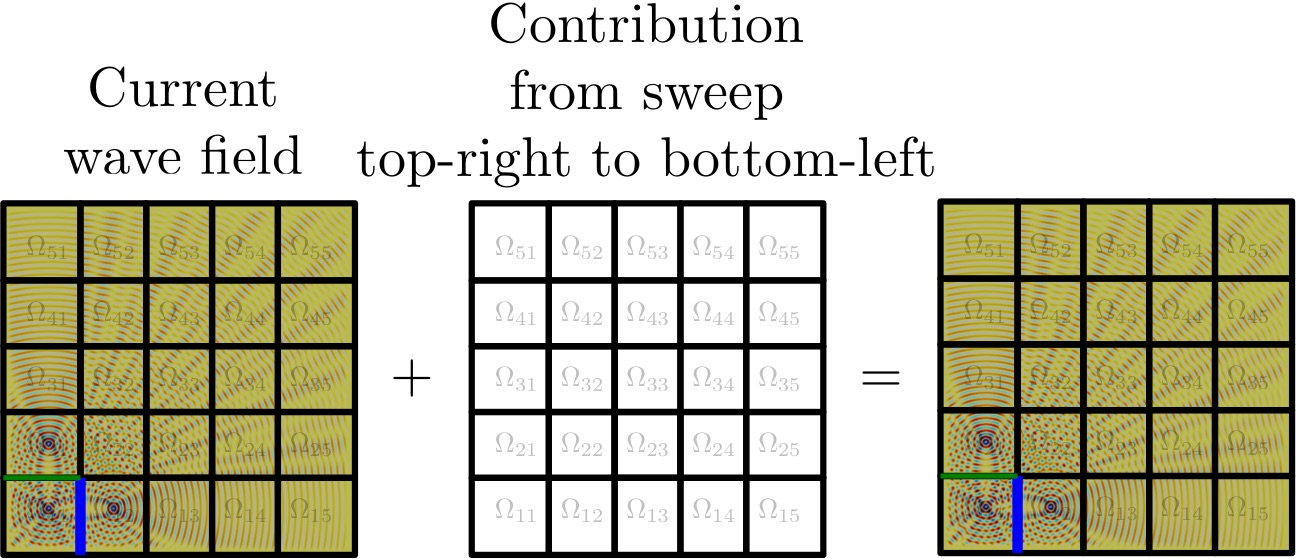}
\subcaption{The sweep from the top-right to the bottom-left corner.}
\end{subfigure}
\caption{The contribution of the sweep from the sweeps over the sweeps diagonal by diagonal in scenario 3. The extracted trace information that still needs to be propagated is shown by the colored lines.}
\label{Fig::DiagonalSweeps1}
\end{figure}

\begin{figure}
\begin{subfigure}[b]{\textwidth}
\centering
\includegraphics[scale=0.25]{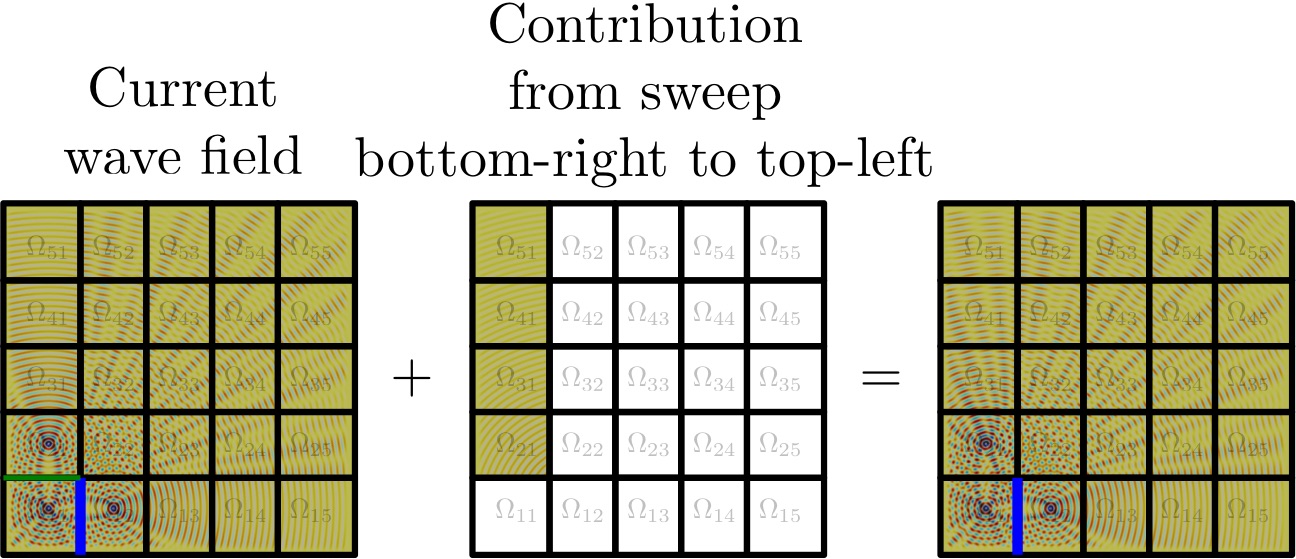}
\subcaption{The sweep from the bottom-right to the top-left corner.}
\end{subfigure}\\
\vspace{1cm}\\
\begin{subfigure}[b]{\textwidth}
\centering
\includegraphics[scale=0.25]{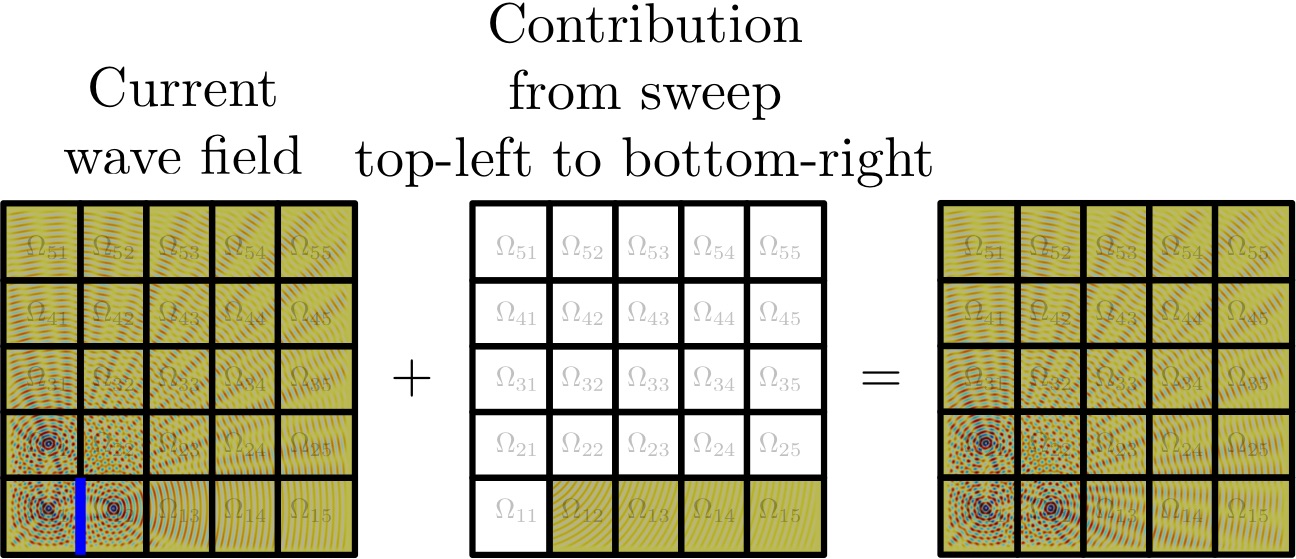}
\subcaption{The sweep from the top-left to the bottom-right corner.}
\end{subfigure}
\caption{The contribution of the sweep from the sweeps over the sweeps diagonal by diagonal in scenario 3. The extracted trace information that still needs to be propagated is shown by the colored lines.}
\label{Fig::DiagonalSweeps2}
\end{figure}

As in Section~\ref{Sec::ArbSrcPoint}, stage 3 is realized by sweeping over the CDD from corner to corner, updating subdomains that are diagonally adjacent, perpendicular to the sweep direction. Computed wavefields in three example diagonals for the sweep from the bottom-left to the top-right corner are illustrated in Figure~\ref{Fig::DiagonalSweep}.
The same technique can also be employed for the other three diagonal sweeps: top-right corner to bottom-left corner, bottom-right to top-left corner, and top-left to bottom-right corner. These sweeps are illustrated in Figures~\ref{Fig::DiagonalSweeps1} and~\ref{Fig::DiagonalSweeps2}.



Using these three stages, we compute the global solution in a total of eight sweeps, independently of the source distribution.  In the next section, we show that applying this procedure to source densities appropriately windowed on four slightly different domain decompositions, the algorithm can be extended to entirely arbitrary source densities in a straightforward way.

\subsubsection{Scenario 4: Arbitrary source distributions}
\label{Sec::ArbSrc}
\begin{figure}[htp]
\centering
\includegraphics[scale=0.18]{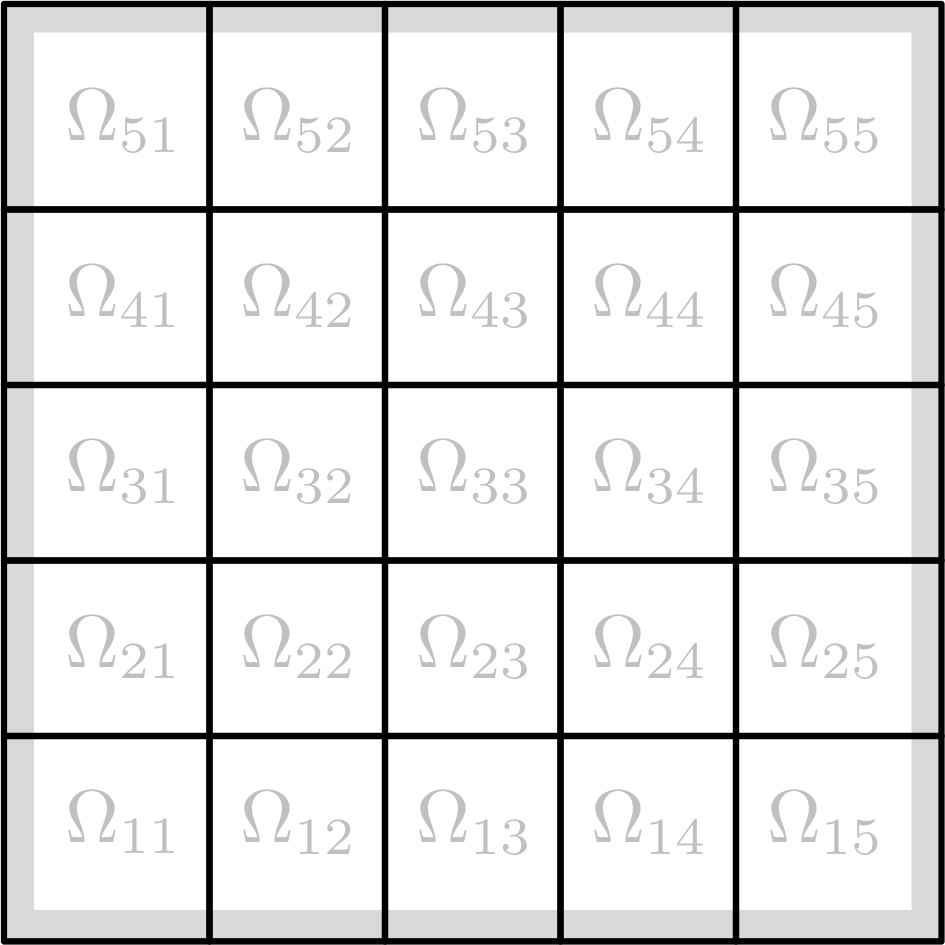}\quad\includegraphics[scale=0.18]{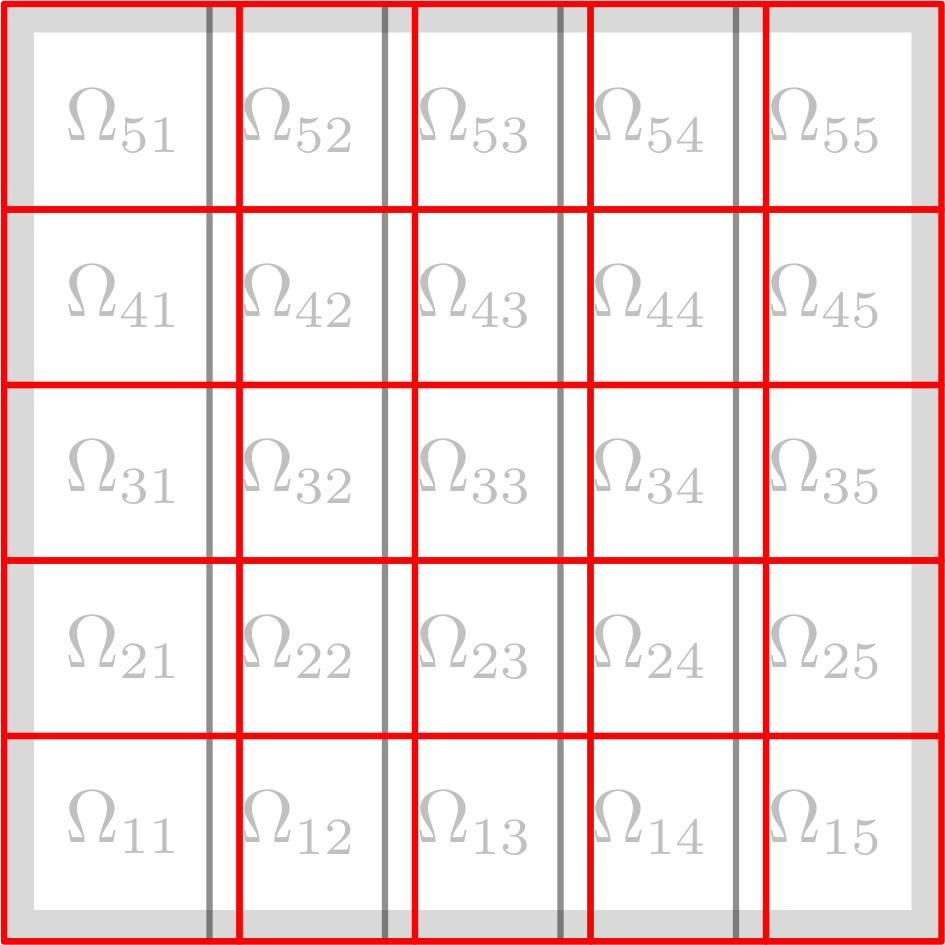}\quad\includegraphics[scale=0.18]{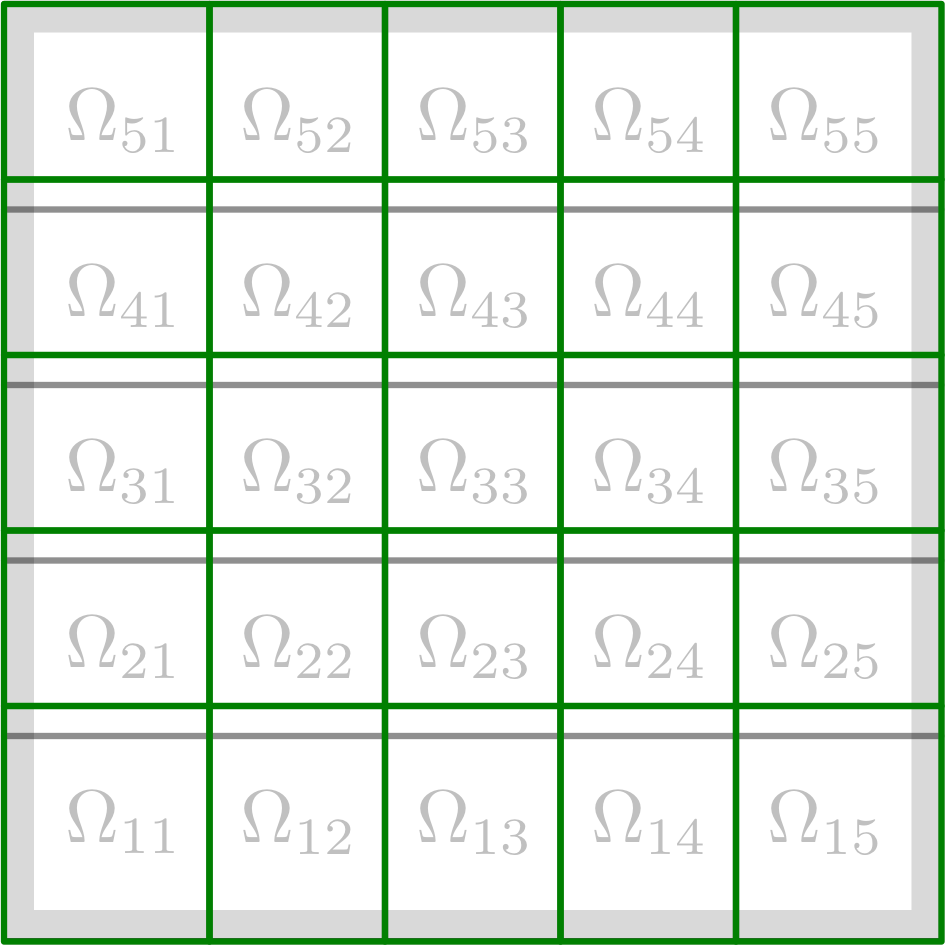}\quad\includegraphics[scale=0.18]{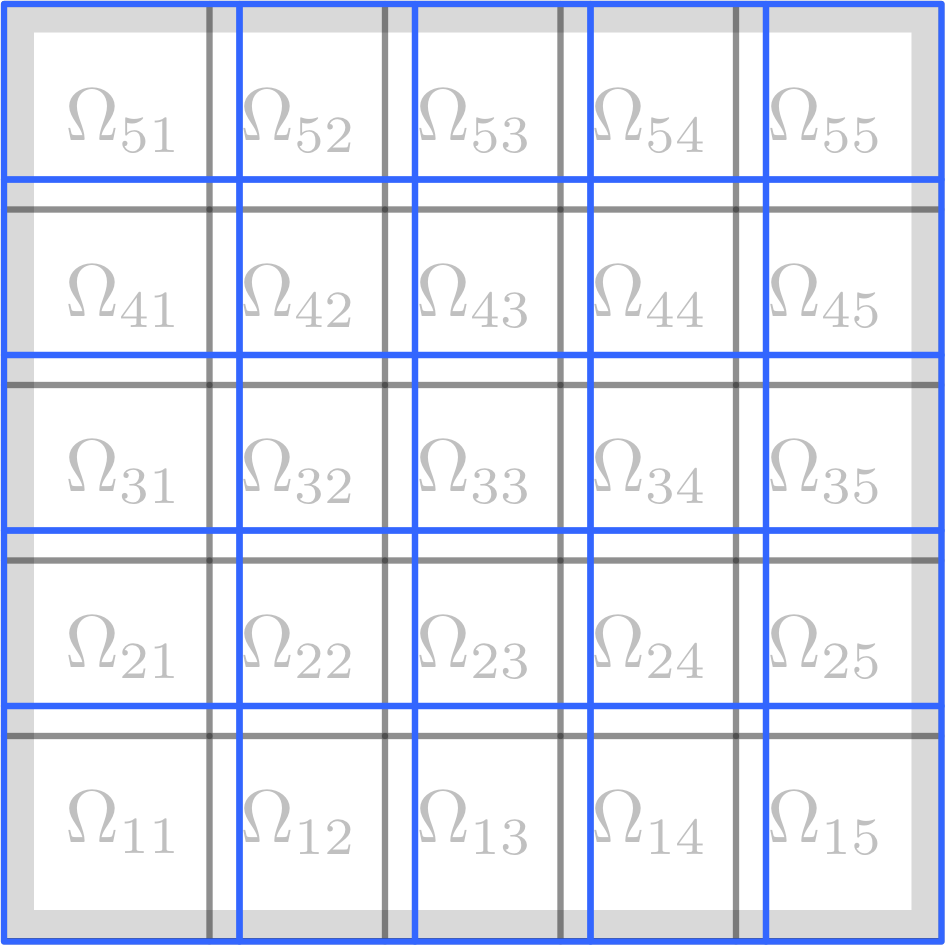}
\caption{The original CDD (black) and the three new, shifted CDDs in red, green, and blue.}
\label{Fig::NewDD}
\end{figure}

\begin{figure}[htp]
\centering
\includegraphics[scale=0.15]{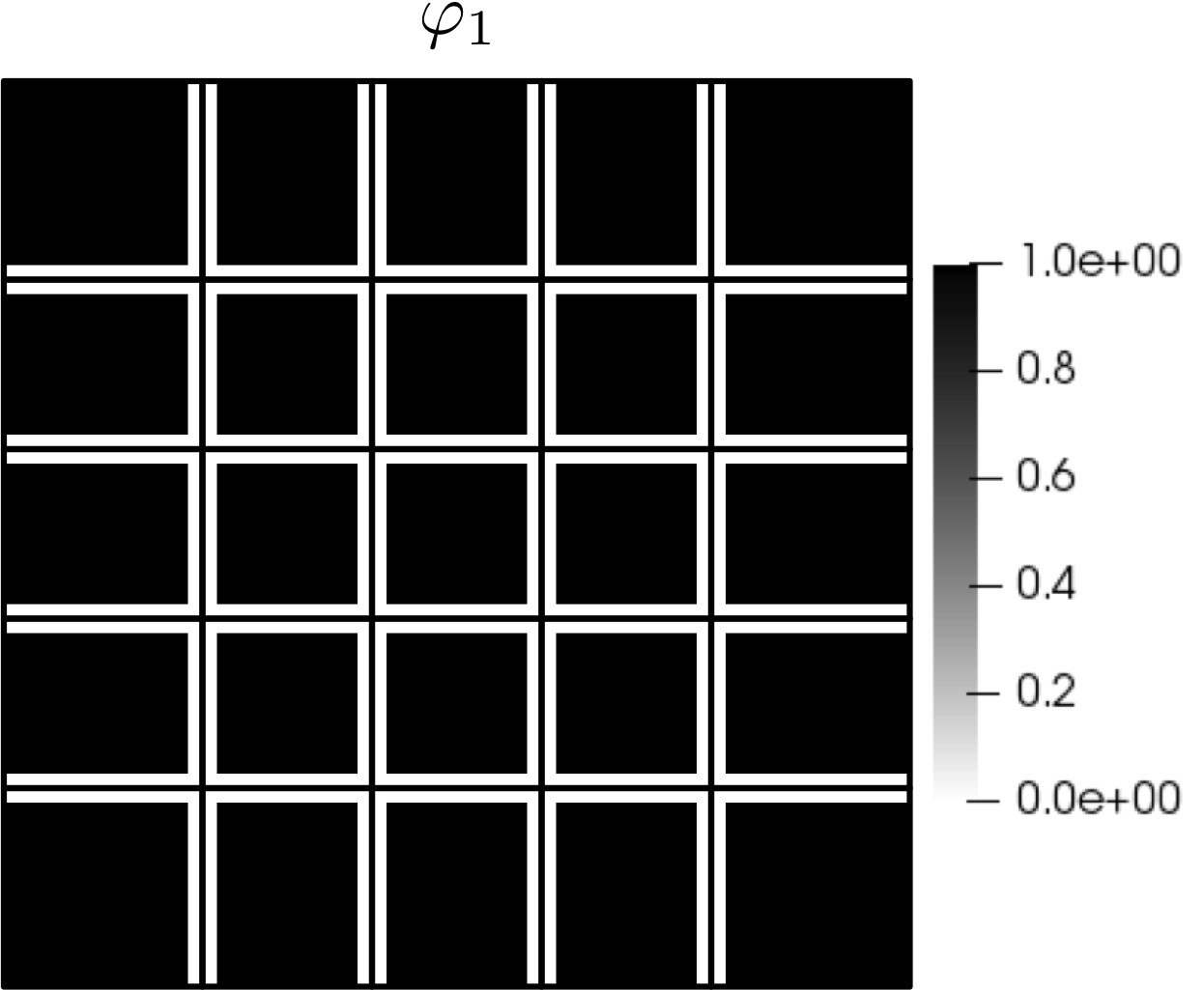}\quad\includegraphics[scale=0.15]{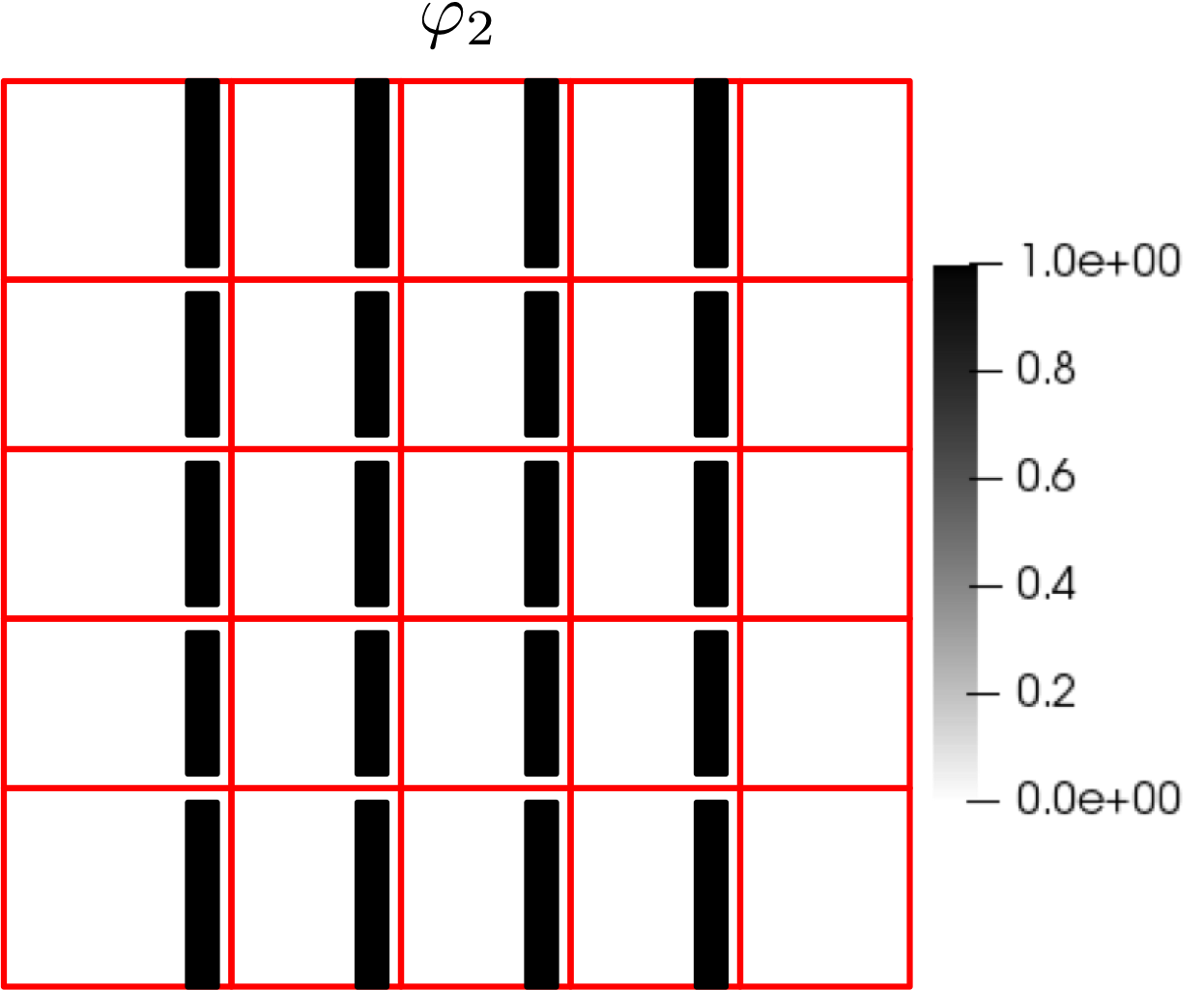}\quad\includegraphics[scale=0.15]{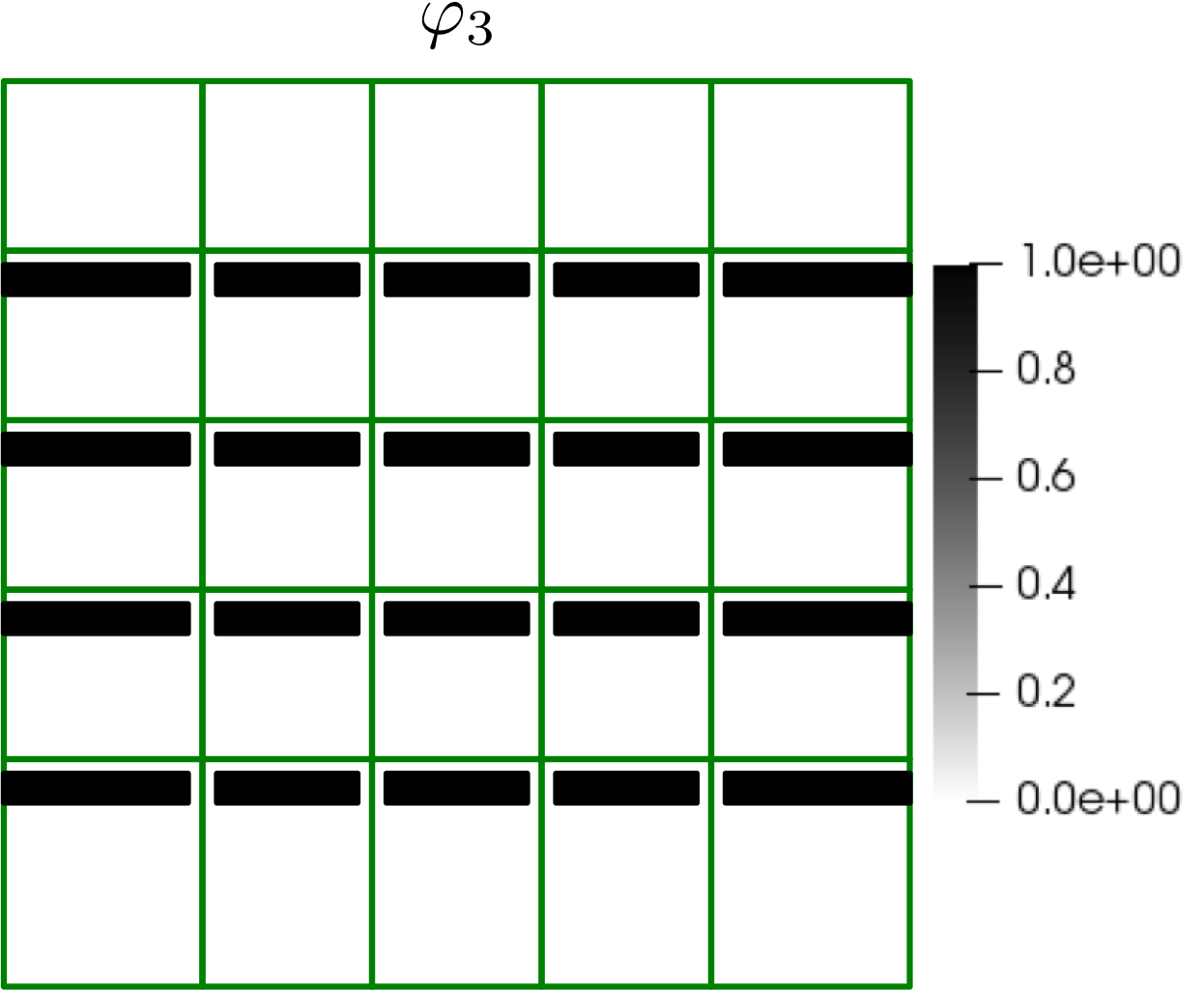}\quad\includegraphics[scale=0.15]{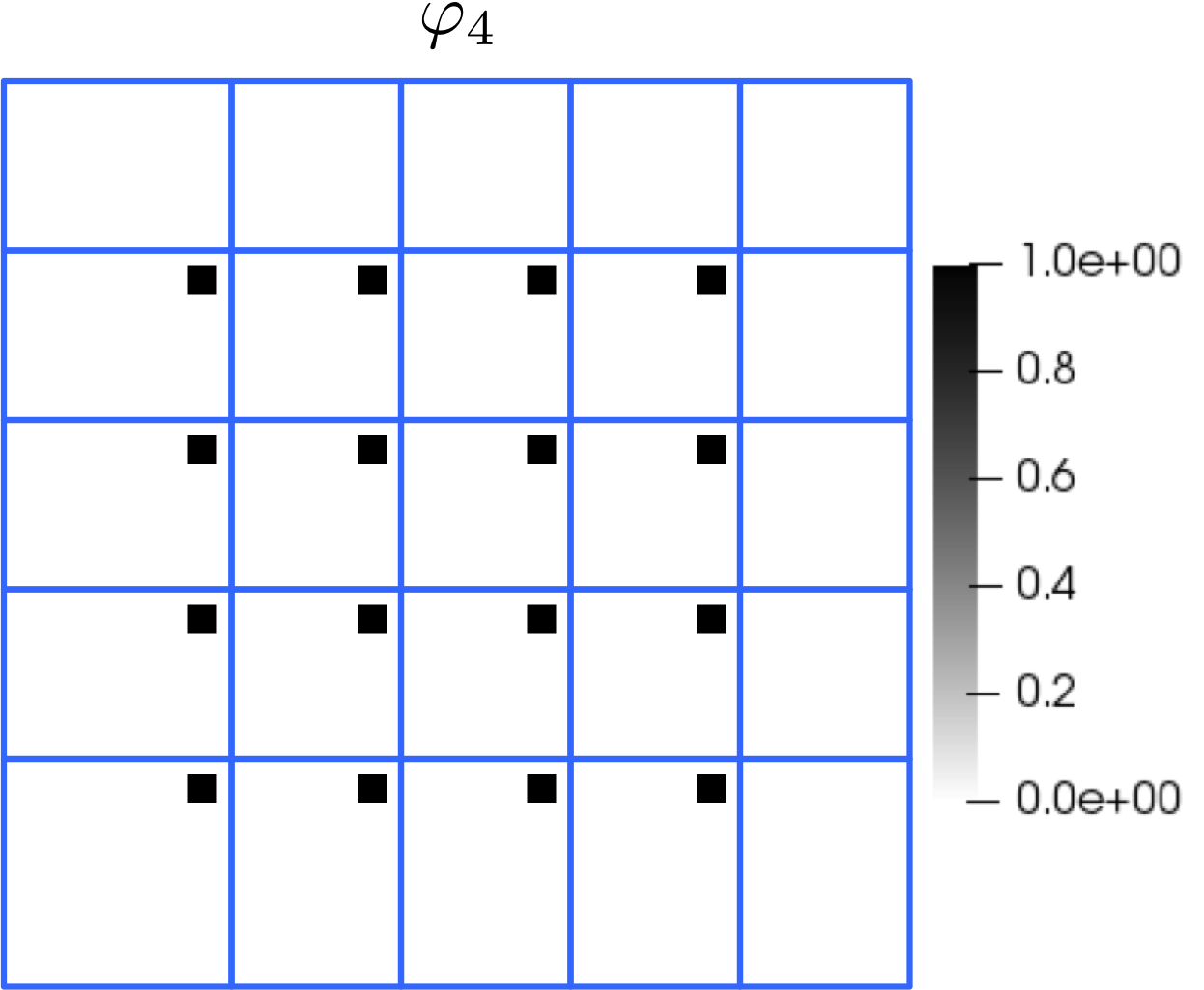}
\caption{The window functions $\varphi_i$ and the CDDs on which they vanish on the skeleton.}
\label{Fig::PhiDD}
\end{figure}

The only remaining restriction on source distribution is that it must not cross the CDD skeleton. This restriction can be overcome 
by considering three more CDDs obtained from horizontal and vertical shifts in the CDD, as illustrated in Figure~\ref{Fig::NewDD}. Using the original CDD and these new shifted CDDs, we define a partition of unity of four window functions $\varphi_i$, $i=1,2,3,4$, each of which vanish on the skeleton of \textit{one} of the four CDD (Figure~\ref{Fig::PhiDD}) and we define the corresponding windowed source densities $f^{\varphi_i}:=\varphi_i f$.
Clearly, $f^{\varphi_i}$ vanishes on the skeleton of one CDD. Using this CDD, we can therefore apply the procedure of Section~\ref{Sec::ArbSrcSkeleton} to obtain the global solution $u^{\varphi_i}$ corresponding to $f^{\varphi_i}$. The global solution for a general $f$ is them simply the sum of those solutions:
\begin{alignat*}{1}
u=\sum_{i=1}^4 u^{\varphi_i}.
\end{alignat*}

\subsection{Discrete formulation}
\label{Sec::LSweeps::Disc}
In this section, we show that the procedure introduced in Section~\ref{Sec::LSweeps::Cont} can be applied on the discrete level to approximately solve~\eqref{LS}. We start by defining the CDD and local problems. We define the CDD so that the skeleton does not intersect with any discretization point. Then, the skeleton clearly divides the global degrees-of-freedom into sets $\boldOmegaij$. We then define a local problem associated with each $\boldOmegaij$. This problem is defined on a superset of $\boldOmegaij$, which we denote $\boldOmegaijepsilon$, and is the discretization of the local problem defined on $\Omega_{ij}^\varepsilon$.  The interior boundaries are first extended by a $(2\delta)$-layer of degrees-of-freedom and then further extended by a PML region.
In the same manner as for the continuous problem, the wave speed for this problem is inherited from the global wave speed. We denote the local system matrices $\boldAij$. 

\begin{figure}[htp]
\centering
\includegraphics[scale=0.5]{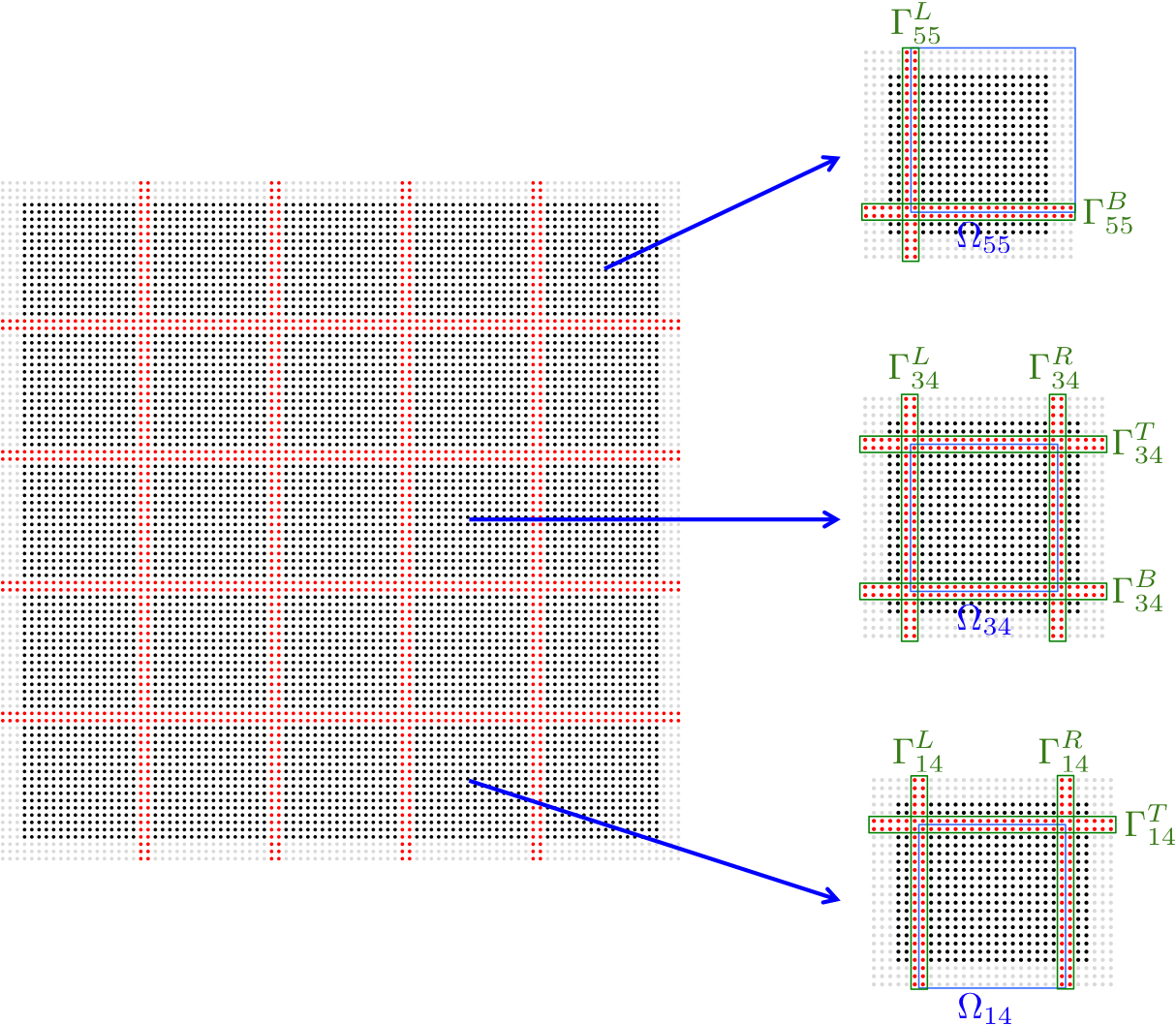}
\caption{The definition of the discrete local problems and the discrete interfaces.}
\label{Fig::DiscLocProbs}
\end{figure}

Furthermore, we define interfaces on these local problems to compute discrete polarized wavefields. For the continuous problem, these interfaces are defined as straight lines along the boundaries of the subdomains $\Omega_{ij}$. For the discrete problem, the necessary interfaces for computing discrete polarized wavefields consist of all degrees-of-freedom $\delta$-adjacent to the polarization interface. Therefore, we define the discrete counterparts of the interfaces $\Gamma_{ij}^{\ell}$, $\ell=L,R,T,B$, to be the $\delta$-adjacent degrees-of-freedom to $\Gamma_{ij}^{\ell}$. We call those sets $\boldGammaijl$, $\ell=B,R,T,L$. Figure~\ref{Fig::DiscLocProbs} illustrates these definitions. By construction, $\boldGammaijl$ defines the trace information that needs to be transferred between subdomains.  In the continuous problem, two neighboring local problems have to coincide in an $\varepsilon$-tube around the polarization interface to accurately extend the wavefield from one subdomain into the other. In the same way, the discrete problem requires two neighboring local problems to coincide in all $\delta$-adjacent degrees-of-freedom to the polarization interface. This justifies the $(2\delta)$-layer between $\boldOmegaij$ and the PML region. One $\delta$-layer is used for $\boldGammaijl$, and one $\delta$-layer is used to ensure that two neighboring local problems coincide in an appropriately sized neighborhood of $\boldGammaijl$.
  
In what follows we will show how the procedures introduced in Section~\ref{Sec::LSweeps::Cont} can be treated on the discrete level. Since the procedure introduced in Section~\ref{Sec::ArbSrcSkeleton} (scenario 3) is a generalization of all previous algorithms, we start with this procedure.
In the first stage, we restrict the global source vector $\boldf$ to the discretization points in $\boldOmegaij$ for each subdomain and define a local source vector $\boldfij$ on $\boldOmegaijepsilon$ such that $\boldfij|_{\boldOmegaij}=\boldf|_{\boldOmegaij}$ and $\boldfij$ is zero everywhere else. In the same way as in the continuous case, this local source vector can then be used to compute discrete local solutions $\bolduij$ in each subdomain:
\begin{alignat*}{1}
\bolduij=\boldAij^{-1}\boldfij.
\end{alignat*}
Using these local solutions, we define a global wavefield $\boldutilde$ such that $\boldutilde|_{\boldOmegaij}=\bolduij|_{\boldOmegaij}$. As with the continuous case, we extract the values of these local solutions on $\boldGammaijl$ and use them in stage 2 to update $\boldutilde$.

The second and third stage of the algorithm are realized by 
applying sweeps over the domain in the same way as in Section~\ref{Sec::ArbSrcSkeleton}. The only difference is that the trace information on $\Gamma_{ij}^{\ell}$ is replaced by its discrete counterparts $\boldGammaijl$, and the computation of continuous polarized wavefields is replaced by the computation of discrete polarized wavefields~\eqref{DiscPol} using the local system matrices $\boldAij$ and the appropriate interfaces $\boldGammaijl$. 

\begin{figure}[htp]
\centering
\includegraphics[scale=0.5]{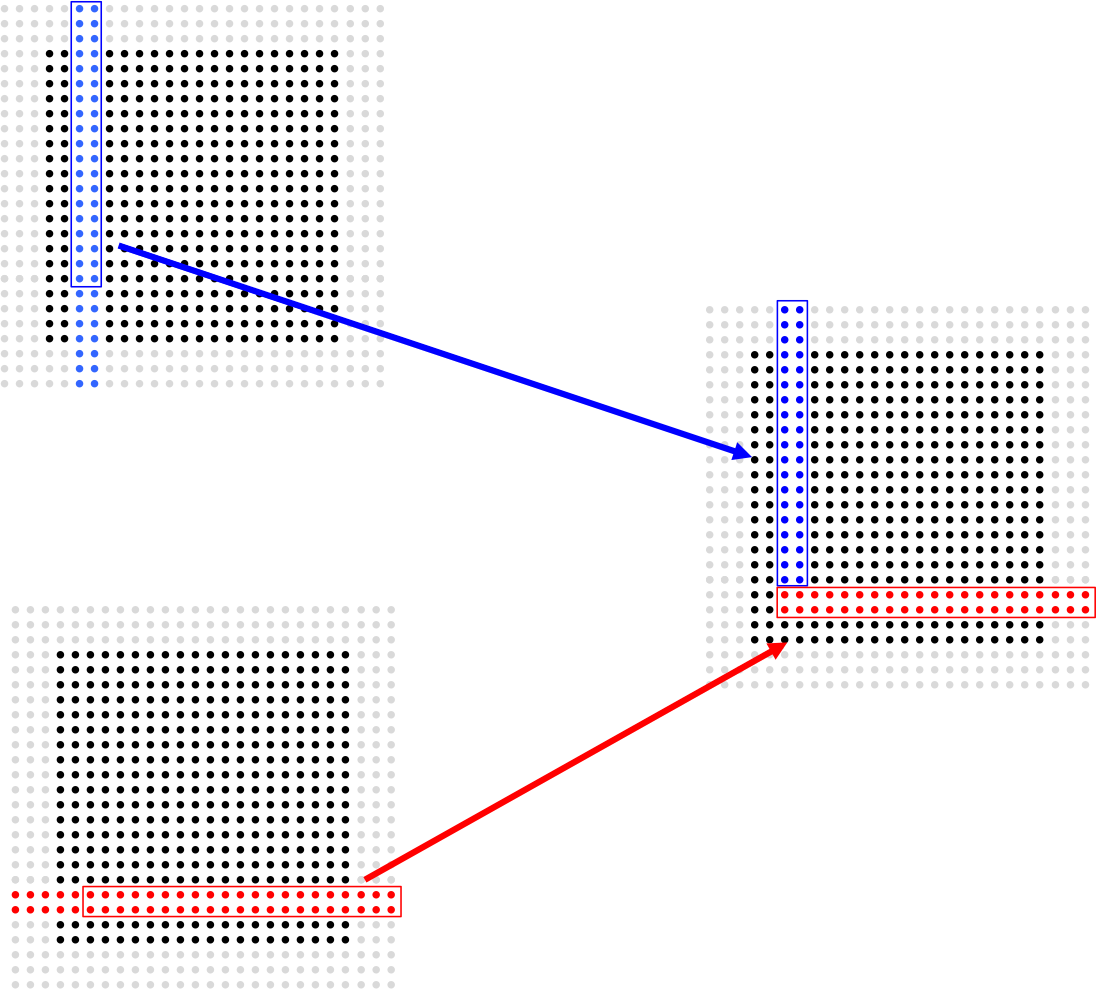}
\caption{The definition of the discrete L-shaped trace.}
\label{Fig::LTrace}
\end{figure}

Until now, we have elided discussion of a subtle point in the discrete algorithm: construction of the L-shaped trace information.   As an example, we consider the discrete counterpart of the trace $\Gamma_{ij}^{BL}$, but any other trace information can be handled analogously. Following the notation introduced in Section~\ref{Sec::LSweeps::Cont}, let us denote the trace information on $\boldGammaijB$ by $\boldsymbol{\lambda_{ij}^B}$ and the trace information on $\boldGammaijL$ by $\boldsymbol{\lambda_{ij}^L}$. Similarly to $\Gamma_{ij}^{BL}$ in Section~\ref{Sec::LSweeps::Cont}, $\boldsymbol{\Gamma_{ij}^{BL}}$ is defined as the L-shaped set of degrees-of-freedom surrounding the upper right quadrant of the set $\boldOmegaij$, as illustrated in Figure~\ref{Fig::LTrace}. The trace information $\boldsymbol{\lambda_{ij}^{BL}}$ is then defined such that $\boldsymbol{\lambda_{ij}^{BL}}|_{\boldGammaijB}=\boldsymbol{\lambda_{ij}^{B}}$ and $\boldsymbol{\lambda_{ij}^{BL}}|_{\boldGammaijL\backslash\boldGammaijB}=\boldsymbol{\lambda_{ij}^{L}}|_{\boldGammaijL\backslash\boldGammaijB}$.
Note that by convention we have defined the trace information of $\boldsymbol{\lambda_{ij}^{BL}}$ to coincide with $\boldsymbol{\lambda_{ij}^{B}}$ in $\boldGammaijB\cap\boldGammaijL$. This is an arbitrary choice, as the trace information on $\boldGammaijB\cap\boldGammaijL$ can be obtained from any of the two traces $\boldsymbol{\lambda_{ij}^{B}}$ or $\boldsymbol{\lambda_{ij}^{L}}$.

Following this procedure, we can compute approximations $\boldutilde$ of the global solution that only differ from $\boldu$ by errors induced by the discretization or PMLs, as long as the global source vector $\boldf$ is zero on the degrees-of-freedom on the skeleton of the CDD, i.e., the union of all $\boldGammaijl$. This procedure can be extended to entirely arbitrary source vectors using window functions, as introduced in Section~\ref{Sec::ArbSrc}.

\section{Implementation and Complexity}
\label{Sec::Complexity}
The procedure introduced in Section~\ref{Sec::LSweeps} produces an approximation of a global solution to~\eqref{LS}. This approach therefore defines the approximate solution operator $\mathbf{\tilde{A}}^{-1}$, such that $\mathbf{\tilde{A}}^{-1}(\boldsymbol{f})\approx \mathbf{A}^{-1}\boldsymbol{f}$. Thus, we use $\mathbf{\tilde{A}}^{-1}$ to precondition~\eqref{LS},
\begin{alignat}{1}
\label{PrecLS}
\mathbf{\tilde{A}}^{-1}\left(\boldA\boldu\right)=\mathbf{\tilde{A}}^{-1}\left(\boldf\right),
\end{alignat}
and use a Krylov subspace method, such as GMRES~\cite{Saad_Schultz:GMRES} or BiCG-stab~\cite{van_der_Vorst:BiCGSTAB}, to solve~\eqref{PrecLS}.

The resulting preconditioned iterative solver has the following properties:
\begin{itemize}
    \item the preconditioner $\mathbf{\tilde{A}}^{-1}$ can be applied with optimal parallel complexity, $\cO(N/p)$, and
    \item Krylov methods applied to~\eqref{PrecLS} converge (empirically) in ${\cal O}(\log\omega)$ iterations;
\end{itemize}
thus making it parallel scalable. 
In this Section we focus on the first property. In particular, we show the implementation of the preconditioner, while analyzing its complexity with respect to the computational and communication cost. For the second property, we will present extensive numerical evidence in Section~\ref{Sec::NumericalExamples}.

We consider a standard  communication model \cite{BallardDemmel:minimizing_communication_in_numerical_linear_algebra,DemmelGrigori:communication_avoiding_rank_revealing_qr_factorization_with_column_pivoting,Poulson_Demanet:a_parallel_butterfly_algorithm,Thakur:optimization_of_collective_communication_operations_in_mpich}. The model assumes that each process is only able to send or receive a single message at a time, though different messages can be sent and received asynchronously. A message of size $M$ can be communicated with $\alpha + \beta M$ complexity. The latency $\alpha$ represents the minimum complexity with which  an arbitrary message from one process to another can be communicated, and is a constant overhead for any communication. The inverse bandwidth $\beta$ is the complexity with which one unit of data can be communicated.

\begin{figure}[htp]
\centering
\includegraphics[scale=0.5]{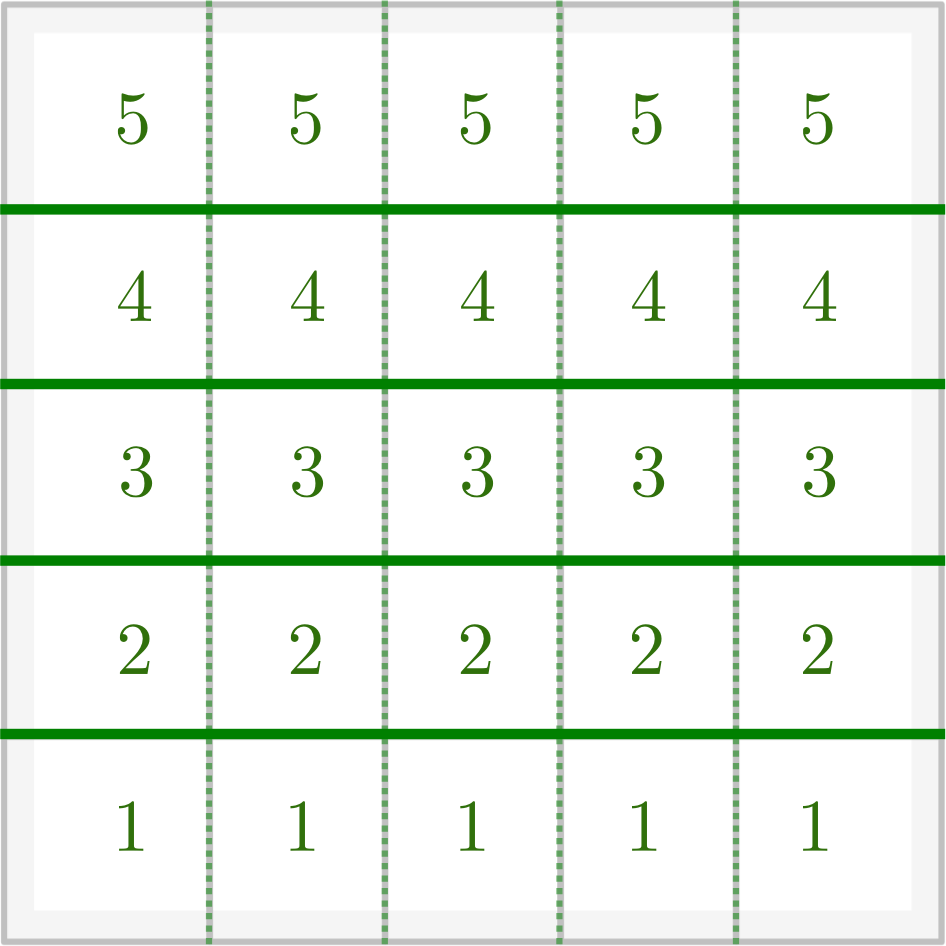}
\caption{The assignment of subdomains to processor in a row-based fashion.}
\label{Fig::ProcAssignment}
\end{figure}

We implement the algorithm within a distributed memory framework using MPI. For simplicity, we assign each row of subdomains to one MPI rank\footnote{ This restriction can be relaxed to exploit asynchronous parallelism models and subdomain pipelining, for extremely large problems.}.  
By assigning multiple subdomains to one rank, the preconditioner can be applied optimally in parallel, as described below.
Consider a diagonal sweep, for example from the bottom-left corner to the top-right corner, for a $q\times q$ CDD.  To maximize parallelism, each subdomain in a diagonal, perpendicular to the sweep direction, has to be processed in parallel. This can be realized by assigning subdomains to MPI ranks in a row-based fashion such that the $i$-th row of the CDD is processed by rank $i\operatorname{mod}p$. Thus, each rank is assigned one or several rows of subdomains and does all of their associated computation. An illustration of the assignment of subdomains to rank is provided in Figure~\ref{Fig::ProcAssignment}. In what follows, we assume this configuration, which allows us to also exploit the maximum parallelism in the application of the preconditioner and obtain a parallel scalable solver. The only restriction in this setup is that the number of ranks $p$ used in this setup is bounded by the number of subdomains in one column, i.e., $p\leq q={\cal O}(n)$.

For analysis purposes, we divide the solver for~\eqref{LS} in three phases:
\begin{enumerate}
    \item \textit{Setup}  the local problems, system matrices, right-hand sides, and the corresponding window functions $\varphi_{i}$;
    \item \textit{Factorize} the local system matrices; and
    \item \textit{Solve} the linear system using a Krylov method. 
\end{enumerate}

\bigskip
\noindent
{\it Phase 1: Setup}\\
For two-dimensional problems, the size of the local problems associated with the subdomain $\Omega_{ij}$ is ${\cal O}(1)$. A crucial requirement to ensure that these local problems are indeed of $\cO(1)$ size is that the thickness of the PML region, in wavelengths normal to $\Omega_{ij}$, is held constant with mesh-refinement, i.e., the number of degrees-of-freedom in the PML does not change and the spatial PML-width scales as $O(1/\omega)$. Consequently, the local system matrices $\boldAij$, right-hand sides, and the window functions $\varphi_{ij}$ can be assembled with ${\cal O}(1)$ computational complexity. 
In a parallel computational environment, each of the ${\cal O}(q^2)$ subdomains can be processed independently and no communication is required between subdomains. Thus all rows are processed with a ${\cal O}(n)$ parallel computational complexity, because $q=\cO(n)$ and $p=\cO(q)$. For two-dimensional problems, the computation in the setup phase is therefore realized with 
optimal $\cO(N/p)$ parallel computational complexity. For three-dimensional problems, the size of the local problems associated with each subdomain is $\cO(n)$, and the same arguments also yield optimal $\cO(N/p)$ parallel complexity.  In either case, there is no communication cost during the setup phase.

\bigskip
\noindent
{\it Phase 2: Factorize}\\
For two-dimensional problems, the system matrices of the local problems have size ${\cal O}(1)$ and therefore can be factorized with $\cO(1)$ computational complexity. As in stage 1, each subdomain can be processed independently with zero communication, thus we can follow the same argument as in the setup phase to justify that all of the $q^2$ subdomains can be factorized with $O(N/p)$ parallel computational complexity. For three-dimensional problems, the system matrices of the local problems have size ${\cal O}(n)$ and are quasi-one-dimensional. Thus, standard sparse direct solvers can be used to factorize the matrices with optimal computational complexity \cite{GeorgeNested_dissection,Duff_Reid:The_Multifrontal_Solution_of_Indefinite_Sparse_Symmetric_Linear,Hockney:cyclic_reduction}, i.e., ${\cal O}(n)$. The same arguments as for the two-dimensional case can be employed to show that, in three dimensions, the $q^2$ subdomains can be factorized with $O(N/p)$ parallel computational complexity.  In either case, there is no communication cost during the factorization phase.

\bigskip
\noindent
{\it Phase 3: Solve }\\
We use the GMRES method to solve the linear system~\eqref{LS}. Each iteration of this method consists of three main operations: the application of $\boldA$, the computation between the Ritz vectors, and the application of the preconditioner $\boldsymbol{\tilde{A}}^{-1}$. For each of the three operations, we assume that the global vectors are provided in a distributed fashion such that each subdomain $\Omega_{ij}$ holds the values of the global vector corresponding to $\boldOmegaij$.

The application of $\boldA$ can be realized by simply applying the local system matrices to a local vector, as long as each subdomain holds all degrees-of-freedom associated with $\boldOmegaij$ and $\boldGammaijl$ for $\ell=B,R,T,L$. This requires some communication because the subdomains only store the degrees-of-freedom in $\boldOmegaij$. Considering the row-based assignments of subdomains to MPI ranks, the sets $\boldGammaijB$ have to be communicated to the subdomain $\Omega_{(i-1)j}$, and the sets $\boldGammaijT$ have to be communicated to the subdomain $\Omega_{(i+1)j}$. In practice, for two-dimensional problems, each row has to communicate ${\cal O}(n)$ values to a neighboring MPI rank. Using the above communication model, this communication has complexity
\begin{alignat*}{1}
{\cal O}(\alpha + \beta n)=\cO(n).
\end{alignat*}
Due to the row-based rank assignment, up to $p$ of the $q$ rows can simultaneously exchange information, resulting in an overall ${\cal O}(qn/p)=\cO(N/p)$ parallel communication complexity. The same arguments can be applied for the communication of the sets $\boldGammaijT$. Thus, for two dimensional problems, all necessary information can be communicated with optimal $\cO(N/p)$ parallel complexity. In three dimensions, each row has to communicate ${\cal O}(n^2)$ values to a neighboring MPI rank. Following the same arguments as before all necessary information can be communicated with optimal $\cO(N/p)$ parallel complexity.
Once all information is communicated, the local system matrices $\boldAij$ are applied to obtain the action of the global system matrix $\mathbf{A}$ on the degrees-of-freedom in $\boldOmegaij$. This application can be realized with $\cO(1)$ and $\cO(n)$ computational complexity per subdomain in two and three dimensions, respectively. Again, each of the ${\cal O}(q^2)$ subdomains can be processed independently resulting in the total $\cO(N/p)$ parallel computational and communication complexity to apply the matrix $\boldA$ for two- and three-dimensional problems.

At iteration $k$, $k$ inner products between the Ritz vectors need to be computed. Each inner product can be computed efficiently in parallel, given that each subdomain contains the local components of the Ritz vectors. Thus, each inner product are computed with $O(N/p)$ parallel computational complexity and $O(p)$ parallel communication complexity. The application of $\boldsymbol{\tilde{A}}^{-1}$ is analyzed in detail in Section~\ref{Sec::Complexity::Preconditioner}, where we show that it can be applied with $\cO(N/p)$ parallel computational and communication complexity. The $k$-th GMRES iteration can therefore be realized with ${\cal O}(kN/p)$ parallel computational complexity and ${\cal O}(N/p)$ parallel communication complexity\footnote{To achieve the parallel communication complexity, we implicitly assume that $\cO(kp)=\cO(N/p)$.}. Due to the effectiveness of the preconditioner, the number of iterations only grows as ${\cal O}(\log\omega)$.  This claim is corroborated by the numerical examples in Section~\ref{Sec::NumericalExamples}. Using restarts, the GMRES method can therefore be applied with a total parallel computational and communication complexity of ${\cal O}\big( (N/p)\log\omega\big)$. The computation and communication complexities are summarized in Tables~\ref{table:complexity2D} and~\ref{table:complexity3D}, for two- and three-dimensional problems.

\begin{table}[H]
\begin{center}
\begin{tabular}{|c|c|c|c|c|}
\hline
 & \multicolumn{2}{|c|}{Computation} & & \\ 
Step 		  &  per subdomain & per rank & Communication    & Total time  \\
\hline
Set up	      & ${\cal O}(1)$ & $\cO(n)=O(N/p)$ & $\cO(1)$	& ${\cal O}(N/p)$    \\
\hline
Factorization	      & ${\cal O}(1)$ & $\cO(n)=\cO(N/p)$ & $\cO(1)$	& ${\cal O}(N/p)$    \\
\hline
Solve         & ${\cal O}(\log \omega )$  &  $\cO(n\log\omega)=\cO\big((N/p)\log\omega\big)$ & $ {\cal O}\big((N/p) \log\omega\big)$   & $ {\cal O}\big((N/p) \log\omega\big)$    \\
\hline
\end{tabular}
\caption{Parallel complexity of the different stages for two-dimensional problems, where $N$ is the total number of discretization points, $n$ is the number of discretization points per dimension, $\omega$ is the frequency following $\omega \sim n$, and $p$ is the number of processors following $p=\cO(n)$. }\label{table:complexity2D}
\end{center}
\end{table}

\begin{table}[H]
\begin{center}
\begin{tabular}{|c|c|c|c|c|}
\hline
 & \multicolumn{2}{|c|}{Computation} & & \\ 
Step 		  &  per subdomain & per rank & Communication    & Total time  \\
\hline
Set up	      & ${\cal O}(n)$ & $\cO(n^2)=O(N/p)$ & $\cO(1)$	& ${\cal O}(N/p)$    \\
\hline
Factorization	      & ${\cal O}(n)$ & $\cO(n^2)=O(N/p)$ & $\cO(1)$	& ${\cal O}(N/p)$    \\
\hline
Solve         & ${\cal O}(n\log \omega )$  &  $\cO(n^2\log\omega)=\cO\big((N/p)\log\omega\big)$ & $ {\cal O}\big((N/p) \log\omega\big)$  & $ {\cal O}\big((N/p) \log\omega\big)$     \\
\hline
\end{tabular}
\caption{Parallel complexity of the different stages for three-dimensional problems, where $N = n^3$ are the total number of discretization points, and the number of discretization points per dimension, $\omega$ is the frequency following $\omega \sim n$, and $p$ is the number of processors following $p=\cO(n)$. }\label{table:complexity3D}
\end{center}
\end{table}

\subsection[]{Application of $\mathbf{\tilde{A}^{-1}}$} 
\label{Sec::Complexity::Preconditioner}

The application of $\boldsymbol{\tilde{A}}^{-1}$ can be divided in four steps:
\begin{enumerate}
    \item[(a)] Computation of the local solutions,
    \item[(b)] Extension of the local solutions into the same row,
    \item[(c)] Extension of the local solutions into the same column, and
    \item[(d)] Extension of the local solutions into the rest of the subdomains.    
\end{enumerate}
As before, we assume that at the beginning and at the end of each step, each rank holds all the values corresponding to the degrees-of-freedom in $\boldOmegaij$ for all of its associated subdomains.

The computation of the local solutions requires to solve a local problem in each subdomain. Each of the local problems can be solved in ${\cal O}(1)$ and $\cO(n)$ computational complexity for two- and three-dimensional problems, respectively\footnote{Similar to the factorization, standard sparse direct solvers are used to solve the quasi-one-dimensional problems arising from three-dimensional problems with optimal computational complexity, i.e., ${\cal O}(n)$.}. All local solutions are computed independently and $p$ local solutions can be computed simultaneously.  Therefore, for both two- and three-dimensional problems, step $1$ is realized with ${\cal O}(N/p)$ parallel computational complexity and there is no communication cost.

The extension of the local solutions into the subdomains contained in the same row requires one left- and one right-sweep in each row. Due to the row-based rank assignment, each of those sweeps must be realized sequentially. However, each row can be processed independently, with no inter-row communication. For example, consider the right-sweep in one row. In each subdomain, a (discrete) polarized wavefield has to be computed (i.e., the local problem has to be solved), 
the local solution has to be updated, and the traces to be transferred to the right have to be extracted. 
For two-dimensional problems, each of those stages can be realized with ${\cal O}(1)$ computational complexity per subdomain because the size of the subdomains is $\cO(1)$. Since each subdomain in the row has to be processed sequentially and there are $q={\cal O}(n)$ subdomains in one row, this procedure can be realized with $\cO(n)$ computational complexity per row. Employing this procedure concurrently in each row yields $\cO(N/p)$ parallel computational complexity and no communication cost. For three-dimensional problems, each stage can be realized with  $\cO(n)$ computational complexity per subdomain\footnote{Again, standard sparse direct solvers are used to solve the quasi-one-dimensional problems arising from three-dimensional problems with optimal complexity, i.e., ${\cal O}(n)$.}.  Because each of the $q=\cO(n)$ subdomains in one row is processed sequentially, this results in a $\cO(n^2)$ computational complexity in every row.  Thus, processing each row concurrently results in the $\cO(N/p)$ parallel computational complexity and no communication cost. Analysis of the leftward sweep follows the same logic and reaches the same conclusions.

Extending the local solutions from one subdomain into the other subdomains within the same column requires one upward- and one downward-sweep. Given that the processor assignment does not allow for each column to be processed independently, the organization of the computation is non-trivial in order to reveal parallelism. As an example, we explain the implementation of an upward sweep. We begin by considering the subdomain in the bottom-left corner, $\Omega_{11}$.  Once $\Omega_{11}$ is processed and its information sent to $\Omega_{21}$, rank 0 can process the subdomain $\Omega_{12}$ and rank 1 can process subdomain $\Omega_{21}$. Following this procedure, the upward sweep is applied by sweeping over the domain, essentially pipelining the rows. Parallelism is achieved by processing diagonally-adjacent sets of subdomains simultaneously.  By construction, each diagonal set has at most $q$ subdomains. For two-dimensional problems, each subdomain can be processed with ${\cal O}(1)$ computational complexity, independent from all other subdomains. Thus, an entire diagonal set can be processed with ${\cal O}(q/p)$ parallel computational complexity. For three-dimensional problems, each subdomain can be processed with $\cO(n)$ computational complexity, resulting in a ${\cal O}(nq/p)$ parallel computational complexity to process each diagonal. Since there are ${\cal O}(q)$ diagonals, and $q=\cO(n)$, this results in the total $\cO(N/p)$ parallel computational complexity to apply the upward sweep for both two- and three-dimensional problems.  

In contrast with previous phases and steps, the upward and downward sweeps also require communication. For each subdomain $\Omega_{ij}$, the traces $\boldGammaijB$ have to be communicated from the subdomain $\Omega_{(i-1)j}$. For two-dimensional problems, the information transfers are performed in parallel, and all $\cO(n)$ values, from all $\boldGammaijB$s, are communicated with $\cO(n/p)$ parallel communication complexity. For three-dimensional problems, the communication volume is $\cO(n^2)$ and, using the same arguments, the parallel communication complexity, per diagonal set, is $\cO(n^2/p)$.  For both two- and three-dimensional problems, considering all $\cO(q)=\cO(n)$ diagonals, we find the optimal parallel $\cO(N/p)$ communication complexity.  The analysis of the downward sweep follows the same argument and achieves the same conclusion.

\begin{figure}[htp]
\centering
\includegraphics[scale=0.55,trim=0.2in 0.2in 0.2in 0.2in,clip]{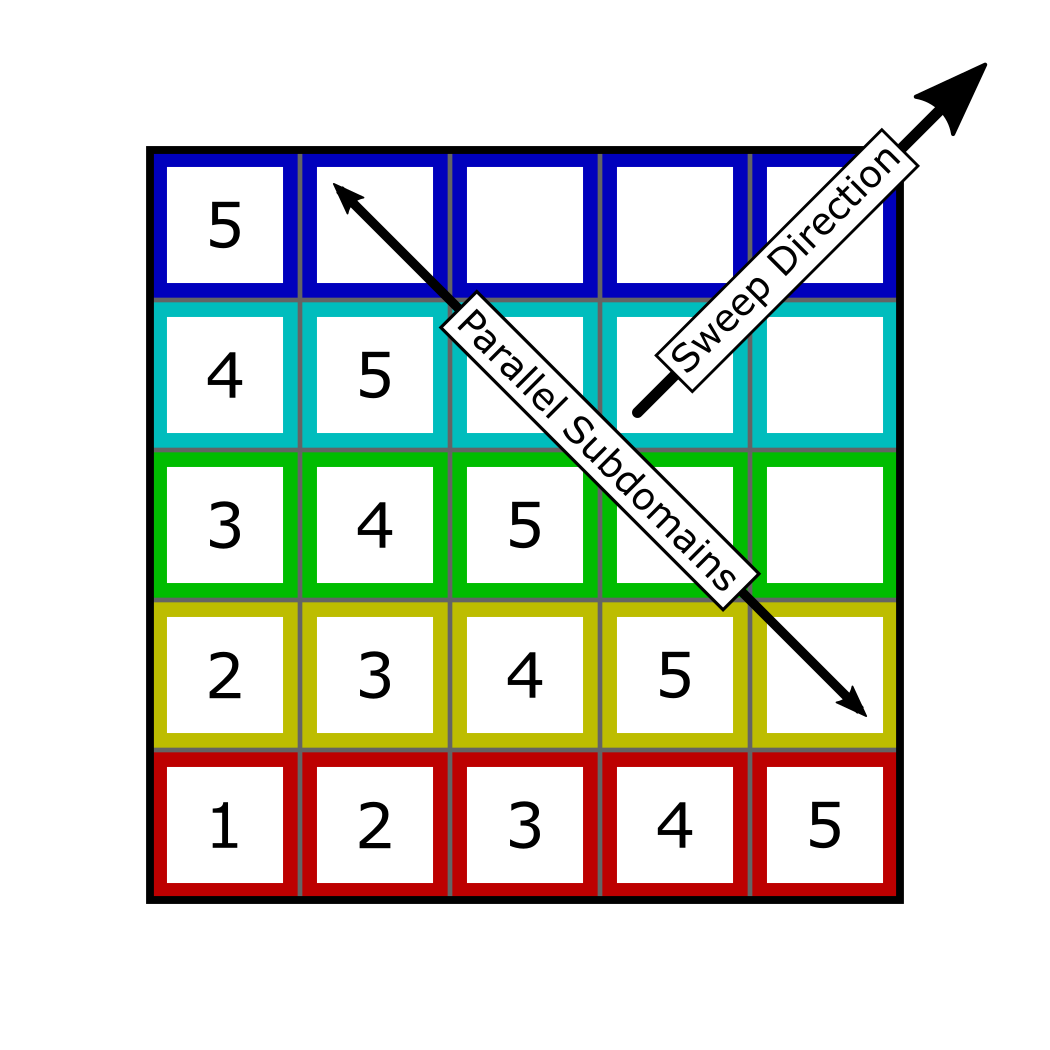}
\caption{Illustration of parallelism in diagonal sweeps.  Subdomains processed in parallel have same label.  Subdomains with like color are assigned to same processor.}
\label{Fig::SweepDirection}
\end{figure}

The extension of the local solutions into the remaining subdomains, those not in the same row or column, requires sweeps over the CDD from corner to corner, along both diagonals. In these diagonal sweeps, each step requires processing of sets of subdomains that are diagonally-adjacent along the direction perpendicular to the sweep direction as illustrated in Fig.~\ref{Fig::SweepDirection}. This process has the same computation and communication patterns as the upward and downward sweeps and thus, following the same analysis, have the same complexities: they are applied with $O(N/p)$ parallel computational complexity and $O(N/p)$ parallel communcation complexity.

In summary, each step has $O(N/p)$ parallel computational complexity and at most $O(N/p)$ parallel communication complexity, thus the preconditioner can be applied with $\cO(N/p)$ parallel computational and communication complexity. The computational complexity and communication complexities, for applying the preconditioner, are summarized in Tables~\ref{table:Stage3Complexity2D} and~\ref{table:Stage3Complexity3D}.

\begin{table}[H]
\begin{center}
\begin{tabular}{|c|c|c|c|c|}
\hline
 & \multicolumn{2}{|c|}{Computation} & & \\
Step 		  &  per subdomain & per rank & Communication    & Total time  \\
\hline
Part 1	      & ${\cal O}(1)$ & $\cO(n)=O(N/p)$ & $\cO(1)$	& ${\cal O}(N/p)$    \\
\hline
Part 2	      & ${\cal O}(1)$ & $\cO(n)=O(N/p)$ & $\cO(1)$	& ${\cal O}(N/p)$    \\
\hline
Part 3	      & ${\cal O}(1)$ & $\cO(n)=O(N/p)$ & $\cO(N/p)$	& ${\cal O}(N/p)$    \\
\hline
Part 4	      & ${\cal O}(1)$ & $\cO(n)=O(N/p)$ & $\cO(N/p)$	& ${\cal O}(N/p)$    \\
\hline
\end{tabular}
\caption{Parallel complexity of the different stages for two-dimensional problems for $p=\cO(n)$.}\label{table:Stage3Complexity2D}
\end{center}
\end{table}

\begin{table}[H]
\begin{center}
\begin{tabular}{|c|c|c|c|c|}
\hline
 & \multicolumn{2}{|c|}{Computation} & & \\
Step 		  &  per subdomain & per rank & Communication    & Total time  \\
\hline
Part 1	      & ${\cal O}(n)$ & $\cO(n^2)=O(N/p)$ & $\cO(1)$	& ${\cal O}(N/p)$    \\
\hline
Part 2	      & ${\cal O}(n)$ & $\cO(n^2)=O(N/p)$ & $\cO(1)$	& ${\cal O}(N/p)$    \\
\hline
Part 3	      & ${\cal O}(n)$ & $\cO(n^2)=O(N/p)$ & $\cO(N/p)$	& ${\cal O}(N/p)$    \\
\hline
Part 4	      & ${\cal O}(n)$ & $\cO(n^2)=O(N/p)$ & $\cO(N/p)$	& ${\cal O}(N/p)$    \\
\hline
\end{tabular}
\caption{Parallel complexity of the different stages for three-dimensional problems when $p=\cO(n)$.}\label{table:Stage3Complexity3D}
\end{center}
\end{table}

\section{Heterogeneities}
\label{Sec::LSweeps::Heterogeneities} 
As developed in  Section~\ref{Sec::LSweeps}, the wavefield computed by the proposed algorithm  can be interpreted as a sum of global solutions, each of which is induced by a global source density supported in a single subdomain only. We denote a source density supported in the subdomain $\Omega_{ij}$ by $\boldfij$ and the global wavefield induced by this source density as $\boldu^\boldfij$. Accordingly, the global wavefield $\boldu$ can be written as
\begin{alignat*}{1}
\boldu=\sum_{i=1}^q\sum_{j=1}^r \boldu^{\boldfij}.
\end{alignat*}
To understand the effect of heterogeneities on the effectiveness of the preconditioner, it is sufficient to consider a wavefield $\boldu^{\boldfij}$ due to  a single point source in detail. We consider two examples. The first example reveals the effect of reflections on the wavefield constructed by the preconditioner. The second example shows the effects of reflections in more pathological media.

\begin{figure}[htp]
\centering
\begin{tabular}{c|c}
\begin{subfigure}[b]{0.45\textwidth}
\centering\includegraphics[scale=0.3]{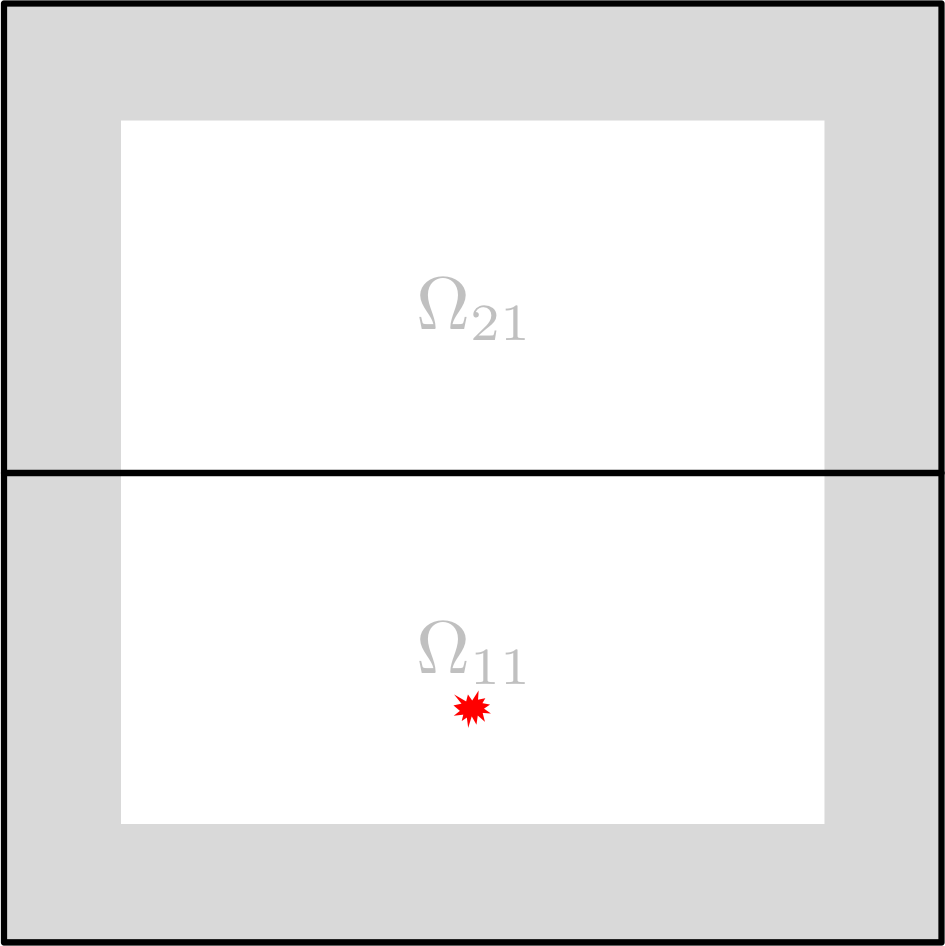}
\subcaption{The CDD.}
\label{Fig::Setup11}
\end{subfigure} & \begin{subfigure}[b]{0.45\textwidth}
\centering\includegraphics[scale=0.3]{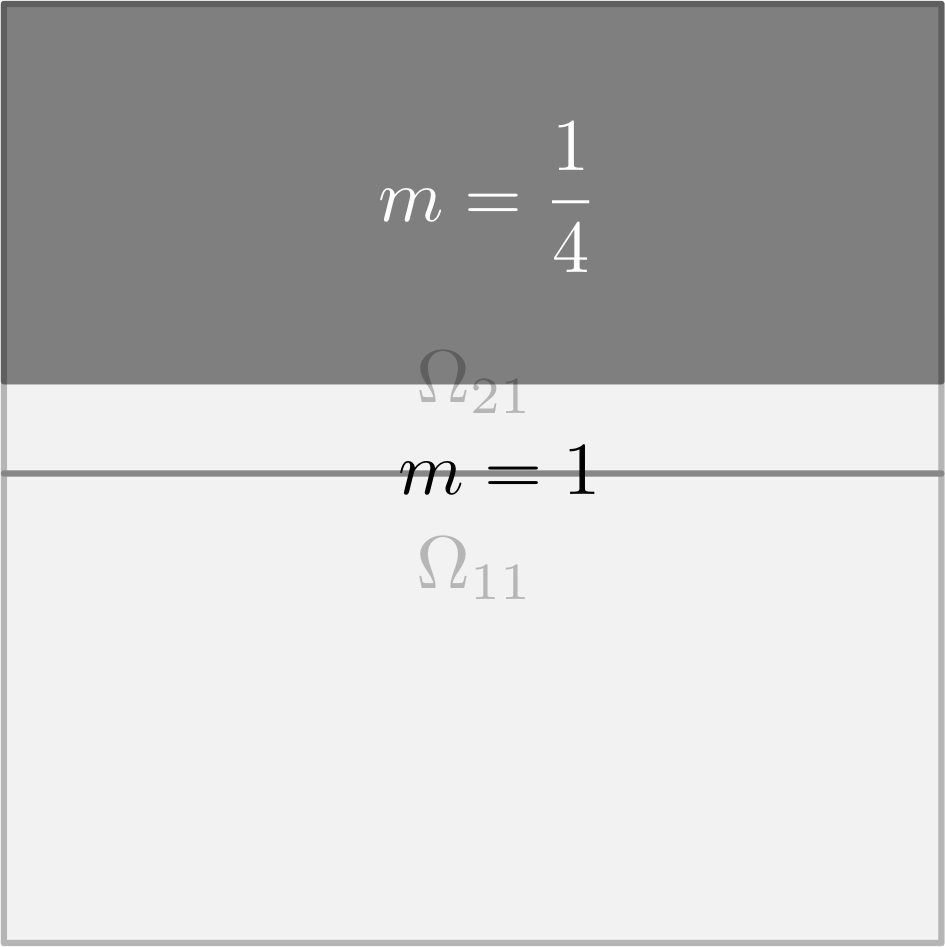}
\subcaption{The wave speed.}
\label{Fig::Setup12}
\end{subfigure} 
\end{tabular}
\caption{Illustration of a $2\times1$ CDD with associated wavespeed and a point source (red star).}
\label{Fig::Setup1}
\end{figure}

For the first example, we consider a simple domain decomposition with two subdomains, $\Omega_{11}$ and $\Omega_{21}$, divided by a horizontal line. We consider a layered wave speed distribution with two layers where the line of discontinuity is a horizontal line in $\Omega_{21}$. The setup is illustrated in Figure~\ref{Fig::Setup1}. 

\begin{figure}[htp]
\centering
\begin{tabular}{c|c}
\begin{subfigure}[b]{0.45\textwidth}
\centering\includegraphics[scale=0.2]{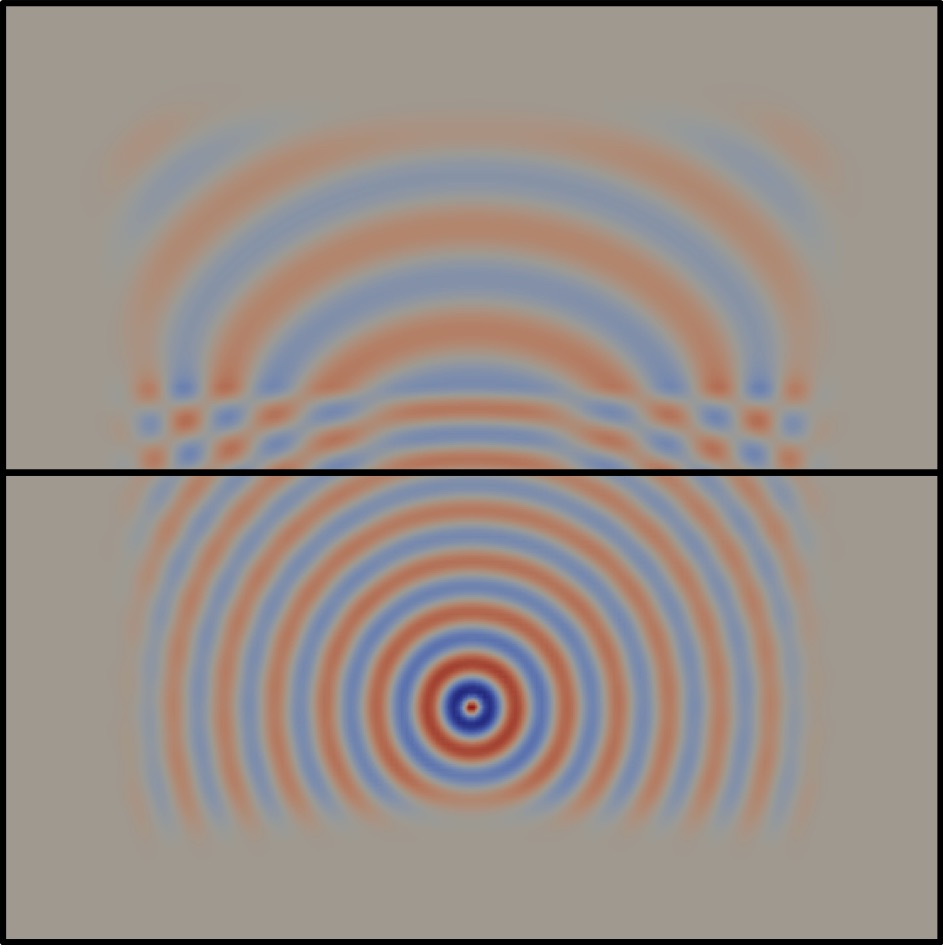}
\subcaption{Before the first GMRES iteration.}
\label{Fig::Sol11}
\end{subfigure} & \begin{subfigure}[b]{0.45\textwidth}
\centering\includegraphics[scale=0.2]{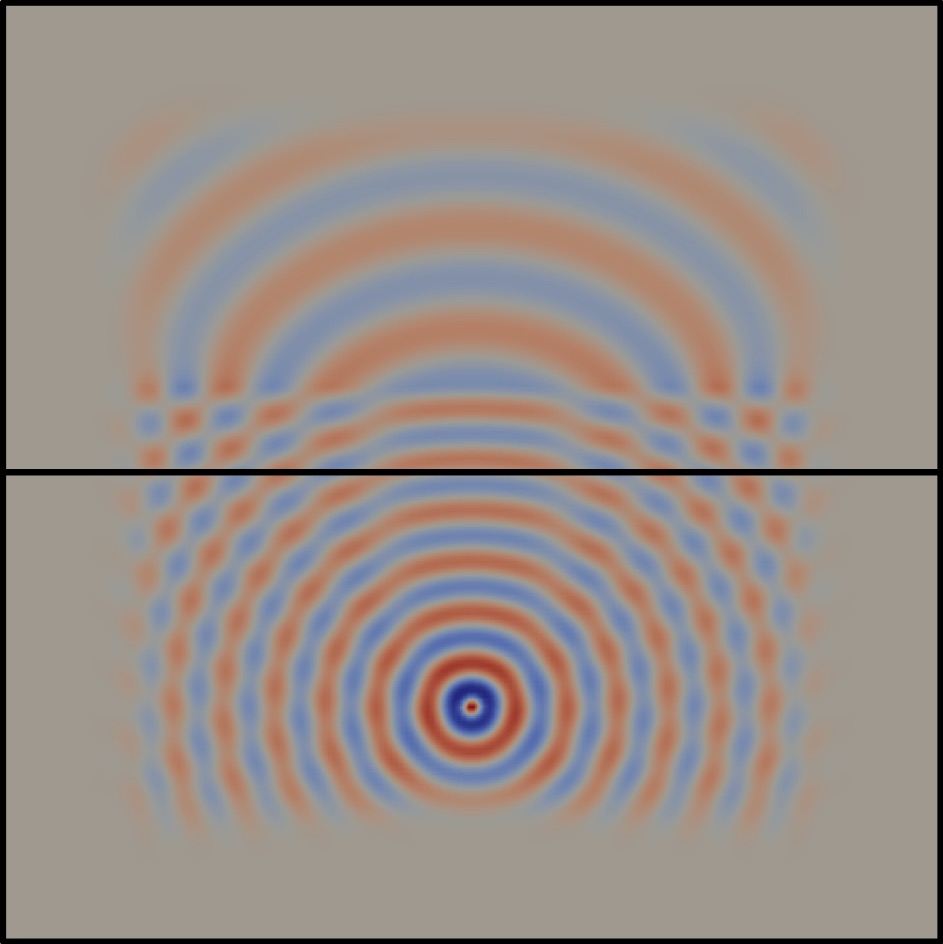}
\subcaption{After the first GMRES iteration).}
\label{Fig::Sol12}
\end{subfigure} 
\end{tabular}
\caption{Illustration of the computed wave fields.}
\label{Fig::Sol1}
\end{figure}

In this example, excluding local partial wavefield computation, only the upward sweep computes a non-zero partial wavefield. The solution obtained from this sweep is shown in Figure~\ref{Fig::Sol1}. The wavefield in the bottom subdomain is a solution to a problem with constant wave speed, and the reflection due to the discontinuity in the wave speed is clearly visible in the field in the top subdomain.  The solutions in each subdomain can be explained in the following way: the local wavefield for $\Omega_{11}$ from stage 1 is computed by solving the local problem in $\Omega_{11}$ for the point source. Therefore the discontinuity in $\Omega_{21}$ is not visible and the wavefield for the constant wave speed is computed. This results in a poor approximation of the global solution in  $\Omega_{11}$. Nevertheless, the top traces $\boldsymbol{\lambda}^T_{11}$ are extracted. In the upward sweep (as part of stage 2), these traces are used in $\Omega_{21}$ to compute a polarized wavefield. This local wavefield is a good approximation of the global solution in $\Omega_{21}$ for two reasons:
\begin{itemize}
    \item The discontinuity in the wave speed is visible to the local problem in $\Omega_{21}$, and
    \item The trace $\boldsymbol{\lambda}_{11}^T$ contains all required information from $\Omega_{11}$ so that the polarized wavefield computed in $\Omega_{21}$ is a good approximation of the global wavefield.
\end{itemize}
The resulting wavefield is therefore a good approximation of the global wavefield only in $\Omega_{21}$. As shown in Figure \ref{Fig::Sol11}, there is a mismatch between the two local solutions, which produces a large residual concentrated at the interface. 

Following the reasoning in \cite{Stolk,Chen_Xiang:a_source_transfer_ddm_for_helmholtz_equations_in_unbounded_domain}, this residual contain the information necessary to propagate the wavefield locally between subdomains. The solutions obtained from applying the preconditioner to the right-hand side, and the solution after the first GMRES iteration are shown in Figure~\ref{Fig::Sol1}. It is clear that after the second iteration, the GMRES method has converged to a good approximation of the wavefield.

This example shows that the preconditioner does not construct a good approximation of the global solution in rough media.  This is because the preconditioner tracks the physical behavior of some waves propagating through the domain, including waves propagating in a straight line, intra-subdomain reflections, and some refracted waves, but does not account for any inter-subdomain reflections. Nevertheless, since the preconditioner is used as part of an iterative solver for the linear system, the inter-subdomain reflections and unresolved refractions will be treated in subsequent iterations. 

\begin{figure}[htp]
\centering
\begin{tabular}{c|c}
\begin{subfigure}[b]{0.45\textwidth}
\centering\includegraphics[scale=0.3]{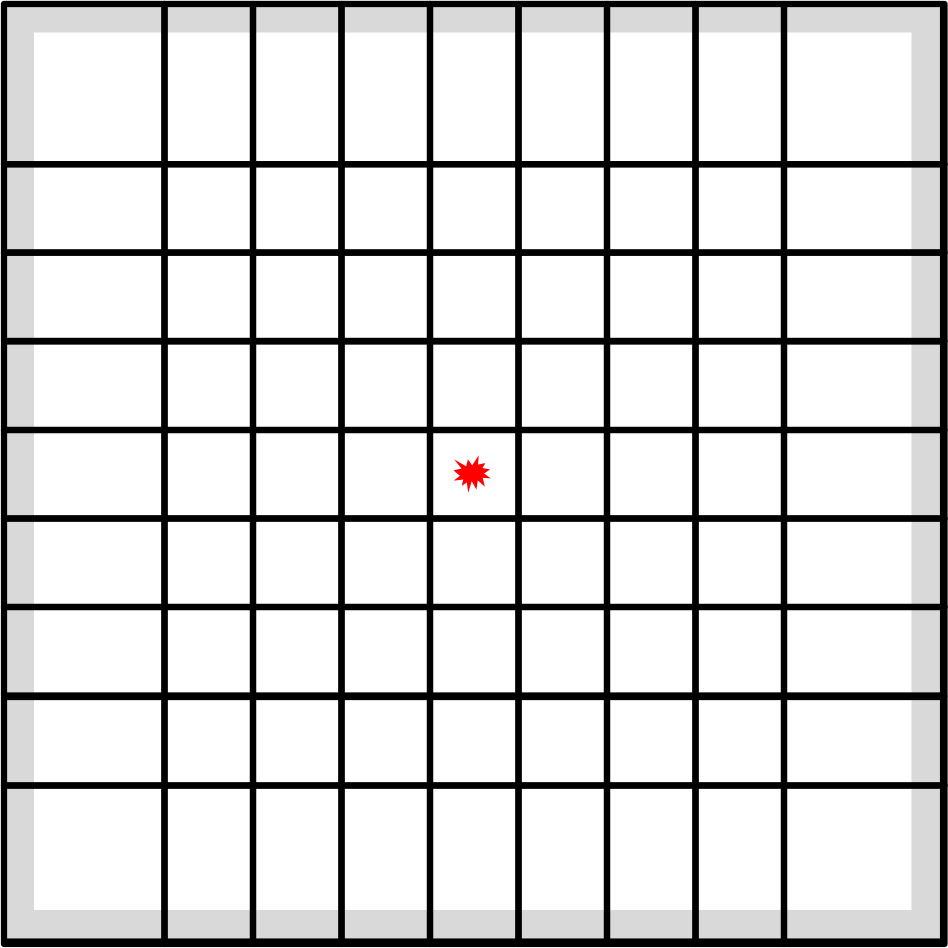}
\subcaption{The CDD.}
\label{Fig::Setup21}
\end{subfigure} & \begin{subfigure}[b]{0.45\textwidth}
\centering\includegraphics[scale=0.3]{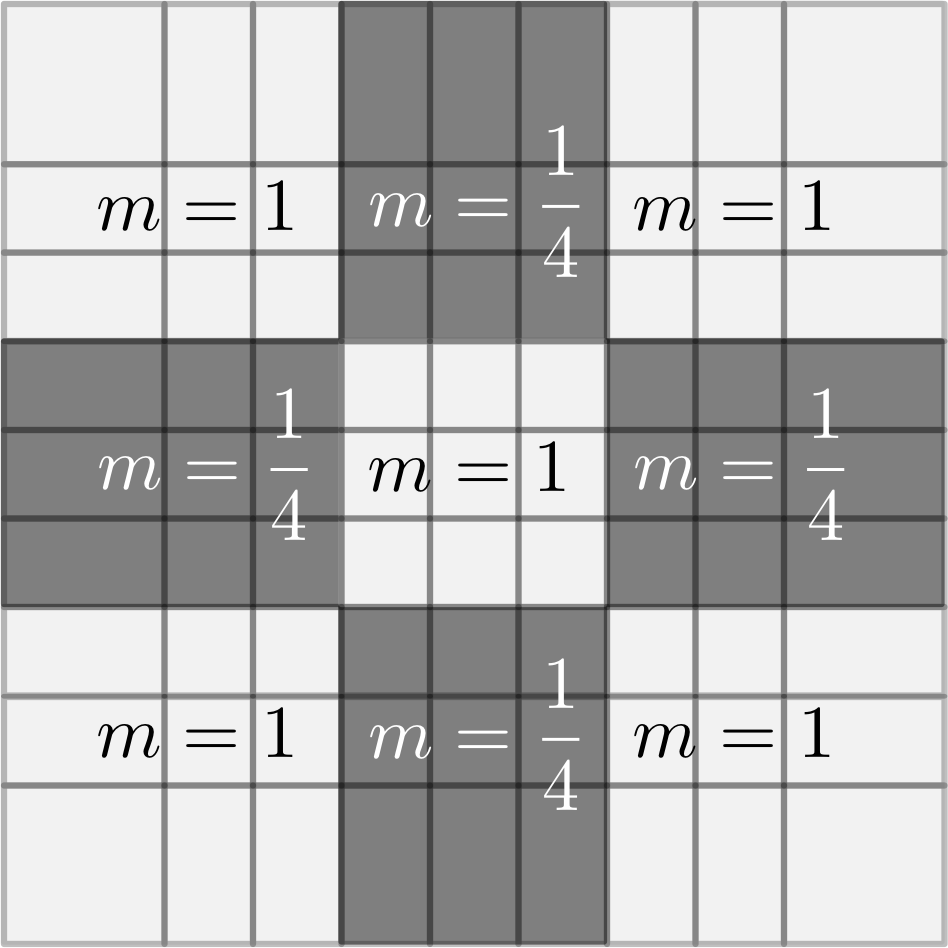}
\subcaption{The wave speed.}
\label{Fig::Setup22}
\end{subfigure} 
\end{tabular}
\caption{Illustration of the setup for the example involving a $9\times9$ CDD. The point source is shown by the red star, the PML region is shown in gray.}
\label{Fig::Setup2}
\end{figure}

\begin{figure}[htp]
\centering
\begin{tabular}{c|c}
\begin{subfigure}[b]{0.45\textwidth}
\centering\includegraphics[scale=0.2]{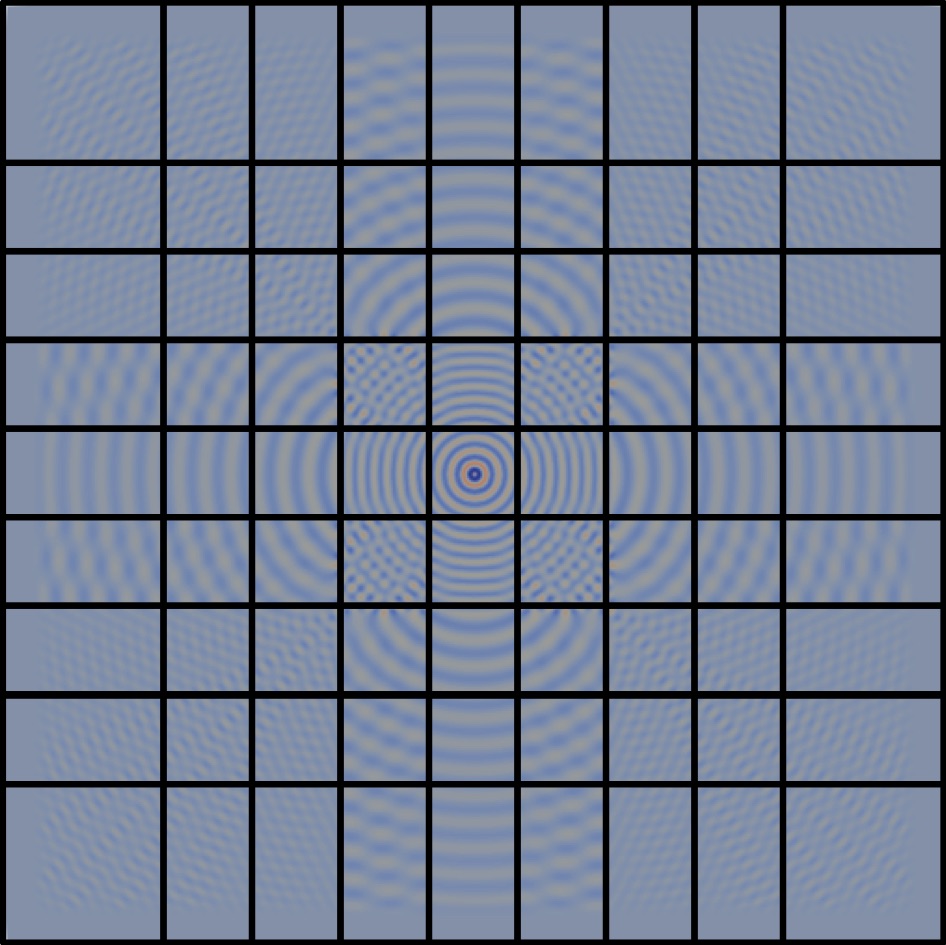}
\subcaption{Iteration 0}
\end{subfigure} & \begin{subfigure}[b]{0.45\textwidth}
\centering\includegraphics[scale=0.2]{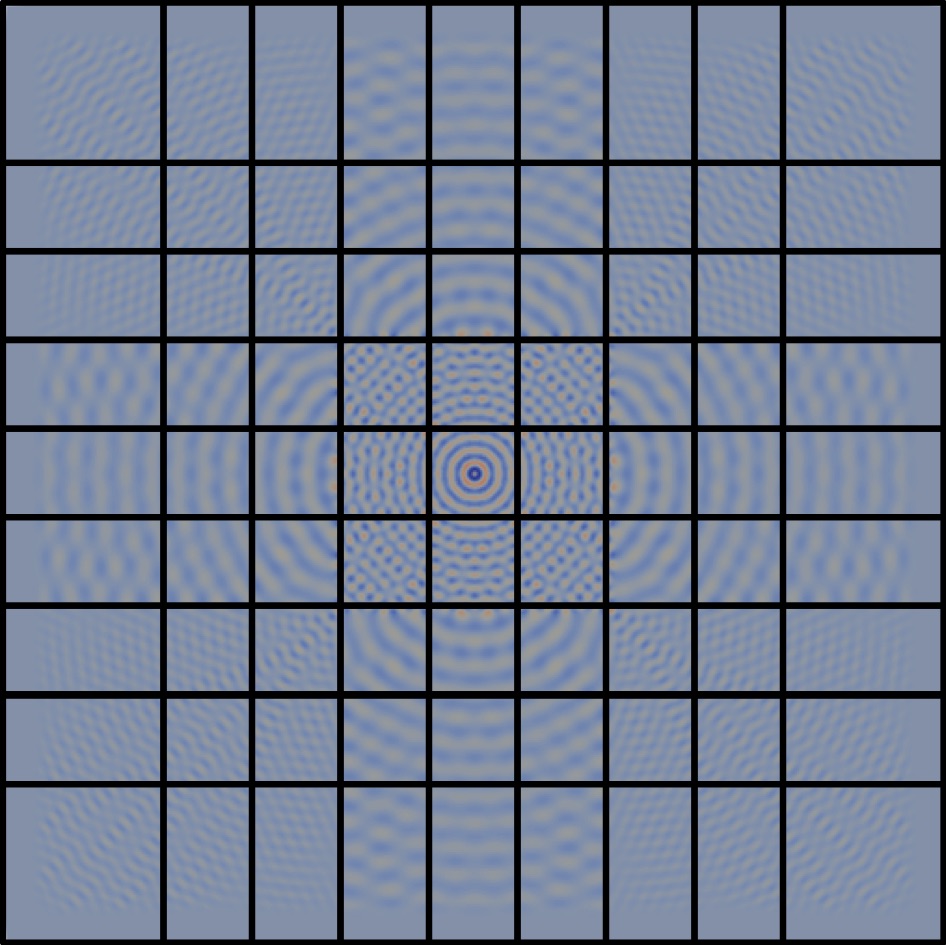}
\subcaption{Iteration 1}
\end{subfigure}\\ 
\hline
\begin{minipage}{0.45\textwidth}
\vspace*{0.5cm}
\end{minipage}&\begin{minipage}{0.45\textwidth}
\end{minipage} \\
\begin{subfigure}[b]{0.45\textwidth}
\centering\includegraphics[scale=0.2]{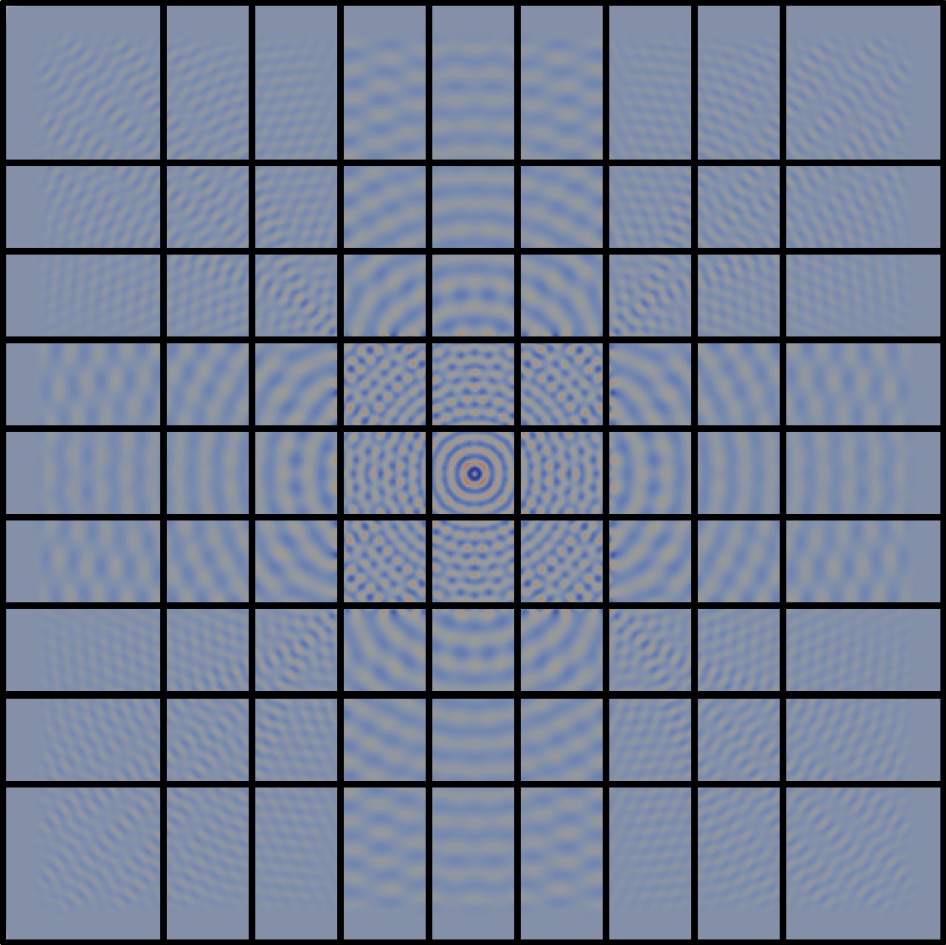}
\subcaption{Iteration 2}
\end{subfigure} & \begin{subfigure}[b]{0.45\textwidth}
\centering\includegraphics[scale=0.2]{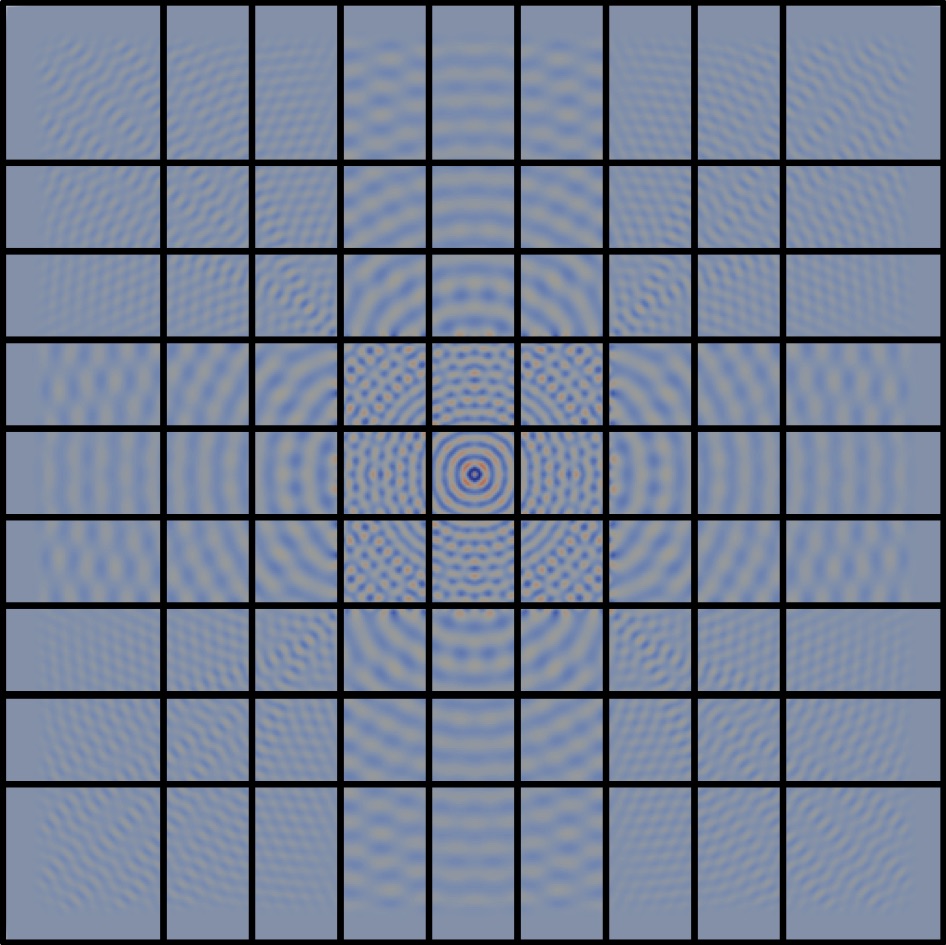}
\subcaption{Iteration 3}
\end{subfigure}
\end{tabular}
\caption{The computed wave fields after applying the preconditioner (Iteration 0) and after the first three GMRES iterations.}
\label{Fig::Sol2}
\end{figure}

To further illustrate this point, let us consider an example for a discontinuous wave speed distribution in a checkerboard pattern shown in Figure~\ref{Fig::Setup2}.
Figure~\ref{Fig::Sol2} shows the reconstructed wavefields after several GMRES iteration and the corresponding residual. Note that after 2 iterations, the procedure has already computed a reasonable wavefield. After that, almost no change is visible in the plot of the wavefield, but only in the residual. It takes 27 GMRES iterations to solve this problem to an accuracy of $10^{-6}$.

The second experiment shows that a single application of the preconditioner is insufficient to propagate all reflected waves. Thus multiple iterations are required to compute accurate global wavefields and consequently, the overall solver is sensitive to the number of reflections induced by the medium.  In fact, once the dominant reflections are resolved by the preconditioner, the number of iterations only grows as $\cO(\log \omega)$ as the frequency (and the degrees-of-freedom) grows. This is corroborated by the numerical example for a wave-guide considered in Section~\ref{Sec::NumericalExamples}.

The proposed solver is therefore scalable, however, the proportionality constant is strongly dependent on the number of reflections induced by the heterogeneous wavespeed. In fact, the numerical examples in Section~\ref{Sec::NumericalExamples} show that even in the case of wave-guides where reflections play a crucial role, the solution strategy surprisingly still only results in iteration counts lower than 40. Of course, the same arguments can be applied to refracted waves in the same way and similar effects can be seen. 

\section{Numerical Examples}
\label{Sec::NumericalExamples}
In this section, we consider several numerical examples to corroborate the claims of this paper. All numerical examples are constructed using variations of a standard setup. The problems are posed on the unit square ($d=2$) or the unit cube ($d=3$) and the the wavespeed is scaled such that the the squared slowness, $m(x) \in [m_0,1]\;\forall x\in\Omega$ for some $m_0 > 0$.   For the characteristic frequency $\omega$, the minimum wavelength is $2\pi/\omega$ and the maximum wavelength is $2\pi/(\sqrt{m_0}\omega)$. 
All problems are discretized using a uniform second-order finite difference approximation with $N$ total discretization points, i.e., there are $n=N^{\frac{1}{d}}$ points in one direction. The characteristic frequency, $\omega$, is chosen such that there are at least 10 points-per-wavelength, i.e., $\omega=2\pi n/10$.

In each case, we decompose $\Omega$ using a CDD with $q=r$ rows and columns. Each subdomain $\Omega_{ij}$ in the CDD is chosen so that, ignoring the discretization points in the PML region, the set $\boldOmegaij$ has $101$ discretization points in each direction. The skeleton of the CDD does not intersect with any discretization point, so $n=101q$ and $N=(101q)^d$.
Accordingly, the maximum frequency is $\omega=2\pi(101q-1)/10$, which is chosen so that there are 10 discretization points per wavelength.
The PML region is chosen so that $\lceil20/\sqrt{m_0}\rceil$ discretization points are in this region. This ensures that the PML region is at least 2 wavelengths thick.  The source distribution of the standard setup consists of four point sources at
\begin{alignat*}{3}
x_1&=(0.125,0.125), &\quad x_2&=(0.125,0.875),\\ 
x_3&=(0.875,0.125), &\quad x_4&=(0.875,0.875).
\end{alignat*}
These point sources are modelled by an approximation with an exponential function:
\begin{alignat*}{1}
f(x):=\frac{n^2}{\pi}\sum_{i=1}^4e^{-n^2|x_i-x|^2}.
\end{alignat*}
In all examples, the global system is solved using a preconditioned GMRES method~\cite{Saad_Schultz:GMRES} with a tolerance of $10^{-6}$, using the zero-vector as an initial guess. All local problems are solved using Pardiso 6.0~\cite{Kourounis_Schenk:Pardiso}.  Source code, models, and experiment configurations for these examples are available online \cite{LSweeps2DCode,LSweeps2DExamples}.

\subsection{Effect of PML-induced and discretization errrors}
\label{Sec::NumericalExamples::Effectiveness}

\begin{table}
    \begin{center}
        \begin{tabular}{|c|c|c|r|r|r|r|r|r|r|}
            \hline
             & & & \multicolumn{7}{|c|}{wavelengths}\\
           $N$ & & & \multicolumn{7}{|c|}{in PML region}\\   
            (without PML) & $\omega/2\pi $& $q=r$ &  1 & 2 & 3 & 4 & 5 & 6  & 7\\
            \hline
            $202 \times 202$      &    20.1 & 2   & 3  & 2 & 1 & 1 & 1 & 1 & 1  \\
            $404 \times 404 $     &    40.3 & 4   & 4  & 3 & 1 & 1 & 1 & 1 & 1 \\
            $808 \times 808 $     &    80.7 & 8   & 5  & 3 & 2 & 1 & 1 & 1 & 1 \\
            $1616 \times 1616 $   &   161.5 & 16  & 6  & 2 & 2 & 1 & 1 & 1 & 1  \\
            $3232 \times 3232 $   &   323.1 & 32   & 7  & 3 & 2 & 2 & 1 & 1 & 1  \\
            $6464 \times 6464 $   &   646.3 & 64  & 8  & 4 & 3 & 3 & 2 & 1 & 1 \\
            $12928 \times 12928 $ &  1292.7 & 128 & 10 & 5 & 3 & 3 & 3 & 2 & 1 \\
            \hline
        \end{tabular}
        \caption{Number of iterations necessary to solve a homogeneous problem using increasingly thick PMLs and different frequencies.}
\label{Tab::PMLeffect}        
    \end{center}
\end{table}

In this section we provide numerical evidence of the claim in Section~\ref{Sec::LSweeps} that, in a homogeneous medium, the accuracy of the preconditioner depends only on the accuracy of the discretization and the quality of the absorbing boundary conditions.
We consider the standard setup for a constant squared slowness $m=1$ and we control the accuracy of the absorbing boundary condition by increasing (or reducing) the PML thickness. 
Table~\ref{Tab::PMLeffect} depicts the dependence of the preconditioner, measured in iterations, on the quality of the absorbing boundary condition, which is tuned by varying the number of wavelengths inside the PML region between 1 and 7. If the accuracy of the absorbing boundary condition is reduced (fewer wavelengths), the preconditioner is less effective. If the accuracy of the absorbing boundary condition is increased (more wavelengths), the preconditioner is more effective. In fact, if the accuracy of the boundary condition is sufficiently high, the preconditioner is perfectly effective resulting in iteration counts independent of $\omega$. It can also be seen that the PML thickness has to be increased with $\omega$ in order to achieve this perfect effectiveness.

Nevertheless, if the thickness of the PML region is kept constant with the problem size, the number of required GMRES iterations to solve the global system only grows logarithmically with the wave number $\omega$. It is surprising that this qualitative property holds independently of the thickness of the PML, even for very thin PML regions, for example one wavelength. 

This shows that the thickness of the PML has a significant influence on the effectiveness of the preconditioner. On the other hand, the thickness of the PML region is also crucial for the scalability of the preconditioner. For example, if the PML region is chosen to contain $n_{\texttt{pml}}$ discretization points with $n_{\texttt{pml}}\ll n$, the sizes of the local problems associated with each subdomain grows as ${\cal O}(n_{\texttt{pml}}^2)$ in 2D and as $\cO((n+n_{\texttt{pml}})n_{\texttt{pml}}^2)$ in 3D. Consequently, the parallel factorization of all subdomains scales as ${\cal O}(n^2n_{\texttt{pml}}^3/p)$ in 2D and ${\cal O}(n^3n_{\texttt{pml}}^6/p)$ in 3D and the application of the L-sweep preconditioner scales in parallel as ${\cal O}(n^2n_{\texttt{pml}}^2\log n_{\texttt{pml}}/p)$ in  2D and ${\cal O}(n^3 n_{\texttt{pml}}^4/p)$ in 3D. Therefore, it is crucial to keep the thickness of the PML region as thin as possible to preserve efficiency. This makes the trade-off between efficiency and effectiveness of the preconditioner apparent. For all subsequent numerical examples considered in this paper, we choose a PML thickness of 2 wavelengths, which is empirically sufficient to achieve satisfactory results. In particular, the logarithmic growth of the number of iterations with respect to $\omega$ with this choice of PML-thickness allows one to achieve the advocated parallel scaling of the solver.

We also performed the same experiment for finite difference discretizations with fewer than 10 discretization points per wavelength. For very low accuracy discretizations  (e.g., four points per wavelengths and more than 300 wavelengths inside the domain), the preconditioner can fail, i.e., the number of iterations to solver the preconditioned linear system grows faster than $\cO(\log\omega)$. This is because the preconditioner relies on physical properties of wave propagation which are not captured by such inaccurate discretizations. 

Considering these experiments, we conjecture that if the thickness of the PML region is kept constant with the frequency $\omega$ and the discretization is sufficiently accurate, the number of iterations required to solve the problem using a GMRES method grows as $O(\log \omega)$. In addition, if the transparent boundary conditions can be modelled perfectly, the iteration count can be reduced to $O(1)$ for reasonably smooth media.

\subsection{Complexity of the solver}
\label{Sec::NumericalExamples::Complexity}
In this section we provide empirical run-times using the proposed methods to support our scalibility claims.

In this experiment, we consider the standard setup with constant wave speed, for $q=2,4,8,\dots,128$ and we measure the factorization time, the time spent in one GMRES iteration, and the total time of the solver, which are summarized in Table~\ref{Tab::Complexity}. All timings are recorded on the NERSC machine 'Cori'. Cori is composed of 2388 Haswell compute nodes, each node containing two 16-core Intel Xeon Processor E5-2698 v3 running at 2.3 GHz and 128 GB of memory. The nodes communicate via a Cray Aries interconnect with Dragonfly topology. The number of MPI ranks are chosen so that each compute node is assigned one MPI rank and each MPI rank is assigned one row of subdomains in the CDD. The number of MPI ranks is therefore $p=q$.

\begin{figure}
    \centering
    \includegraphics[trim = 0mm 0mm 0mm 0mm, clip,width=0.7\textwidth]{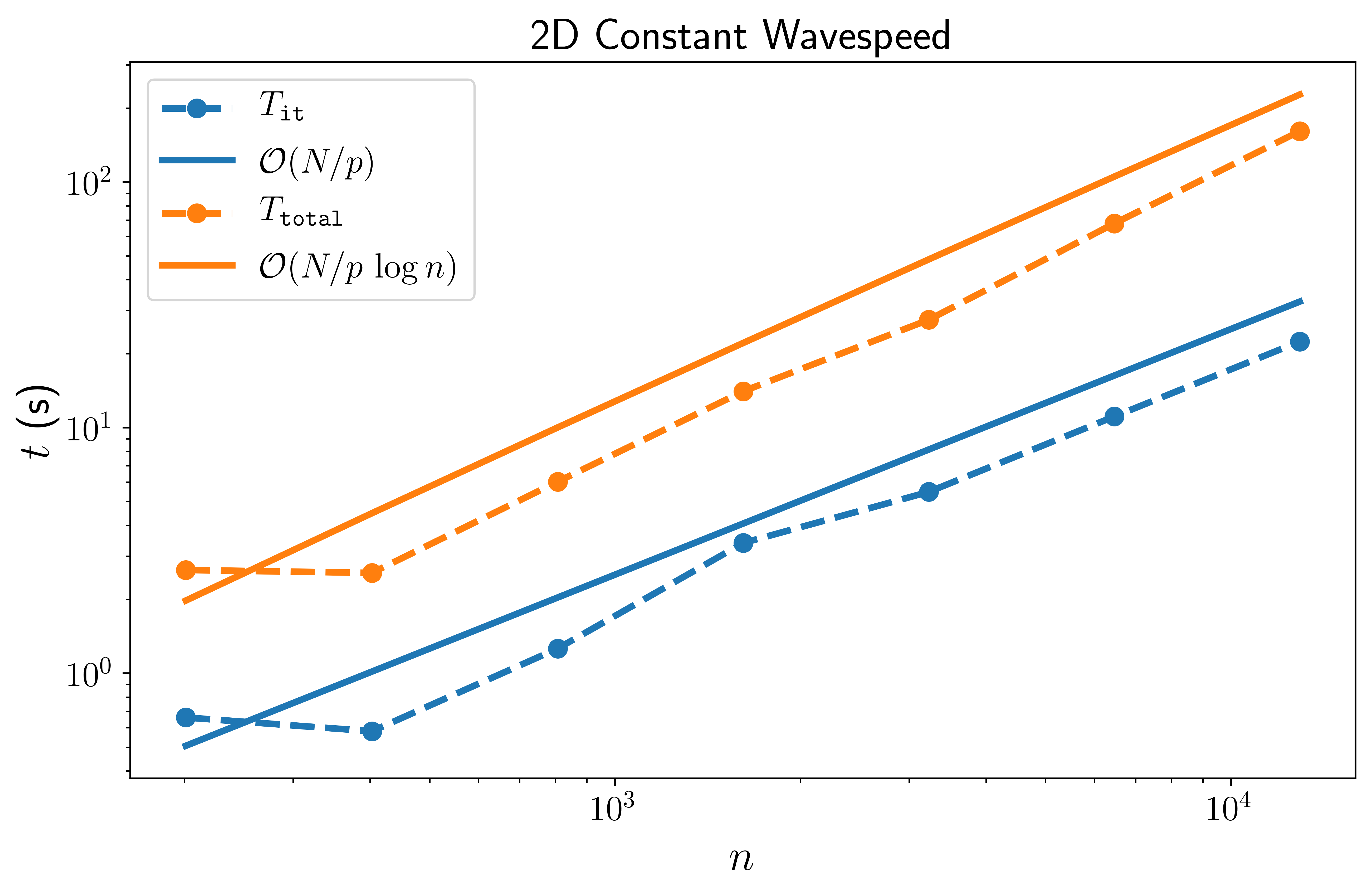}
    \caption{Timings for $2$D experiment with constant wavespeed. Solid lines are empirical complexities for one application of the preconditioner (blue) and for the over-all solve (orange).  Dashed lines are with measured results for one application of the preconditioner (blue) and for the over-all solve (orange).
    }
    \label{Fig::Scaling2D}
\end{figure}

\begin{table}
    \begin{center}
        \begin{tabular}{|c|r|r|r|r|r|r|r|}
            \hline
            $N$ & $\omega/2\pi $& $p$ & $T_{\texttt{fact}}$ & $N_{\texttt{it}}$ & $T_{\texttt{it}}$ & $T_{\texttt{total}}$\\
            \hline
              $202 \times   202$    &      20.1 & 2   &  1.09 & 2  &  0.66  &    2.63  \\
              $404 \times   404 $   &      40.3 & 4   &  1.00 & 3  &  0.58  &    2.56  \\
              $808 \times   808 $   &      80.7 & 8   &  1.41 & 3  &  1.26  &    6.02  \\
             $1616 \times  1616 $   &     161.5 & 16  &  2.80 & 2  &  3.39  &   14.05  \\
             $3232 \times  3232 $   &     323.1 & 32  &  4.41 & 3  &  5.47  &   27.47  \\
             $6464 \times  6464 $   &     646.3 & 64  &  8.34 & 4  & 11.09  &   67.74  \\
            $12928 \times 12928 $   &    1292.7 & 128 & 15.66 & 5  & 22.39  &  160.88  \\
            \hline
        \end{tabular}
        \caption{ Timings for the solver on a homogeneous wave speed. $N$ is the number of degrees-of-freedom, $\omega$ is the frequency, $T_{\texttt{fact}}$ is the time spent to factorize the local problems, $N_{\texttt{it}}$ is the number of GMRES iterations needed to solve the problem, $T_{\texttt{it}}$ is the average time spent per GMRES iteration, and $T_{\texttt{total}}$ is the total time spent to solve the linear system. All timings are measured in seconds. }
        \label{Tab::Complexity}
    \end{center}
\end{table}
Table~\ref{Tab::Complexity} contains the results of this experiment, from which we observe the scalability of the preconditioner.  The timings clearly reflect the ${\cal O}(N/p)$ scaling as claimed in the prequel. Furthermore, the number of iterations $N_{\texttt{it}}$ grows as ${\cal O}(\log \omega)$ which further proves the effectiveness of the preconditioner. Combining both leads to the total time, $T_{\texttt{total}}$, for the solver which reflects the claimed $\cO\big((N/p)\log\omega\big)$ scaling, as illustrated in Figure~\ref{Fig::Scaling2D}.

\subsection{Smooth wave speeds}
\label{Sec::NumericalExamples::SmoothWavespeed}
\begin{figure}[htp]
\centering
\includegraphics[scale=0.3175]{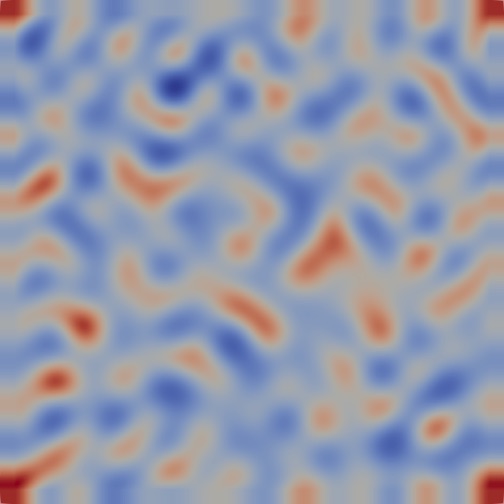}\hspace{2cm}\includegraphics[scale=0.3]{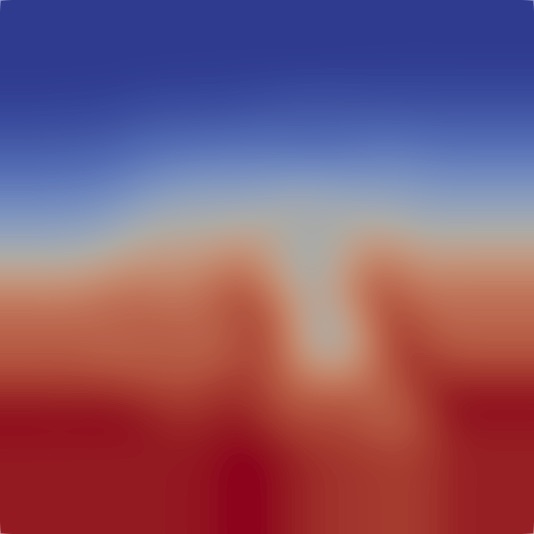}
\caption{The smooth wave speed distributions. Left: Random smooth wave speed; Right: Smooth BP model background}
\label{Fig::SmoothWavespeeds}
\end{figure}

In this section, we demonstrate the impact of \textit{refracted} waves, induced by smooth media, on the performance of the preconditioner and iterative solver.
As established in Section~\ref{Sec::LSweeps::Heterogeneities}, the effectivness of the L-sweeps preconditioner strongly depends on the presence of reflections and refractions in the solution.
To evaluate the impact of refracted waves, we consider the standard setup for two different smooth wave speeds, which are shown in Figure~\ref{Fig::SmoothWavespeeds}. One is a random smooth wavefield, the other is a typical background wave speed used in seismic imaging. The latter is obtained by smoothing the background of the BP model~\cite{BP_model}.
 
Table~\ref{Tab::RandomWavespeed} shows the iteration counts obtained for both wave speeds.
\begin{table}
    \begin{center}
        \begin{tabular}{|c|c|c|r|r|}
        \hline
           $N$ & & & Random smooth & Smooth BP model\\   
            (without PML) & $\omega/2\pi $& $q=r$ &  wave speed  & background\\
            \hline
            $202   \times 202$    &    20.1 & 2   & 3  & 1\\
            $404   \times 404 $   &    40.3 & 4   & 3  & 2\\
            $808   \times 808 $   &    80.7 & 8   & 5  & 4\\
            $1616  \times 1616 $  &   161.5 & 16  & 7  & 5\\
            $3232  \times 3232 $  &   323.1 & 32  & 9 & 6\\
            $6464  \times 6464 $  &   646.3 & 64  & 10 & 7\\
            $12928 \times 12928 $ &  1292.7 & 128 & 12 & 8 \\    
            \hline
        \end{tabular}
        \caption{Iteration counts for smooth wave speeds for different frequencies.}
\label{Tab::RandomWavespeed}        
    \end{center}
\end{table}
Table~\ref{Tab::RandomWavespeed} shows the effectiveness of the preconditioner when solving problems involving a smooth wave speed. In particular, we observe that the number of iterations scales as ${\cal O}(\log\omega)$, albeit with higher constants compared to the constant case.

\subsection{BP model}
\label{NumericalExamples::BP}
In this section, we consider the wave speed of a standard geophysical benchmark problem, a subset of the BP model~\cite{BP_model}. It is well-known that this problem involves many reflections in the resulting wavefields. We therefore use it to study the effect of reflections on the L-sweeps preconditioner.

\begin{figure}
    \centering
\begin{tabular}{cc}
\begin{subfigure}[b]{0.45\textwidth}
    \centering
\includegraphics[scale=0.3]{BPsmoothBackground.jpg}
\caption{The smooth background model (BG1).}
\label{Fig::WaveSpeedBP::BG1}
\end{subfigure} & 
\begin{subfigure}[b]{0.45\textwidth}
    \centering
\includegraphics[scale=0.298125]{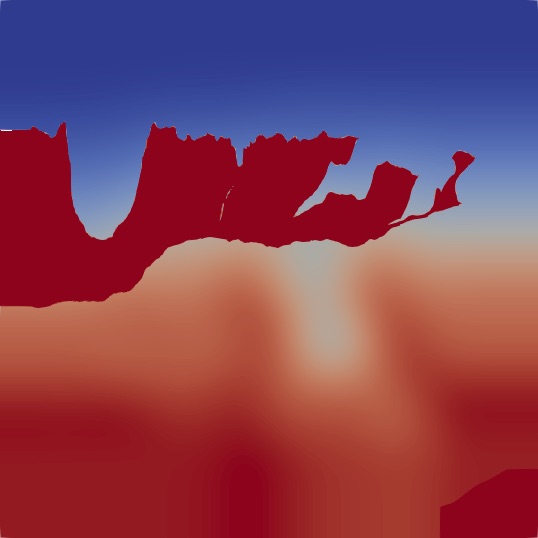}
\caption{BG1, with salt body.}
\label{Fig::WaveSpeedBP::BG1Salt}
\end{subfigure}\\
\vspace{0.75cm}\\
\begin{subfigure}[b]{0.45\textwidth}
    \centering
\includegraphics[scale=0.3]{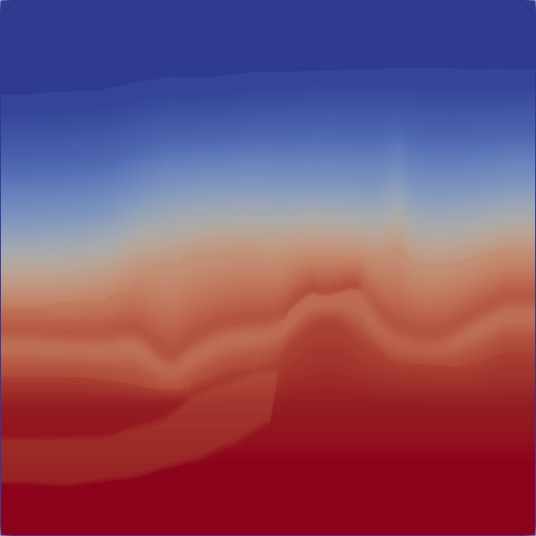}
\caption{The `true' BP model background velocity (BG2).}
\label{Fig::WaveSpeedBP::BG2}
\end{subfigure} & \begin{subfigure}[b]{0.45\textwidth}
    \centering
\includegraphics[scale=0.298125]{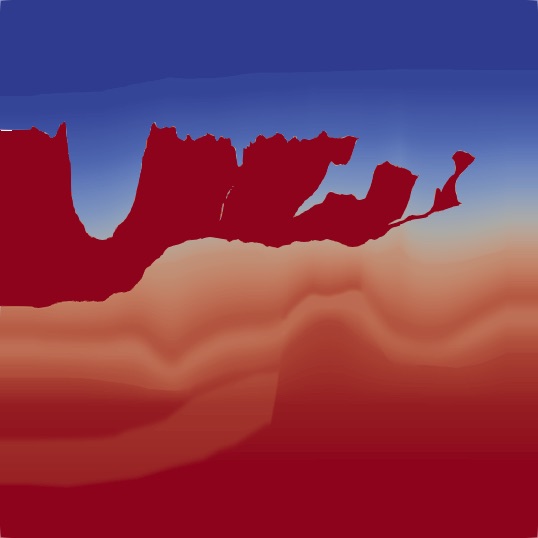}
\caption{BG2 with salt body.}
\label{Fig::WaveSpeedBP::BG2Salt}
\end{subfigure}\\
\end{tabular}
\vspace{0.75cm}\\
\begin{subfigure}[b]{0.45\textwidth}
    \centering
\includegraphics[scale=0.298125]{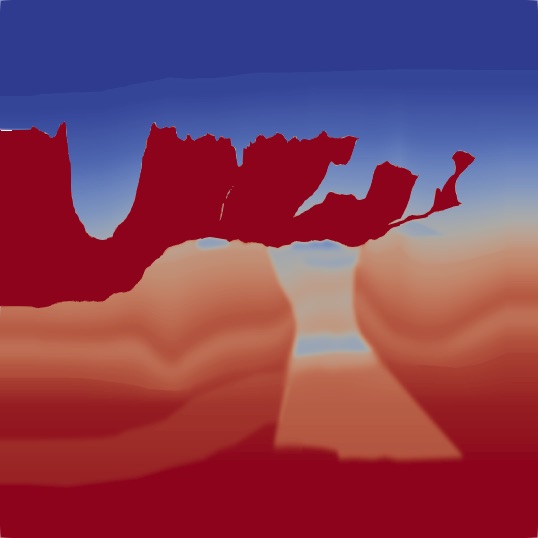}
\caption{The `true' BP model.}
\label{Fig::WaveSpeedBP::BPTrue}
\end{subfigure}
    \caption{Wave speed distributions for the study of the BP model.}
    \label{Fig::WaveSpeedsBP}
\end{figure}

First, we establish that reflected waves arise primarily due to the salt body, the dark red region in Figures~\ref{Fig::WaveSpeedBP::BG1Salt}, \ref{Fig::WaveSpeedBP::BG2Salt}, and \ref{Fig::WaveSpeedBP::BPTrue}. To this end, we first consider two smooth background wave speeds, the one from Section~\ref{Sec::NumericalExamples::SmoothWavespeed}, which we will call BG1 and illustrate in Figure~\ref{Fig::WaveSpeedBP::BG1} and the true background velocity provided in~\cite{BP_model}, which we call BG2 and illustrate in Figure~\ref{Fig::WaveSpeedBP::BG2}. Then, we superimpose the salt body on BG1 and BG2, building to examine the effect of the high contrasts due to the salt body on the performance of the solver. Finally, for completeness, we also run the test on the true BP model, as obtained from~\cite{BP_model}, illustrated in Figure~\ref{Fig::WaveSpeedBP::BPTrue}.

\begin{figure}
     \hspace*{-0.5cm}
    \includegraphics[scale=0.1075]{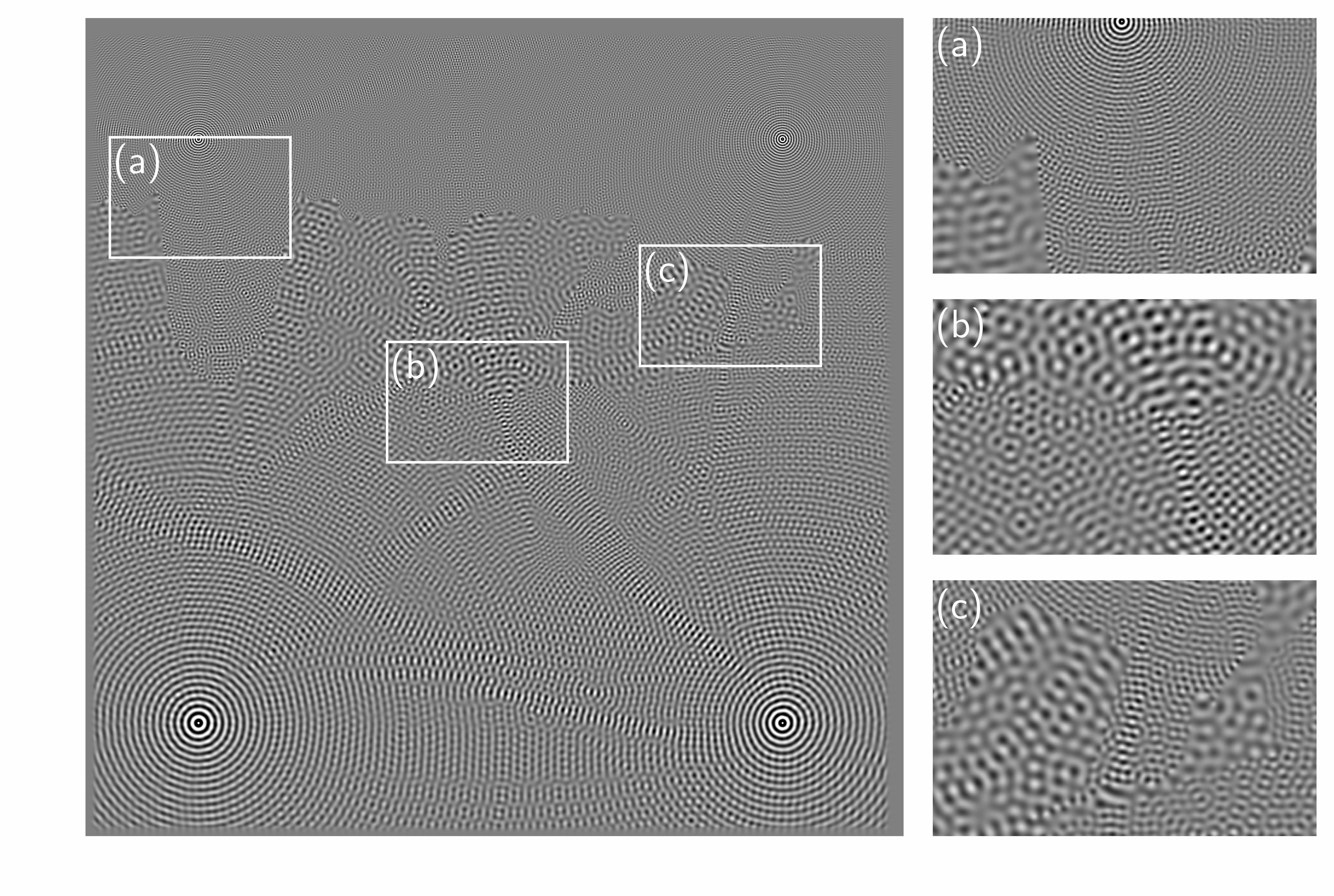}
    \caption{The computed wave field for the BP model and $\omega/2\pi=323.1$.  Zoom regions (a), (b), and (c) highlight solution quality in areas of high-contrast in the wavespeed (at the salt boundary).}
    \label{Fig::BPWavefield}
\end{figure}

Table~\ref{Tab::BP} shows the iteration counts observed for each of the 5 experiments.
\begin{table}
    \begin{center}
        \begin{tabular}{|c|c|c|r|r|r|r|r|}
        \hline
           $N$ & & &  & & \multicolumn{1}{|c|}{BG1} & \multicolumn{1}{|c|}{BG2} & \multicolumn{1}{|c|}{BP}\\   
            (without PML) & $\omega/2\pi $& $q=r$ &  BG1  & BG2 & with salt & with salt & model\\
            \hline
            $202 \times 202$     &   20.1 & 2   &  1 &  4 &  7 &  6 &  7\\
            $404 \times 404$     &   40.3 & 4   &  2 &  4 &  9 &  9 &  9\\
            $808 \times 808$     &   80.7 & 8   &  4 &  6 & 12 & 12 & 12\\
            $1616 \times 1616$   &  161.4 & 16  &  5 &  6 & 15 & 15 & 15\\
            $3232 \times 3232$   &  323.1 & 32  &  6 &  7 & 17 & 17 & 16\\
            $6464 \times 6464$   &  646.3 & 64  &  7 &  7 & 19 & 19 & 19\\
            $12928 \times 12928$ & 1292.7 & 128 &  8 &  8 & 21 & 21 & 20\\            
            \hline
        \end{tabular}
        \caption{Iteration counts for the BP model.}
\label{Tab::BP}        
    \end{center}
\end{table}
For both BG1 and BG2, with no salt, we observe similar iteration counts and conclude that  reflections do not play a big role in the resulting wavefields. However, after superimposing salt body on top of of BG1 and BG2, we clearly see that the number of iterations increase significantly. This is explained because reflections are introduced in the wavefield due to the high contrast between the salt and the background.  We also observe that the exact BP model can be solved in almost the same iteration counts as BG1 and BG2 with the salt body superimposed.  Thus, even for the exact BP model, the reflections induced by the salt body dominate the performance of the preconditioner. Finally, let us note that even though reflections in the wavefield result in higher iteration counts, with increasing frequency $\omega$, the number of iterations always grows only as ${\cal O}(\log\omega)$ and they are never higher than 21 iterations, even for problems involving 1000 wavelengths inside the domain. This shows the effectiveness of the preconditioner, even in the presence of reflections as they appear in practically relevant problems. The wave field computed for the BP model and $\omega/2\pi=323.1$ is shown in Figure \ref{Fig::BPWavefield}.

\subsection{Wave-guide}
\label{NumericalExamples::WaveGuide}
\begin{figure}[htp]
\centering
\includegraphics[scale=0.4]{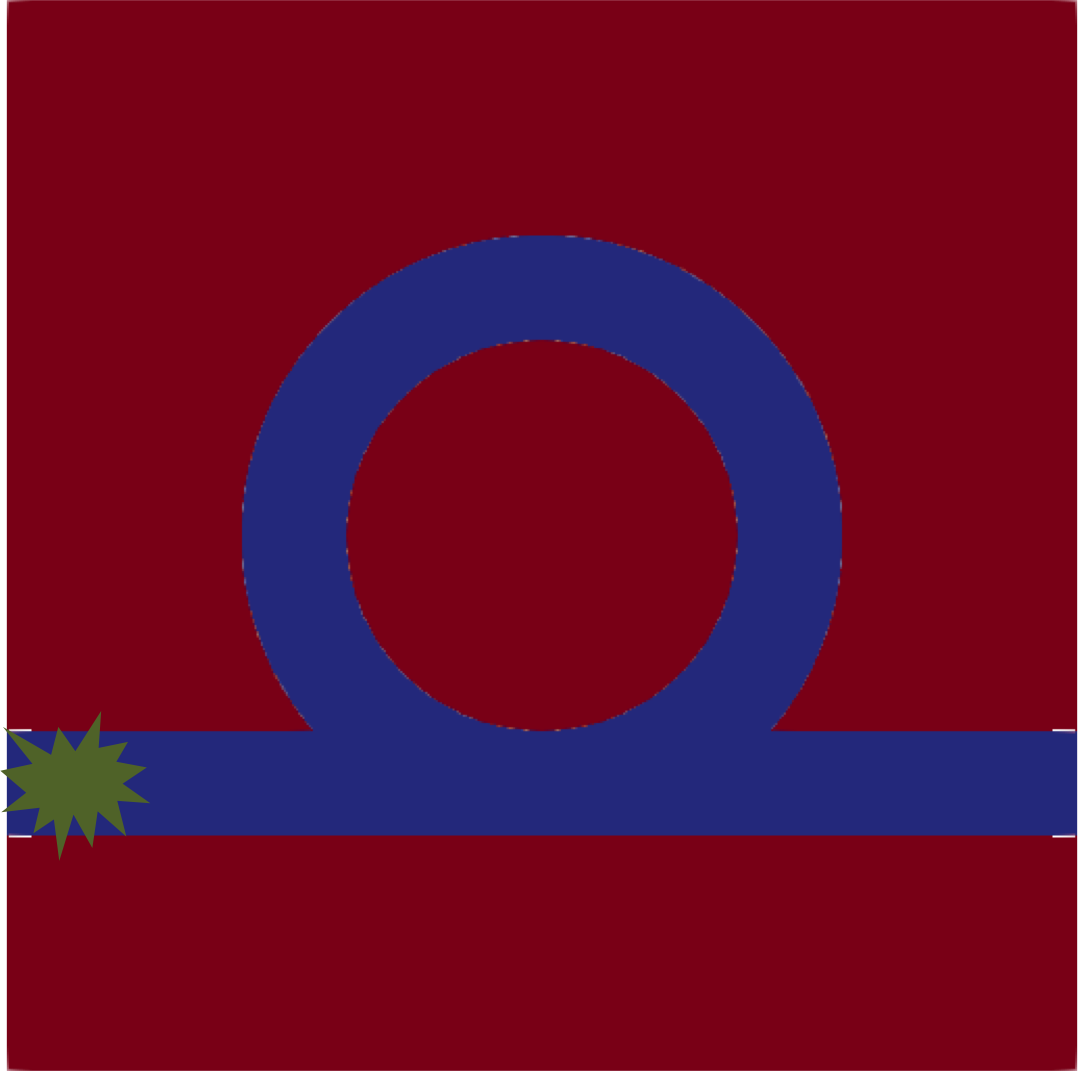}
\caption{The wave speed distribution of the wave-guide. The point source at the entrance is shown by the green star.} 
\label{Fig::WaveGuide}
\end{figure}

To conclude the studies of two-dimensional problems, we consider a 
wave-guide with a point source at the entrance, a problem that is designed to stress the capabilities of the preconditioner because solutions contain  many reflections.  The wave-guide is illustrated in Figure~\ref{Fig::WaveGuide}. Other than the change in the source and velocity configuration, the experimental setup follows the standard setup.

\begin{figure}
    \centering
     \hspace*{-0.5cm}
    \includegraphics[scale=0.1075]{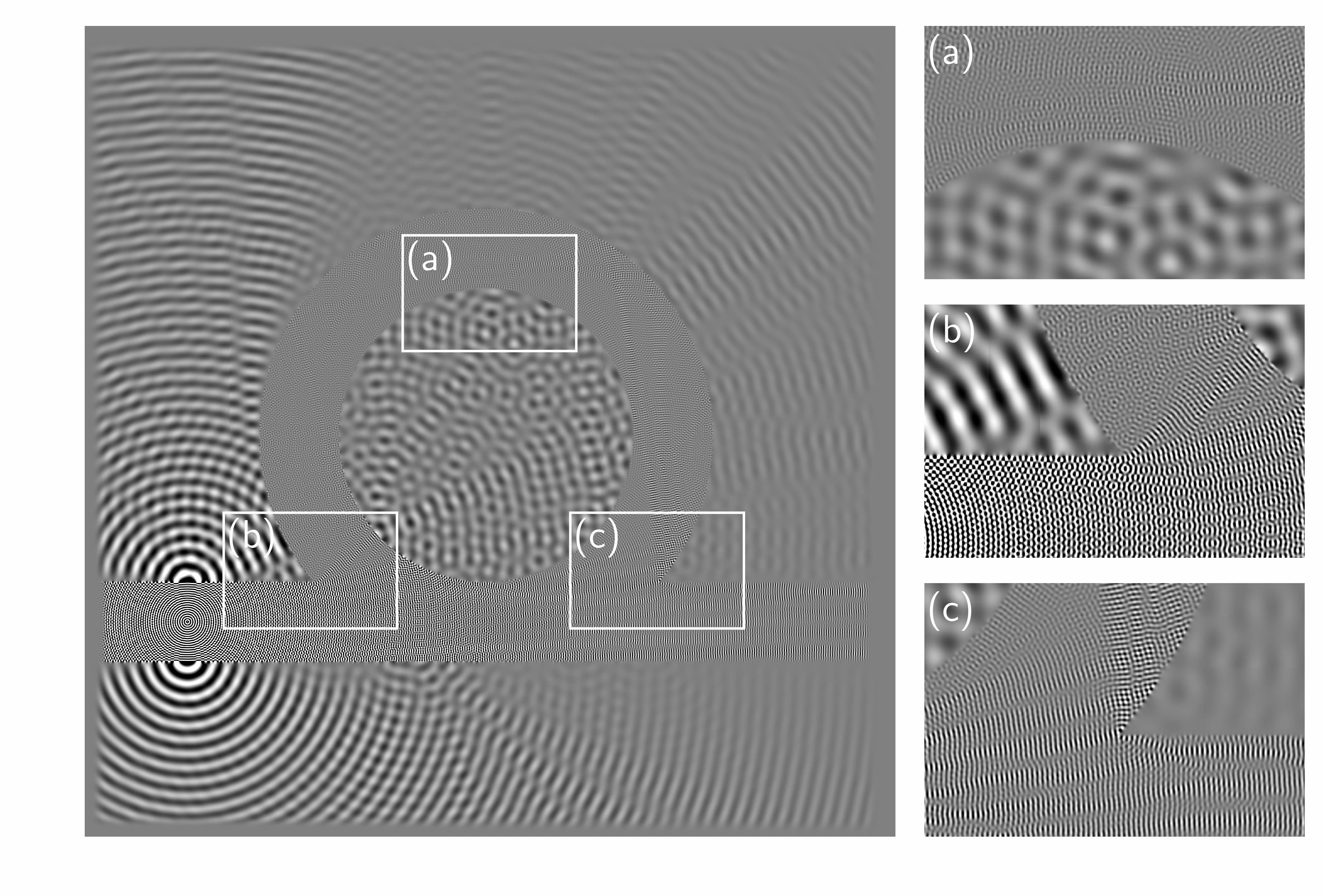}
    \caption{The computed wave field for a contrast ratio of 6 and $\omega/2\pi=323.1$.  Zoom regions (a), (b), and (c) highlight solution quality in areas of high-contrast in the wavespeed.}
    \label{Fig::WaveguideWavefield}
\end{figure}

To stress the preconditioner, we vary the contrast ratio between the background and wave-guide.
Table~\ref{Tab::Waveguide} shows the results of the experiment.
\begin{table}
    \begin{center}
        \begin{tabular}{|c|c|c|c|r|r|r|r|r|r|}
        \hline
           $N$ & & & \multicolumn{5}{|c|}{Contrast ratio} \\
            (without PML) & $\omega/2\pi $& $m=n$ &  2 & 3 & 4 & 5 & 6 \\
            \hline
            $202 \times 202$    &   20.1  & 2   & 18 & 24 & 24 & 25 & 26 \\
            $404 \times 404 $   &   40.3  & 4   & 28 & 29 & 29 & 28 & 30  \\
            $808 \times 808 $   &   80.7  & 8   & 30 & 32 & 34 & 33 & 33  \\
            $1616 \times 1616 $ &   161.5 & 16  & 31 & 33 & 33 & 34 & 35  \\
            $3232 \times 3232 $ &   323.1 & 32  & 32 & 34 & 36 & 36 & 37  \\
            $6464 \times 6464 $ &   646.3 & 64  & 32 & 34 & 35 & 36 & 36  \\
            \hline
        \end{tabular}
        \caption{Iteration counts for wave-guide problem for different contrast ratios of the wave speed between the wave-guide and the background.}
\label{Tab::Waveguide}        
    \end{center}
\end{table}
Once all reflections are resolved for each problem, the number of iterations grows logarithmically with the frequency $\omega$. In fact, the growth is so slow that it is comparable to the problem involving constant wave speeds. In addition, while the number of iterations is clearly higher for this problem, it is interesting to observe that all problems can still be solved in under 40 iterations. It is noteworthy that the number of iterations is apparently independent the contrast ratio.  The wave field computed for a contrast ratio of 6 and $\omega/2\pi=323.1$ is shown in Figure \ref{Fig::WaveguideWavefield}.

%

\subsection{3D Example}
As a final example, we consider a three-dimensional problem. The effectiveness of the preconditioner has been thoroughly established in the previous two-dimensional numerical examples. These results directly translate over to three-dimensional problems. Therefore, we only consider the scalability of the solver for a three-dimensional problem for a constant wave speed.

To this end, we consider the standard setup in a 3D setting. This means that instead of constantly sized two-dimensional problems, the local problems in each subdomain are beam-shaped three-dimensional problems. These beam shaped local problems are obtained from the local 2D problems by extending them in the third dimension such that the global problem has $n$ degrees of freedom in each direction where, as before, $q=r$ is the number of subdomains in one direction of the CDD. The only other differences in problem setup are that we consider a discretization of $6$ points per wavelength, assume that each subdomain contain 2 wavelengths, and choose a PML thickness of only one wavelength. Not counting the PML region, this results in problems with $n=13q$ degrees of freedom in one direction. All of these changes are for computational expedience.

The timings are measured in the same computational environment as in Section~\ref{Sec::NumericalExamples::Complexity}. As before, one row of subdomains is assigned to one node, i.e., $p=q$. Table~\ref{Tab::Complexity_3D} shows the resulting time $T_{\texttt{fact}}$ needed for the factorization of all the local problems, the number of iterations $N_{\texttt{it}}$ required in the GMRES method, the average time $T_{\texttt{it}}$ required in each iteration, and the total time $T_{\texttt{total}}$ required for the solver.

\begin{table}
    \begin{center}
        \begin{tabular}{|c|r|r|r|r|r|r|}
            \hline
            $N$ & & & & & & \\
            (without PML) & $\omega/2\pi $& $p$ & $T_{\texttt{fact}}$ & $N_{\texttt{it}}$ & $T_{\texttt{it}}$ & $T_{\texttt{total}}$\\
            \hline
            $26\times 26\times 26$      &      4.17  & 2   &   .04 & 4  &  1.34  &    6.52  \\
            $52\times 52\times 52$      &      8.50  & 4   &  5.54 & 6  &  5.30  &   37.17  \\
            $78\times 78\times 78$      &      12.83 & 6   & 12.42 & 6  & 12.80  &   89.76  \\
            $104\times 104\times 104$   &      17.17 & 8   & 22.91 & 6  &  23.27  &   163.62 \\
            $13-\times 130\times 130$   &      21.50 & 10  &  37.53 & 7  & 36.47  &    292.33  \\
            $156\times 156\times 156$   &     25.83 & 12  &  52.47 & 7  & 51.62  &    417.08  \\
            $182\times 182\times 182$   &     30.17 & 14  &  71.71 & 8  & 68.92  &   627.23  \\
            $208\times 208\times 208$   &     34.50  & 16  &  96.14 & 7  & 91.65 &   743.37 \\
            $234\times 234\times 234$   &     38.83  & 18  &  124.64 & 8  & 116.08 &   1050.31 \\
            $260\times 260\times 260$   &     43.17 & 20  &  211.87 & 7  & 177.21 &   1438.12 \\
            $312\times 312\times 312$   &     51.83  & 24  &  314.93 & 8  &  263.16 &   2457.40 \\
            $416\times 416\times 416$   &     69.17  & 32  &  418.36 & 9  &  377.60 &   3992.63 \\
            \hline
        \end{tabular}\caption{ Timings for the solver on a $3$D homogeneous wave speed. $N$ is the number of degrees-of-freedom, $\omega$ is the frequency, $T_{\texttt{fact}}$ is the time spent to factorize the local problems, $N_{\texttt{it}}$ is the number of GMRES iterations needed to achieve a relative residual of $10^{-7}$, $T_{\texttt{it}}$ is the average time spent per GMRES iteration, and $T_{\texttt{total}}$ is the total time spent to solve the linear system. All timings are measured in seconds. }
        \label{Tab::Complexity_3D}
    \end{center}
\end{table}

\begin{figure}
    \centering
    \includegraphics[trim = 0mm 0mm 0mm 0mm, clip,width=0.7\textwidth]{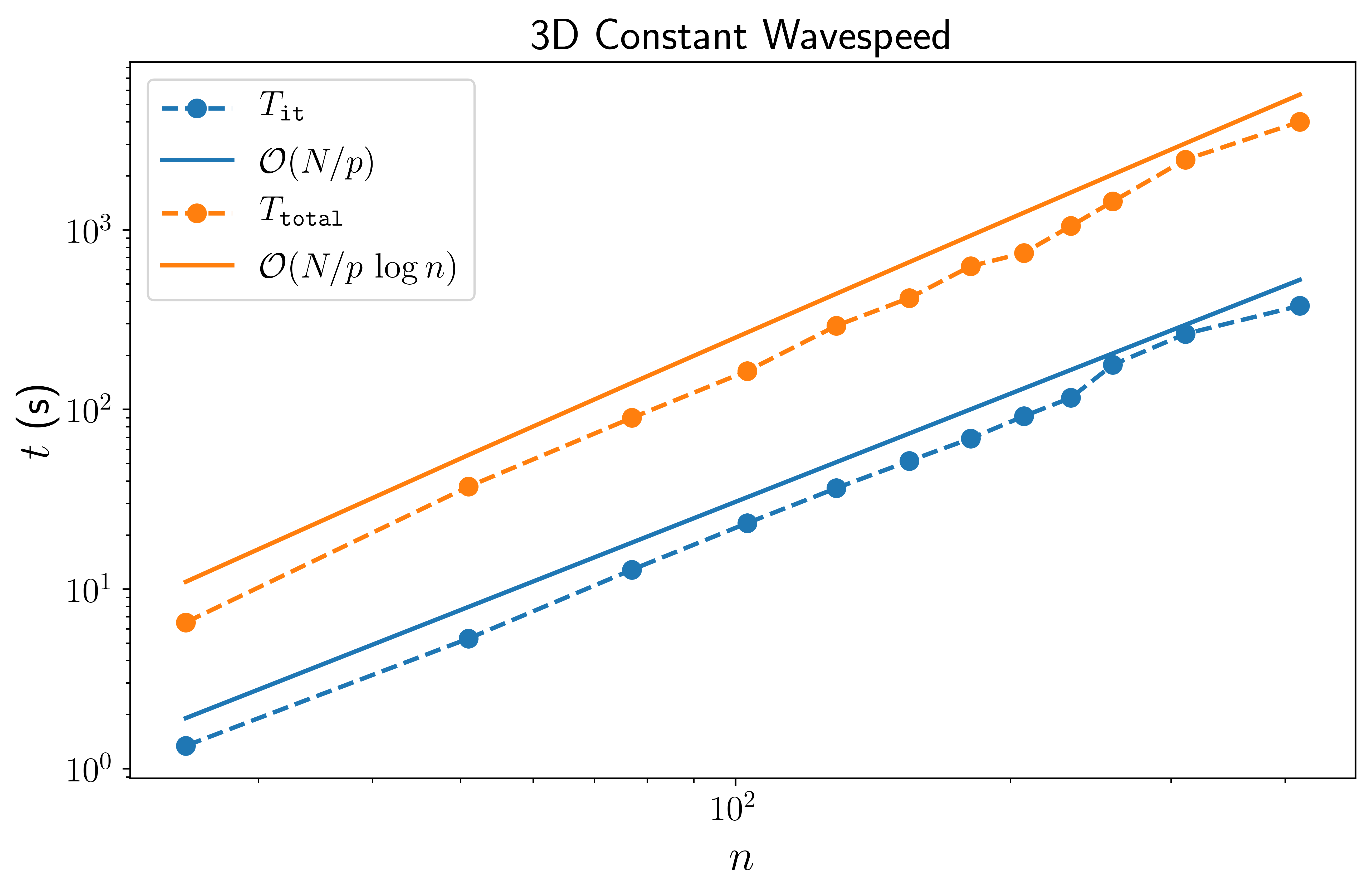}
    \caption{Timings for the $3$D experiment with constant wavespeed. Solid lines are empirical complexities for one application of the preconditioner (blue) and for the over-all solve (orange).  Dashed lines are with measured results for one application of the preconditioner (blue) and for the over-all solve (orange).}
    \label{Fig::Scaling3D}
\end{figure}

It can be clearly seen from Table~\ref{table:complexity3D} that $T_{\texttt{fact}}$ and $T_{\texttt{it}}$ both show the claimed $\cO(N/p)$ scaling. In addition, as in the two-dimensional cases, the number of iterations $N_{\texttt{it}}$ grows as $\cO(\log\omega)$ resulting in the almost optimal $\cO\big((N/p)\log\omega\big)$ scaling for the total time $T_{\texttt{total}}$ needed for the solver as shown in  Figure~\ref{Fig::Scaling3D}. This corroborates the claimed complexity for three-dimensional problems.

\section{Discussion}
\label{Sec::Discussion}
We have introduced the first solver for the high-frequency Helmholtz equation that scales as $\cO\big((N/p)\log \omega\big)$ a distributed-memory parallel computational environment, for a single right-hand side.
The new solver constructs a preconditioner based on a CDD. This approach reveals parallelism which is exploited to obtain the optimal parallel complexity. The performance of the preconditioner is similar to the well-established method of polarized traces, i.e., the number of iterations grows as $\cO(\log\omega)$.

We have shown that the preconditioner is effective once all reflections and refractions in the wavefield are resolved.
While this is feasible for many practical applications, the performance of the proposed solver deteriorates in the presence of excessive reflected or refracted waves.  This limitation arises, for example, on problems containing resonant cavities.
We do not think that sweeping strategies are the correct approach to solve these problems, and due to the relevance of these problems, we view this case as a topic for further research.

Future work of the current method includes improvements for three-dimensional problems. As it is presented in this paper, the method results in a $\cO(n^2)=\cO(N^{\frac{2}{3}})$ parallel complexity for three-dimensional problems. However, all beam-shaped local problems are treated by one MPI rank only. It is well established that these quasi-one-dimensional problems can be treated in a distributed memory parallel computational environment by the use of multi-frontal methods \cite{Rouet_Li_Ghysels:A_distributed-memory_package_for_dense_Hierarchically_Semi-Separable_matrix_computations_using_randomization,Amestoy_Duff:MUMPS,li_demmel03:SuperLU_DIST}. In this case, further parallelism can be exploited and an almost  $\cO(n)=\cO(N^{\frac{1}{3}})$ parallel complexity can be achieved for three-dimensional problems.

Using this extension, it can be shown that the high-frequency Helmholtz equation can be optimally parallelized in all spatial directions but one. The remaining direction is due to the inherently serial nature of the sweeps. In this regard, we conjecture that this work achieves the boundary of the scalability of sweeping preconditioners for a single right-hand side and that further improvements require new ideas. Nevertheless, in the presence of $\cO(n)$ right-hand sides, sweeping preconditioners can be further parallelized by pipelining over the right-hand sides. In that regard, the parallel complexity can be further improved to an average $O(1)$, up to logarithmic factors, parallel complexity per right-hand side. All of these extensions are currently under investigation.

Finally, while we only consider a row-based processor assignment in this work, in highly heterogeneous parallel computational environments other processor assignments, in particular processor assignments tailored to the heterogeneous computational environment, can be applied straightforwardly with our approach. This may lead to further improvements of the efficiency of the algorithm. An investigation of these aspects would also be very interesting.

\section{Acknowledgement}
\label{Sec::Acknowledgement}
The authors thank Total SA for support and for permission to release the example code. LD is also supported by AFOSR grant FA9550-17-1-0316. 
The BP model is provided courtesy of BP and Fr\'ederic Billette. We thank NERSC for computation resources. MT thanks the Institute of Applied Mathematics at Graz University of Technology for hosting him for part of this research. 

\begin{appendices}
\section{Perfectly Matched Layers (PMLs) and Finite Difference Discretizations}\label{App::PML}
For a given squared slowness $m$ and source density $f$ defined on $\Omega_{\texttt{extended}}$, following the technique of PMLs, a new complex-valued squared slowness and source density can be derived so that the PMLs act as absorbing boundary conditions.
The coefficients in $\Omega_{\texttt{extended}}$ can be derived from the variable transformation~\cite{Johnson:PML,Berenger:PML}
\begin{alignat*}{1}
\frac{\partial}{\partial x_k}\rightarrow \alpha_k(x)\frac{\partial}{\partial x_k}
\end{alignat*}
for $k=1,..,d$ where $\alpha_k$ is the complex-valued function
\begin{alignat*}{1}
\alpha_k(x)=\frac{1}{1+i\frac{\sigma_k(x)}{\omega}},
\end{alignat*}
where $i=\sqrt{-1}$ is the imaginary unit and $\sigma_k$ is the PML profile function. 

To describe the PML profile, assume that $\Omega_{\texttt{bulk}}$ is the unit square or the unit cube. Then the extended domain $\Omega_{\texttt{extended}}$ is the square ($d=2$) or cube ($d=3$), $[-\delta_{\texttt{pml}},1+\delta_{\texttt{pml}}]^d$, with the PML-width $\delta_{\texttt{pml}}>0$. The PML profile $\sigma_k$ is then chosen to be the cubic function
\begin{alignat*}{1}
\sigma_k(x):=\left\{\begin{array}{ll}
\frac{C}{\delta_{\texttt{pml}}}\left(\frac{-x_k}{\delta_{\texttt{pml}}}\right)^3 \quad & \text{for }x_k<0\\
0 \quad & \text{for }0\leq x_k \leq 1 \\
\frac{C}{\delta_{\texttt{pml}}}\left(\frac{x_k-1}{\delta_{\texttt{pml}}}\right)^3 \quad & \text{for }x_k>1\\
\end{array}\right. 
\end{alignat*}
where $C>0$ is the absorption constant chosen to be $C=\ln \omega$.
Employing this variable transformation in equation~\eqref{HelmholtzExtended} gives rise to the diagonal matrix $\Lambda$ with diagonal entries $\alpha_k^2/\beta$, and the the complex-valued squared slowness and source density can be written as
\begin{alignat*}{1}
\frac{m}{\prod_{k=1}^d\alpha_k},\quad \text{and}\quad \frac{f}{\prod_{k=1}^d\alpha_k},
\end{alignat*}
respectively.

The resulting system is usually called the symmetric formulation of the Helmholtz equation, which in 2D takes the form: 
\begin{equation} \label{eq:Hemholtz_Symmetric}
-\left( \nabla \cdot \Lambda \nabla   + \frac{\omega ^2 m(\x)}{\alpha_1(\x) \alpha_2(\x)} \right) u(\x) = \frac{f(\x)}{\alpha_1(\x) \alpha_2(\x)},
\end{equation}
and
\begin{equation}
    \Lambda(\x) = \left [   \begin{array}{cc}
                            s_1(\x)  & 0 \\
                            0    & s_2(\x)
                        \end{array}
              \right],
 \end{equation}
 $s_1 = \alpha_1/\alpha_2$, $s_2 = \alpha_2/\alpha_1$,
with homogeneous Dirichlet boundary conditions.

We discretize $\Omega_{\texttt{bulk}}$ as an equispaced regular grid of stepsize $h$, and of dimensions $n_x \times n_y$. For the extended domain $\Omega_{\texttt{extended}}$, we extend this grid by $ n_{\texttt{pml}}=\delta_{\texttt{pml}}/h$ points in each direction, obtaining a grid of size $ (2n_{\texttt{pml}}+n_x) \times (2n_{\texttt{pml}}+n_z$). Define $\x_{p,q} = (x_p, y_q) = (ph,qh)$.
We use the 5-point stencil Laplacian to discretize~\eqref{eq:Hemholtz_Symmetric}. For the interior points $\x_{i,j} \in \Omega$, we have
\begin{equation} \label{eq:discrete_helmholtz_equation_by_point}
    \left( \boldA \u \right)_{p,q} = -\frac{1}{h^2} \left ( \u_{p-1,q}-2\u_{p,q} +\u_{p+1,q} \right) -  \frac{1}{h^2} \left ( \u_{p,q-1}-2\u_{p,q} +\u_{p,q+1} \right) - \omega^2 \frac{m(\x_{p,q})}{ s_1(\x_{p,q})s_2(\x_{p,q})}.
\end{equation}
In the PML, we discretize $\partial_x (s_1 \partial_x u ) $ as
\begin{equation}
\frac{ s_1(\x_{p+1/2,q})(\u_{p+1,q} - \u_{p,q}) - s_1(\x_{p-1/2,q})( \u_{p,q} - \u_{p-1,q} )   }{h^2},
\end{equation}
and analogously for $\partial_y(s_2 \partial_y u )$.

\section{Proof of the annihilation condition}\label{App::AnCond}
Let $x\in\Omega_1$, then using the definition of $U$ and the fact that $u$ and $G(x,.)$ vanish on $\partial\Omega$, it holds
\begin{alignat*}{1}
U(x)=&\int_{\Gamma}\Lambda(y)\nabla u(y)G(x,y)ds_y-\int_{\Gamma}u(y)\left[n(y)\cdot\left(\Lambda(y)\nabla_yG(x,y)\right)\right]ds_y\\
=&\int_{\partial\Omega_2}\Lambda(y)\nabla u(y)G(x,y)ds_y-\int_{\partial\Omega_2}u(y)\left[n(y)\cdot\left(\Lambda(y)\nabla_yG(x,y)\right)\right]ds_y.
\end{alignat*}
Then, employing the Divergence Theorem yields
\begin{alignat*}{1}
U(x)=&\int_{\Omega_2}u(y)\left[-\operatorname{div}_y\left(\Lambda(y)\nabla_yG(x,y)\right)\right]dy-\int_{\Omega_2}\left[-\operatorname{div}\left(\Lambda(y)\nabla u(y)\right)\right]G(x,y)dy\\
=&\int_{\Omega_2}u(y)\underbrace{\left[-\operatorname{div}_y\left(\Lambda(y)\nabla_yG(x,y)\right)-m(y)G(x,y)\right]}_{=0}dy-\int_{\Omega_2}\underbrace{\left[-\operatorname{div}\left(\Lambda(y)\nabla u(y)\right)-m(y)u(y)\right]}_{=0}G(x,y)dy
\end{alignat*}
proving the Annihilation condition~\eqref{AnnihilationCondition}.

\section{Pseudocode}
\label{App::Code}
\subsection{Core algorithms}
\begin{algorithm}[H]
\caption{Extraction of traces}
\begin{algorithmic}[1]
\Function{extractTraces}{$u$, $i$, $j$}
\State Define the Dirichlet traces $\lambda^{B}$, $\lambda^{R}$, $\lambda^{T}$, $\lambda^{L}$ and set them to zero
\State Define the Neumann traces $\mu^{B}$, $\mu^{R}$, $\mu^{T}$, $\mu^{L}$ and set them to zero

\If{$i>1$}
\State Extract the traces $\lambda^{B}$ and $\mu^{B}$ from $u$ according to~\eqref{Traces}
\EndIf
\If{$i<m$}
\State Extract the traces $\lambda^{T}$ and $\mu^{T}$ from $u$ according to~\eqref{Traces}
\EndIf
\If{$j>1$}
\State Extract the traces $\lambda^{L}$ and $\mu^{L}$ from $u$ according to~\eqref{Traces}
\EndIf
\If{$j<n$}
\State Extract the traces $\lambda^{R}$ and $\mu^{R}$ from $u$ according to~\eqref{Traces}
\EndIf

\Return $\lambda^B$, $\lambda^R$, $\lambda^T$, $\lambda^L$, $\mu^B$, $\mu^R$, $\mu^T$, $\mu^L$
\EndFunction
\end{algorithmic}
\end{algorithm}

\begin{algorithm}[H]
\caption{Computation of a local solution}
\begin{algorithmic}[1]
\Function{computeLocSol}{$f$, $i$, $j$}
\State Define $f_{ij}$ for any $x\in\Omega_{ij}^\varepsilon$ such that
\begin{alignat*}{1}
f_{ij}(x):=\left\{\begin{array}{ll}
f(x) &\quad x\in\Omega_{ij}\\
0 &\quad \text{otherwise}
\end{array}\right.
\end{alignat*}
\State Solve the local problem on $\Omega_{ij}^\varepsilon$ for $f_{ij}$ and obtain the solution $u_{ij}$

\Return $u_{ij}$
\EndFunction
\end{algorithmic}
\end{algorithm}

\begin{algorithm}[H]
\caption{Computation of a local polarized wavefield}
\label{Alg::PolSol}
\begin{algorithmic}[1]
\Function{computePolSol}{$\lambda$, $\mu$, $\Gamma$, $i$, $j$}
\State With $\lambda$, $\mu$ and $\Gamma$, compute $u_{ij}$ in $\Omega_{ij}^\varepsilon$ using~\eqref{RepForm}

\Return $u_{ij}$
\EndFunction
\end{algorithmic}
\end{algorithm}

\subsection{Algorithm for the local solutions (stage 1)}
\begin{algorithm}[H]
\caption{Compute all local solutions}
\label{Alg::CombHVTraces}
\begin{algorithmic}[1]
\Function{addLocSols}{$\boldsymbol{u}$, $f$}
\State Define an $m\times n$ matrix $\boldsymbol{\lambda}^L$ of left Dirichlet traces and set them to zero
\State Define an $m\times n$ matrix $\boldsymbol{\mu}^L$ of left Neumann traces and set them to zero
\State Define an $m\times n$ matrix $\boldsymbol{\lambda}^R$ of right Dirichlet traces and set them to zero
\State Define an $m\times n$ matrix $\boldsymbol{\mu}^R$ of right Neumann traces and set them to zero
\State Define an $m\times n$ matrix $\boldsymbol{\lambda}^B$ of bottom Dirichlet traces and set them to zero
\State Define an $m\times n$ matrix $\boldsymbol{\mu}^B$ of bottom Neumann traces and set them to zero
\State Define an $m\times n$ matrix $\boldsymbol{\lambda}^T$ of top Dirichlet traces and set them to zero
\State Define an $m\times n$ matrix $\boldsymbol{\mu}^T$ of top Neumann traces and set them to zero
\For{$i=1,\dots,m$}
\For{$j=1,\dots,n$}
\State $\boldsymbol{u}[i,j]$ = \Call{computeLocSol}{$f$, $i$, $j$}
\State $(\lambda^{B,\texttt{loc}},\lambda^{R,\texttt{loc}},\lambda^{T,\texttt{loc}},\lambda^{L,\texttt{loc}}, \mu^{B,\texttt{loc}},\mu^{R,\texttt{loc}},\mu^{T,\texttt{loc}},\mu^{L,\texttt{loc}})$ = \Call{extractTraces}{$\boldsymbol{u}[i,j]$, $i$, $j$}
\State Set the bottom Dirichlet trace $\boldsymbol{\lambda}^B[i,j]$ = $\lambda^{B,\texttt{loc}}$
\State Set the right Dirichlet trace $\boldsymbol{\lambda}^R[i,j]$ = $\lambda^{R,\texttt{loc}}$
\State Set the top Dirichlet trace $\boldsymbol{\lambda}^T[i,j]$ = $\lambda^{T,\texttt{loc}}$
\State Set the left Dirichlet trace $\boldsymbol{\lambda}^L[i,j]$ = $\lambda^{L,\texttt{loc}}$
\State Set the bottom Neumann trace $\boldsymbol{\mu}^B[i,j]$ = $\mu^{B,\texttt{loc}}$
\State Set the right Neumann trace $\boldsymbol{\mu}^R[i,j]$ = $\mu^{R,\texttt{loc}}$
\State Set the top Neumann trace $\boldsymbol{\mu}^T[i,j]$ = $\mu^{T,\texttt{loc}}$
\State Set the left Neumann trace $\boldsymbol{\mu}^L[i,j]$ = $\mu^{L,\texttt{loc}}$
\EndFor
\EndFor

\Return $\boldsymbol{u}$, $\boldsymbol{\lambda}^B$, $\boldsymbol{\lambda}^R$, $\boldsymbol{\lambda}^T$, $\boldsymbol{\lambda}^L$, $\boldsymbol{\mu}^B$, $\boldsymbol{\mu}^R$, $\boldsymbol{\mu}^T$, $\boldsymbol{\mu}^L$
\EndFunction
\end{algorithmic}
\end{algorithm}

\subsection{Algorithms for the horizontal and vertical sweeps (stage 2)}
\begin{algorithm}[H]
\caption{Global wavefield from the up sweeps}
\label{Alg::UpSweep}
\begin{algorithmic}[1]
\Function{addUpSweeps}{$\boldsymbol{u}$, $\boldsymbol{\lambda}^T$, $\boldsymbol{\mu}^T$}
\State Define an $m\times n$ matrix $\boldsymbol{\lambda}^L$ of left Dirichlet traces and set them to zero
\State Define an $m\times n$ matrix $\boldsymbol{\mu}^L$ of left Neumann traces and set them to zero
\State Define an $m\times n$ matrix $\boldsymbol{\lambda}^R$ of right Dirichlet traces and set them to zero
\State Define an $m\times n$ matrix $\boldsymbol{\mu}^R$ of right Neumann traces and set them to zero

\For{$d=1\dotsm+n$}
\For{$i=1\dotsm$}
\State $j=d-i+1$
\If{$i>1$}
\State $u^{\texttt{loc}}$ = \Call{computePolSol}{$\boldsymbol{\lambda}^T[i-1,j]$, $\boldsymbol{\mu}^T[i-1,j]$, $\Gamma_{ij}^B$, $i$, $j$}
\State $(\lambda^{B,\texttt{loc}},\lambda^{R,\texttt{loc}},\lambda^{T,\texttt{loc}},\lambda^{L,\texttt{loc}}, \mu^{B,\texttt{loc}},\mu^{R,\texttt{loc}},\mu^{T,\texttt{loc}},\mu^{L,\texttt{loc}})$ = \Call{extractTraces}{$u^{\texttt{loc}}$, $i$, $j$}
\State Update wave field $\boldsymbol{u}[i,j]\pluseq u^{\texttt{loc}}$
\State Update top Dirichlet trace $\boldsymbol{\lambda}^T[i,j]\pluseq\lambda^{T,\texttt{loc}}$
\State Update top Neumann trace $\boldsymbol{\mu}^T[i,j]\pluseq\mu^{T,\texttt{loc}}$
\State Set right Dirichlet trace $\boldsymbol{\lambda}^R[i,j]=\lambda^{R,\texttt{loc}}$
\State Set right Neumann trace $\boldsymbol{\mu}^R[i,j]=\mu^{R,\texttt{loc}}$
\State Set left Dirichlet trace $\boldsymbol{\lambda}^L[i,j]=\lambda^{L,\texttt{loc}}$
\State Set left Neumann trace $\boldsymbol{\mu}^L[i,j]=\mu^{L,\texttt{loc}}$
\EndIf
\EndFor
\EndFor

\Return $\boldsymbol{u}$, $\boldsymbol{\lambda}^{R}$, $\boldsymbol{\lambda}^{L}$, $\boldsymbol{\mu}^{R}$, $\boldsymbol{\mu}^{L}$
\EndFunction
\end{algorithmic}
\end{algorithm}

\begin{algorithm}[H]
\caption{Global wavefield from the down sweeps}
\label{Alg::DownSweep}
\begin{algorithmic}[1]
\Function{addDownSweeps}{$\boldsymbol{u}$, $\boldsymbol{\lambda}^B$, $\boldsymbol{\mu}^B$}
\State Define an $m\times n$ matrix $\boldsymbol{\lambda}^L$ of left Dirichlet traces and set them to zero
\State Define an $m\times n$ matrix $\boldsymbol{\mu}^L$ of left Neumann traces and set them to zero
\State Define an $m\times n$ matrix $\boldsymbol{\lambda}^R$ of right Dirichlet traces and set them to zero
\State Define an $m\times n$ matrix $\boldsymbol{\mu}^R$ of right Neumann traces and set them to zero

\For{$d=1,\dots,m+n$}
\For{$i=1,\dots,m$}
\State $j=d-i+1$
\State $i=m-i+1$
\State $j=n-j+1$
\If{$i<n$}
\State $u^{\texttt{loc}}$ = \Call{computePolSol}{$\boldsymbol{\lambda}^{B}[i+1,j]$, $\boldsymbol{\mu}^{B}[i+1,j]$, $\Gamma_{ij}^T$, $i$, $j$}
\State $(\lambda^{B,\texttt{loc}},\lambda^{R,\texttt{loc}},\lambda^{T,\texttt{loc}},\lambda^{L,\texttt{loc}}, \mu^{B,\texttt{loc}},\mu^{R,\texttt{loc}},\mu^{T,\texttt{loc}},\mu^{L,\texttt{loc}})$ = \Call{extractTraces}{$u^{\texttt{loc}}$, $i$, $j$}
\State Update wave field $\boldsymbol{u}[i,j]\pluseq u^{\texttt{loc}}$
\State Update bottom Dirichlet trace $\boldsymbol{\lambda}^B[i,j]\pluseq\lambda^{B,\texttt{loc}}$
\State Update bottom Neumann trace $\boldsymbol{\mu}^B[i,j]\pluseq\mu^{B,\texttt{loc}}$
\State Set right Dirichlet trace $\boldsymbol{\lambda}^R[i,j]=\lambda^{R,\texttt{loc}}$
\State Set right Neumann trace $\boldsymbol{\mu}^R[i,j]=\mu^{R,\texttt{loc}}$
\State Set left Dirichlet trace $\boldsymbol{\lambda}^L[i,j]=\lambda^{L,\texttt{loc}}$
\State Set left Neumann trace $\boldsymbol{\mu}^L[i,j]=\mu^{L,\texttt{loc}}$
\EndIf
\EndFor
\EndFor

\Return $\boldsymbol{u}$, $\boldsymbol{\lambda}^{R}$, $\boldsymbol{\lambda}^{L}$, $\boldsymbol{\mu}^{R}$, $\boldsymbol{\mu}^{L}$
\EndFunction
\end{algorithmic}
\end{algorithm}

\begin{algorithm}[H]
\caption{Global wavefield from the right sweeps}
\label{Alg::RightSweep}
\begin{algorithmic}[1]
\Function{addRightSweeps}{$\boldsymbol{u}$, $\boldsymbol{\lambda}^R$, $\boldsymbol{\mu}^R$}
\State Define an $m\times n$ matrix $\boldsymbol{\lambda}^B$ of bottom Dirichlet traces and set them to zero
\State Define an $m\times n$ matrix $\boldsymbol{\mu}^B$ of bottom Neumann traces and set them to zero
\State Define an $m\times n$ matrix $\boldsymbol{\lambda}^T$ of top Dirichlet traces and set them to zero
\State Define an $m\times n$ matrix $\boldsymbol{\mu}^T$ of top Neumann traces and set them to zero

\For{$i=1,\dots,m$}
\For{$j=2,\dots,n$}
\State $u^{\texttt{loc}}$ = \Call{computePolSol}{$\lambda^{R}[i,j-1]$, $\mu^{R}[i,j-1]$, $\Gamma_{ij}^L$, $i$, $j$}
\State $(\lambda^{B,\texttt{loc}},\lambda^{R,\texttt{loc}},\lambda^{T,\texttt{loc}},\lambda^{L,\texttt{loc}}, \mu^{B,\texttt{loc}},\mu^{R,\texttt{loc}},\mu^{T,\texttt{loc}},\mu^{L,\texttt{loc}})$ = \Call{extractTraces}{$u^{\texttt{loc}}$, $i$, $j$}
\State Update wave field $\boldsymbol{u}[i,j]\pluseq u^{\texttt{loc}}$
\State Update right Dirichlet trace $\boldsymbol{\lambda}^R[i,j]\pluseq\lambda^{R,\texttt{loc}}$
\State Update right Neumann trace $\boldsymbol{\mu}^R[i,j]\pluseq\mu^{R,\texttt{loc}}$
\State Set bottom Dirichlet trace $\boldsymbol{\lambda}^B[i,j]=\lambda^{B,\texttt{loc}}$
\State Set bottom Neumann trace $\boldsymbol{\mu}^B[i,j]=\mu^{B,\texttt{loc}}$
\State Set top Dirichlet trace $\boldsymbol{\lambda}^T[i,j]=\lambda^{T,\texttt{loc}}$
\State Set top Neumann trace $\boldsymbol{\mu}^T[i,j]=\mu^{T,\texttt{loc}}$
\EndFor
\EndFor

\Return $\boldsymbol{u}$, $\boldsymbol{\lambda}^{B}$, $\boldsymbol{\lambda}^{T}$, $\boldsymbol{\mu}^{B}$, $\boldsymbol{\mu}^{T}$
\EndFunction
\end{algorithmic}
\end{algorithm}

\begin{algorithm}[H]
\caption{Global wavefield from an left sweep in row $i$}
\label{Alg::LeftSweep}
\begin{algorithmic}[1]
\Function{addLeftSweeps}{$\boldsymbol{u}$, $\boldsymbol{\lambda}^L$, $\boldsymbol{\mu}^L$}
\State Define an $m\times n$ matrix $\boldsymbol{\lambda}^B$ of bottom Dirichlet traces and set them to zero
\State Define an $m\times n$ matrix $\boldsymbol{\mu}^B$ of bottom Neumann traces and set them to zero
\State Define an $m\times n$ matrix $\boldsymbol{\lambda}^T$ of top Dirichlet traces and set them to zero
\State Define an $m\times n$ matrix $\boldsymbol{\mu}^T$ of top Neumann traces and set them to zero

\For{$i=1,\dots,m$}
\For{$j=n-1,\dots,1$}
\State $u^{\texttt{loc}}$ = \Call{computePolSol}{$\lambda^{L}[i,j+1]$, $\mu^{L}[i,j+1]$, $\Gamma_{ij}^R$, $i$, $j$}
\State $(\lambda^{B,\texttt{loc}},\lambda^{R,\texttt{loc}},\lambda^{T,\texttt{loc}},\lambda^{L,\texttt{loc}}, \mu^{B,\texttt{loc}},\mu^{R,\texttt{loc}},\mu^{T,\texttt{loc}},\mu^{L,\texttt{loc}})$ = \Call{extractTraces}{$u^{\texttt{loc}}$, $i$, $j$}
\State Update wave field $\boldsymbol{u}[i,j]\pluseq u^{\texttt{loc}}$
\State Update left Dirichlet trace $\boldsymbol{\lambda}^L[i,j]\pluseq\lambda^{L,\texttt{loc}}$
\State Update left Neumann trace $\boldsymbol{\mu}^L[i,j]\pluseq\mu^{L,\texttt{loc}}$
\State Set bottom Dirichlet trace $\boldsymbol{\lambda}^B[i,j]=\lambda^{B,\texttt{loc}}$
\State Set bottom Neumann trace $\boldsymbol{\mu}^B[i,j]=\mu^{B,\texttt{loc}}$
\State Set top Dirichlet trace $\boldsymbol{\lambda}^T[i,j]=\lambda^{T,\texttt{loc}}$
\State Set top Neumann trace $\boldsymbol{\mu}^T[i,j]=\mu^{T,\texttt{loc}}$
\EndFor
\EndFor

\Return $\boldsymbol{u}$, $\boldsymbol{\lambda}^{B}$, $\boldsymbol{\lambda}^{T}$, $\boldsymbol{\mu}^{B}$, $\boldsymbol{\mu}^{T}$
\EndFunction
\end{algorithmic}
\end{algorithm}

\subsection{Algorithms for diagonal sweeps (stage 3)}
\begin{algorithm}[H]
\caption{Global wavefield from a diagonal sweep from the bottom-left to the top-right corner}
\label{Alg::DiagSweepBL2TR}
\begin{algorithmic}[1]
\Function{addBL2TRSweep}{$\boldsymbol{u}$, $\boldsymbol{\lambda}^{L}$, $\boldsymbol{\lambda}^{B}$, $\boldsymbol{\mu}^{L}$, $\boldsymbol{\mu}^{B}$}
\For{$d=1,\dots,m+n$}
\For{$i=1,\dots,m$}
\If{$i>1$ and $j>1$}
\State Define the L-shaped line $\Gamma^{BL}$by combining $\Gamma_{ij}^B$ and $\Gamma_{ij}^L$
\State Define $\lambda^{B}:=\boldsymbol{\lambda}^T[i-1,j]$ and $\mu^{B}:=\boldsymbol{\mu}^T[i-1,j]$
\State Define $\lambda^{L}:=\boldsymbol{\lambda}^R[i,j-1]$ and $\mu^{L}:=\boldsymbol{\mu}^R[i,j-1]$
\State Define $\lambda^{BL}$ and $\mu^{BL}$ similar to~\eqref{CombinedTracesDirichlet} and~\eqref{CombinedTracesNeumann}.
\State $u^{\texttt{loc}}$ = \Call{computePolSol}{$\lambda^{BL}$, $\mu^{BL}$, $\Gamma^{BL}$, $i$, $j$}
\State $(\lambda^{B,\texttt{loc}},\lambda^{R,\texttt{loc}},\lambda^{T,\texttt{loc}},\lambda^{L,\texttt{loc}}, \mu^{B,\texttt{loc}},\mu^{R,\texttt{loc}},\mu^{T,\texttt{loc}},\mu^{L,\texttt{loc}})$ = \Call{extractTraces}{$u^{\texttt{loc}}$, $i$, $j$}
\State Update wave field $\boldsymbol{u}[i,j]\pluseq u^{\texttt{loc}}$
\State Update right Dirichlet trace $\boldsymbol{\lambda}^R[i,j]\pluseq\lambda^{R,\texttt{loc}}$
\State Update right Neumann trace $\boldsymbol{\mu}^R[i,j]\pluseq\mu^{R,\texttt{loc}}$
\State Update top Dirichlet trace $\boldsymbol{\lambda}^T[i,j]\pluseq\lambda^{T,\texttt{loc}}$
\State Update top Neumann trace $\boldsymbol{\mu}^T[i,j]\pluseq\mu^{T,\texttt{loc}}$
\EndIf
\EndFor
\EndFor

\Return $\boldsymbol{u}$
\EndFunction
\end{algorithmic}
\end{algorithm}

\begin{algorithm}[H]
\caption{Global wavefield from a diagonal sweep from the top-right to the bottom-left corner}
\label{Alg::DiagSweepTR2BL}
\begin{algorithmic}[1]
\Function{addTR2BLSweep}{$\boldsymbol{u}$, $\boldsymbol{\lambda}^{R}$, $\boldsymbol{\lambda}^{T}$, $\boldsymbol{\mu}^{R}$, $\boldsymbol{\mu}^{T}$}
\For{$d=1,\dots,m+n$}
\For{$i=m,\dots,1$}
\State $j=n+m-1-(d+i)$
\If{$i<n$ and $j<n$}
\State Define the L-shaped line $\Gamma^{TR}$by combining $\Gamma_{ij}^T$ and $\Gamma_{ij}^R$
\State Define $\lambda^{T}:=\boldsymbol{\lambda}^B[i+1,j]$ and $\mu^{T}:=\boldsymbol{\mu}^B[i+1,j]$
\State Define $\lambda^{R}:=\boldsymbol{\lambda}^L[i,j+1]$ and $\mu^{R}:=\boldsymbol{\mu}^L[i,j+1]$
\State Define $\lambda^{TR}$ and $\mu^{TR}$ similar to~\eqref{CombinedTracesDirichlet} and~\eqref{CombinedTracesNeumann}.
\State $u^{\texttt{loc}}$ = \Call{computePolSol}{$\lambda^{TR}$, $\mu^{TR}$, $\Gamma^{TR}$, $i$, $j$}
\State $(\lambda^{B,\texttt{loc}}, \lambda^{R,\texttt{loc}}, \lambda^{T,\texttt{loc}}, \lambda^{L,\texttt{loc}},  \mu^{B,\texttt{loc}}, \mu^{R,\texttt{loc}}, \mu^{T,\texttt{loc}}, \mu^{L,\texttt{loc}})$ = \Call{extractTraces}{$u^{\texttt{loc}}$, $i$, $j$}
\State Update wave field $\boldsymbol{u}[i,j]\pluseq u^{\texttt{loc}}$
\State Update left Dirichlet trace $\boldsymbol{\lambda}^L[i,j]\pluseq\lambda^{L,\texttt{loc}}$
\State Update left Neumann trace $\boldsymbol{\mu}^L[i,j]\pluseq\mu^{L,\texttt{loc}}$
\State Update bottom Dirichlet trace $\boldsymbol{\lambda}^B[i,j]\pluseq\lambda^{B,\texttt{loc}}$
\State Update bottom Neumann trace $\boldsymbol{\mu}^B[i,j]\pluseq\mu^{B,\texttt{loc}}$
\EndIf
\EndFor
\EndFor

\Return $\boldsymbol{u}$
\EndFunction
\end{algorithmic}
\end{algorithm}

\begin{algorithm}[H]
\caption{Global wavefield from a diagonal sweep from the bottom-right to the top-left corner}
\label{Alg::DiagSweepBR2TL}
\begin{algorithmic}[1]
\Function{addBR2TLSweep}{$\boldsymbol{u}$, $\boldsymbol{\lambda}^{L}$, $\boldsymbol{\lambda}^{T}$, $\boldsymbol{\mu}^{L}$, $\boldsymbol{\mu}^{T}$}
\For{$d=1,\dots,m+n$}
\For{$i=1,\dots,m$}
\State $j=n-(d+i)$
\If{$i>1$ and $j<n$}
\State Define the L-shaped line $\Gamma^{BR}$by combining $\Gamma_{ij}^B$ and $\Gamma_{ij}^R$
\State Define $\lambda^{B}:=\boldsymbol{\lambda}^T[i-1,j]$ and $\mu^{B}:=\boldsymbol{\mu}^T[i-1,j]$
\State Define $\lambda^{R}:=\boldsymbol{\lambda}^L[i,j+1]$ and $\mu^{R}:=\boldsymbol{\mu}^L[i,j+1]$
\State Define $\lambda^{BR}$ and $\mu^{BR}$ similar to~\eqref{CombinedTracesDirichlet} and~\eqref{CombinedTracesNeumann}.
\State $u^{\texttt{loc}}$ = \Call{computePolSol}{$\lambda^{BR}$, $\mu^{BR}$, $\Gamma^{BR}$, $i$, $j$}
\State $(\lambda^{B,\texttt{loc}},\lambda^{R,\texttt{loc}},\lambda^{T,\texttt{loc}},\lambda^{L,\texttt{loc}}, \mu^{B,\texttt{loc}},\mu^{R,\texttt{loc}},\mu^{T,\texttt{loc}},\mu^{L,\texttt{loc}})$ = \Call{extractTraces}{$u^{\texttt{loc}}$, $i$, $j$}
\State Update wave field $\boldsymbol{u}[i,j]\pluseq u^{\texttt{loc}}$
\State Update left Dirichlet trace $\boldsymbol{\lambda}^L[i,j]\pluseq\lambda^{L,\texttt{loc}}$
\State Update left Neumann trace $\boldsymbol{\mu}^L[i,j]\pluseq\mu^{L,\texttt{loc}}$
\State Update top Dirichlet trace $\boldsymbol{\lambda}^T[i,j]\pluseq\lambda^{T,\texttt{loc}}$
\State Update top Neumann trace $\boldsymbol{\mu}^T[i,j]\pluseq\mu^{T,\texttt{loc}}$
\EndIf
\EndFor
\EndFor

\Return $\boldsymbol{u}$
\EndFunction
\end{algorithmic}
\end{algorithm}

\begin{algorithm}[H]
\caption{Global wavefield from a diagonal sweep from the top-left to the bottom-right corner}
\label{Alg::DiagSweepTL2BR}
\begin{algorithmic}[1]
\Function{addTL2BRSweep}{$\boldsymbol{u}$, $\boldsymbol{\lambda}^{R}$, $\boldsymbol{\lambda}^{B}$, $\boldsymbol{\mu}^{R}$, $\boldsymbol{\mu}^{B}$}
\For{$d=1,\dots,m+n$}
\For{$i=m,\dots,1$}
\State $j=(d+i)-m$
\If{$i<m$ and $j>1$}
\State Define the L-shaped line $\Gamma^{TL}$by combining $\Gamma_{ij}^T$ and $\Gamma_{ij}^L$
\State Define $\lambda^{T}:=\boldsymbol{\lambda}^B[i+1,j]$ and $\mu^{T}:=\boldsymbol{\mu}^B[i+1,j]$
\State Define $\lambda^{L}:=\boldsymbol{\lambda}^R[i,j-1]$ and $\mu^{L}:=\boldsymbol{\mu}^R[i,j-1]$
\State Define $\lambda^{TL}$ and $\mu^{TL}$ similar to~\eqref{CombinedTracesDirichlet} and~\eqref{CombinedTracesNeumann}.
\State $u^{\texttt{loc}}$ = \Call{computePolSol}{$\lambda^{TL}$, $\mu^{TL}$, $\Gamma^{TL}$, $i$, $j$}
\State $(\lambda^{B,\texttt{loc}},\lambda^{R,\texttt{loc}},\lambda^{T,\texttt{loc}},\lambda^{L,\texttt{loc}}, \mu^{B,\texttt{loc}},\mu^{R,\texttt{loc}},\mu^{T,\texttt{loc}},\mu^{L,\texttt{loc}})$ = \Call{extractTraces}{$u^{\texttt{loc}}$, $i$, $j$}
\State Update wave field $\boldsymbol{u}[i,j]\pluseq u^{\texttt{loc}}$
\State Update right Dirichlet trace $\boldsymbol{\lambda}^R[i,j]\pluseq\lambda^{R,\texttt{loc}}$
\State Update right Neumann trace $\boldsymbol{\mu}^R[i,j]\pluseq\mu^{R,\texttt{loc}}$
\State Update bottom Dirichlet trace $\boldsymbol{\lambda}^B[i,j]\pluseq\lambda^{B,\texttt{loc}}$
\State Update bottom Neumann trace $\boldsymbol{\mu}^B[i,j]\pluseq\mu^{B,\texttt{loc}}$
\EndIf
\EndFor
\EndFor

\Return $\boldsymbol{u}$
\EndFunction
\end{algorithmic}
\end{algorithm}

\subsection{The algorithm of scenario 3}
\begin{algorithm}[H]
\caption{Computation of the wavefield for an arbitrary source density that does not intersect the skeleton of the CDD}
\label{Alg::Scenario3}
\begin{algorithmic}[1]
\Function{computeScenario3}{$f$}
\State Define an $m\times n$ matrix $\boldsymbol{u}$ of local wave fields

\State ($\boldsymbol{u}$, $\boldsymbol{\lambda}^{B,\texttt{loc}}$, $\boldsymbol{\lambda}^{R,\texttt{loc}}$, $\boldsymbol{\lambda}^{T,\texttt{loc}}$, $\boldsymbol{\lambda}^{L,\texttt{loc}}$, $\boldsymbol{\mu}^{B,\texttt{loc}}$, $\boldsymbol{\mu}^{R,\texttt{loc}}$, $\boldsymbol{\mu}^{T,\texttt{loc}}$, $\boldsymbol{\mu}^{L,\texttt{loc}}$) = \Call{addLocSols}{$\boldsymbol{u}$, $f$}
\State ($\boldsymbol{u}$, $\boldsymbol{\lambda}^{R,\texttt{up}}$, $\boldsymbol{\lambda}^{L,\texttt{up}}$, $\boldsymbol{\mu}^{R,\texttt{up}}$, $\boldsymbol{\mu}^{L,\texttt{up}}$) = \Call{addUpSweeps}{$\boldsymbol{u}$, $\boldsymbol{\lambda}^{T,\texttt{loc}}$, $\boldsymbol{\mu}^{T,\texttt{loc}}$}
\State ($\boldsymbol{u}$, $\boldsymbol{\lambda}^{R,\texttt{down}}$, $\boldsymbol{\lambda}^{L,\texttt{down}}$, $\boldsymbol{\mu}^{R,\texttt{down}}$, $\boldsymbol{\mu}^{L,\texttt{down}}$) = \Call{addDownSweeps}{$\boldsymbol{u}$, $\boldsymbol{\lambda}^{B,\texttt{loc}}$, $\boldsymbol{\mu}^{B,\texttt{loc}}$}
\State ($\boldsymbol{u}$, $\boldsymbol{\lambda}^{B,\texttt{left}}$, $\boldsymbol{\lambda}^{T,\texttt{left}}$, $\boldsymbol{\mu}^{B,\texttt{left}}$, $\boldsymbol{\mu}^{T,\texttt{left}}$) = \Call{addLeftSweeps}{$\boldsymbol{u}$, $\boldsymbol{\lambda}^{L,\texttt{loc}}$, $\boldsymbol{\mu}^{L,\texttt{loc}}$}
\State ($\boldsymbol{u}$, $\boldsymbol{\lambda}^{B,\texttt{right}}$, $\boldsymbol{\lambda}^{T,\texttt{right}}$, $\boldsymbol{\mu}^{B,\texttt{right}}$, $\boldsymbol{\mu}^{T,\texttt{right}}$) = \Call{addRightSweeps}{$\boldsymbol{u}$, $\boldsymbol{\lambda}^{R,\texttt{loc}}$, $\boldsymbol{\mu}^{R,\texttt{loc}}$}

\State $\boldsymbol{u}$ = \Call{addBL2TRSweep}{$\boldsymbol{u}$, $\boldsymbol{\lambda}^{R,\texttt{up}}$, $\boldsymbol{\lambda}^{T,\texttt{right}}$, $\boldsymbol{\mu}^{R,\texttt{up}}$, $\boldsymbol{\mu}^{T,\texttt{right}}$}
\State $\boldsymbol{u}$ = \Call{addTR2BLSweep}{$\boldsymbol{u}$, $\boldsymbol{\lambda}^{L,\texttt{down}}$, $\boldsymbol{\lambda}^{B,\texttt{left}}$, $\boldsymbol{\mu}^{L,\texttt{down}}$, $\boldsymbol{\mu}^{B,\texttt{left}}$}
\State $\boldsymbol{u}$ = \Call{addBR2TLSweep}{$\boldsymbol{u}$, $\boldsymbol{\lambda}^{L,\texttt{up}}$, $\boldsymbol{\lambda}^{T,\texttt{left}}$, $\boldsymbol{\mu}^{L,\texttt{up}}$, $\boldsymbol{\mu}^{T,\texttt{left}}$}
\State $\boldsymbol{u}$ = \Call{addTL2BRSweep}{$\boldsymbol{u}$, $\boldsymbol{\lambda}^{R,\texttt{down}}$, $\boldsymbol{\lambda}^{B,\texttt{right}}$, $\boldsymbol{\mu}^{R,\texttt{down}}$, $\boldsymbol{\mu}^{B,\texttt{right}}$}

\State Define $u=0$ in $\Omega$
\For{$i=1,\dots,m$}
\For{$j=1,\dots,n$}
\State $u|_{\Omega_{ij}}$=$\boldsymbol{u}[i,j]|_{\Omega_{ij}}$
\EndFor
\EndFor

\Return $u$
\EndFunction
\end{algorithmic}
\end{algorithm}
\end{appendices}

\bibliographystyle{plain}
\bibliography{LSweeps}

\end{document}